\documentclass[11pt]{article}
\usepackage[top=1.20in,bottom=1.20in,left=1.00in,right=1.00in]{geometry}
\usepackage{authblk}
\usepackage{amsmath,amsthm,amssymb,url,bm}
\usepackage{graphicx}
\usepackage{caption}
\usepackage{subcaption}
\usepackage{spverbatim}
\usepackage{listings,wasysym}
\usepackage{listings}
\usepackage{enumerate}
\usepackage{algorithm}
\usepackage{xcolor}
\usepackage[numbers, square]{natbib} 
\RequirePackage[colorlinks,citecolor=blue,urlcolor=blue]{hyperref}

\usepackage{accents}
\usepackage{import}
\usepackage{mathtools}
\usepackage[normalem]{ulem}
\usepackage[noend]{algpseudocode}
\graphicspath{{../figs/}}

\theoremstyle{plain}
\newtheorem{thm}{Theorem}[section]
\newtheorem{lem}[thm]{Lemma}

\theoremstyle{definition}

\newtheorem{theorem}{Theorem}[section]
\newtheorem{definition}[theorem]{Definition}

\theoremstyle{remark}
\newtheorem{rem}{Remark}

\newcommand{\one}{\textbf{1}}

\newcommand{\es}{\emptyset}
\newcommand{\R}{\mathbb{R}}
\newcommand{\E}{\mathbb{E}}
\renewcommand{\P}{\mathbb{P}}

\newcommand{\G}{\mathcal{G}}
\newcommand{\T}{\mathcal{T}}

\newcommand{\lb}{\left(}
\newcommand{\rb}{\right)}
\newcommand{\td}{\tilde}
\newcommand{\Ep}{\td{E}}
\newcommand{\bp}{\td{b}}
\newcommand{\sigmap}{\td{\sigma}}
\newcommand{\etap}{\td{\eta}}
\newcommand{\vct}[1]{\bm{#1}}
\newcommand{\mtx}[1]{\bm{#1}}
\renewcommand{\vec}[1]{{\boldsymbol{#1}}}
\newcommand{\A}{\mtx{A}}
\newcommand{\I}{\mtx{I}}
\newcommand{\J}{\mtx{J}}
\renewcommand{\L}{\mtx{L}}

\newcommand{\diag}{\mathrm{diag}}

\newcommand{\colspan}{\operatorname{colspan}}
\newcommand{\op}{\mathrm{op}}

\newcommand{\sign}{\mathrm{sign}}

\newcommand{\ubar}[1]{\underaccent{\bar}{#1}}
\newcommand{\nm}[1]{\left\lVert#1\right\rVert}
\newcommand{\mnorm}[1]{\|#1\|_{2\rightarrow\infty}}
\newcommand{\maxnorm}[1]{\|#1\|_{\max}}
\newcommand{\tti}{2\rightarrow\infty}
\DeclareMathOperator{\argmax}{argmax}

\usepackage{stackengine}
\stackMath
\newcommand\tsup[2][2]{%
 \def\useanchorwidth{T}%
  \ifnum#1>1%
    \stackon[-.5pt]{\tsup[\numexpr#1-1\relax]{#2}}{\scriptscriptstyle\sim}%
  \else%
    \stackon[.5pt]{#2}{\scriptscriptstyle\sim}%
  \fi%
}

\begin{document}

\def\spacingset#1{\renewcommand{\baselinestretch}%
{#1}\small\normalsize} \spacingset{1}
\spacingset{1.5}

\title{Consistency of Spectral Clustering on Hierarchical Stochastic Block Models}
\author[1]{\rm Lihua Lei}
\author[2]{\rm Xiaodong Li}
\author[2]{\rm Xingmei Lou}
\affil[1]{Department of Statistics, Stanford University}
\affil[2]{Department of Statistics, University of California, Davis}

%\author{Lihua Lei \thanks{Department of Statistics, Stanford University, Stanford, CA 94305}, ~~Xiaodong Li\thanks{Department of Statistics, University of California Davis, Davis CA 95616}, ~~and Xingmei Lou\thanks{Department of Statistics, University of California Davis, Davis CA 95616}}
\date{}
\maketitle

\begin{abstract}
We study the hierarchy of communities in real-world networks under a generic stochastic block model, in which the connection probabilities are structured in a binary tree. Under such model, a standard recursive bi-partitioning algorithm is dividing the network into two communities based on the Fiedler vector of the unnormalized graph Laplacian and repeating the split until a stopping rule indicates no further community structures. We prove the strong consistency of this method under a wide range of model parameters, which include sparse networks with node degrees as small as $O(\log n)$. In addition, unlike most of existing work, our theory covers multiscale networks where the connection probabilities may differ by orders of magnitude, which comprise an important class of models that are practically relevant but technically challenging to deal with. Finally we demonstrate the performance of our algorithm on synthetic data and real-world examples.
\end{abstract}


\section{Introduction}\label{sec:intro}
Community structures of real-world networks are typically hierarchical. In a coauthor network, it is not clear-cut whether all statisticians should be viewed as a single community --- at a high level of granularity, they can be combined with mathematicians, physicists, computer scientists, and so on as quantitative scientists, while at a low level of granularity, they can be further split into finer groups based on research areas. When the desired level of granularity is unknown a priori, a hierarchy is a more informative representation of the relational information than a single partition of the network. For example, to design a clusterwise randomized experiment, an A/B test designer can trade off the amount of interference (e.g., the number of edges between clusters) and effective sample size (e.g., the number of clusters) by searching over the hierarchy. 

Agglomerative community detection algorithms \citep[e.g.][]{girvan2002community} are intrinsically hierarchical because they are able to produce a dendrogram characterizing the hierarchy of communities. However, bottom-up algorithms are sensitive to noise when amalgamating small clusters at the beginning of the run. As a consequence, theoretical guarantees are hard to come by for sparse and noisy networks. In contrast, divisive community detection algorithms, such as spectral clustering, has been proved to recover the community structure theoretically under various sparse network models \citep[e.g.][]{mcsherry2001spectral,dasgupta2006spectral,rohe2011spectral,balakrishnan2011noise,lei2015consistency,jin2015fast,abbe2015exact,li2018hierarchical}. However, the network models for theoretical analysis rarely encode the hierarchy; often the communities are treated as logically separate units. Algorithmically, most divisive algorithms which have been analyzed in theory are unable to produce a dendrogram. It is somewhat disappointing that hierarchical algorithms typically have no theoretical guarantees (under practically reasonable assumptions) while those justified in theory are often non-hierarchical algorithmically. 

To mitigate this gap, one would need to consider a \emph{hierarchical clustering algorithm} and study its statistical properties under a \emph{hierarchical network model}. There have been attempts on this route, focusing on recursive divisive clustering algorithms, which recursively partition the network based on a top-down algorithm. A handful of such algorithms have been shown to recover the hierarchy under dense network models with average degree polynomial in $n$, where $n$ is the number of nodes \citep{balakrishnan2011noise, dasgupta2016cost, lyzinski2016community}. However, real-world networks are typically much sparser. \cite{dasgupta2006spectral} analyzed a recursive spectral algorithm under a network model with average degree $O(\log^6 n)$. Apart from the artificial exponent, the algorithm involves multiple tuning parameters with no recommended default values, making it hard to implement in practice. Recently, \cite{li2018hierarchical} proposed the Binary Tree Stochastic Block Model (BTSBM) which encodes a binary hierarchy among primitive communities in the spirit of \cite{clauset2008hierarchical} and \cite{balakrishnan2011noise}. They analyzed a recursive spectral clustering algorithm, which splits the network into two clusters based on the signs of the components in an eigenvector of the adjacency matrix, under the BTSBM and showed it consistently recovers the hierarchy when the average degree scales as $O(\log^{2+\epsilon} n)$ for $\epsilon > 0$. The refined analysis by \cite{lei2019unified} brought it down to the critical regime with an $O(\log n)$ average degree. Despite being able to handle sparse networks, their analyses are based on a restrictive model which assumes a balanced hierarchy and strict layerwise homogeneity in connection probabilities. 

In this paper, we analyze a Laplacian-based recursive bi-partitioning algorithm. It recursively splits the network into two based on the signs of the Fiedler vector \citep{fiedler1975property}, the eigenvector of the unnormalized Laplacian (formally defined in Section \ref{sec:model}) corresponding to the second smallest eigenvalue. The procedure is repeated iteratively until a stopping rule indicates that there are no further communities in any subgraphs. \cite{li2018hierarchical} suggested various stopping rules that work reasonably well empirically and provided theoretical justification for certain ones. As shown by \cite{li2018hierarchical}, to get $K$ communities, this algorithm is computationally more efficient than the $K$-way spectral clustering algorithm which splits the network into $K$ communities at once, especially for a large $K$.   

For the theoretical analysis, we consider a more general hierarchical SBM that has been proposed in \cite{clauset2008hierarchical} and \cite{balakrishnan2011noise} to allow for an unbalanced hierarchy and heterogeneous connection probabilities. We prove that the proposed algorithm consistently recovers the hierarchy for sparse networks in the critical regime where the average degree scales as $O(\log n)$. Notably, as opposed to \cite{li2018hierarchical} and most works on non-hierarchical SBMs, we do not require all connection probabilities to be on the same scale; instead, we allow the connection probability between communities closer on the hierarchy to be orders of magnitude larger than that between communities farther apart. It makes our analysis more realistic since real-world networks are often multiscale. Meanwhile, the theoretical analysis for clustering multiscale networks becomes much more challenging. To highlight our main theoretical contribution, we will not investigate the performance of stopping rules in this paper.

Our theory is built on (1) that the eigenvectors of the population Laplacian can identify the hierarchy; and (2) an entrywise perturbation bound showing that the Fiedler vector of the observed Laplacian approximates the population version with a high accuracy. The first part generalizes the result of \cite{balakrishnan2011noise} on the Fiedler vector to the entire eigenstructure. The second part rests on the recent development of $\ell_{\tti}$ eigenvector perturbation theory \citep{abbe2017entrywise,eldridge2017unperturbed,mao2017estimating,cape2019signal,lei2019unified,damle2020uniform}. The challenge is twofold: dependence between entries in the Laplacian and multiscale connection probabilities. We tackle the first with the technique developed by \cite{lei2019unified}, which, unlike most other perturbation bounds for random matrices, allows certain dependency structure among the entries. However, this technique alone is not enough to handle multiscale networks. We overcome this challenge by introducing novel techniques, and we will elaborate on them in Section \ref{sec:theory}.

\subsection{Notation}
We use $[n]$ to denote the set $\{1, \ldots, n\}$ and $\mtx{e}_{j}$ to denote the $j$-th canonical basis where the $j$-th element equals to $1$ and all other elements equal to $0$ (with the dimension depending on the context). Vectors and matrices are boldfaced while scalars are not. We denote by $\I_n$ the $n\times n$ identity matrix and by $\one_{n}$ the $n\times 1$ column vector with all entries $1$.
% and by $\one_{n\times m}$ the $n\times m$ matrix with all entries $1$.
For any vector $\vec{v}$, let $\|\vec{v}\|_{p}$ denote its $\ell_p$ norm. 
For any matrix $\mtx{M}$, let $\mtx{M}_{k}^{T}$ denote the $k$-th row of $\mtx{M}$, $\|\mtx{M}\|$ its spectral norm, and $\|\mtx{M}\|_{\mathrm{F}}$ its Frobenius norm.
Further, we denote by $\lambda_{1}(\mtx{M}), \ldots, \lambda_{n}(\mtx{M})$ the eigenvalues of $\mtx{M}$ in descending order, with $\vec{u}_1(\mtx{M}), \ldots, \vec{u}_{n}(\mtx{M})$ being the corresponding eigenvectors.
% Further we denote by $\lambda_{\max}(\mtx{M})$ (resp. $\lambda_{\min}(\mtx{M})$) the largest (resp. the smallest) eigenvalue of $\mtx{M}$ in absolute values, by $\kappa(\mtx{M})$ the condition number $\lambda_{\max}(\mtx{M}) / \lambda_{\min}(\mtx{M})$. 
% For any vector $\vec{v}$, $\diag(\vec{v})$ denotes the diagonal matrix with $\vec{v}$ as the diagonal entries. For any square matrix $\mtx{M}$, $\diag(\mtx{M})$ denotes the matrix obtained by setting all off-diagonal entries of $\mtx{M}$ to $0$. 
For two sequences of real numbers $\{x_n\}$ and $\{y_n\}$, we write $x_n = o(y_n)$ or $y_n=\omega(x_n)$ if $\lim _{n \rightarrow \infty} \frac{x_n}{y_n}=0$, $x_n = O(y_n)$ or $x_n \lesssim y_n$ if $|x_n|<C|y_n|$ for some constant $C$. Likewise, $x_n = \Omega(y_n)$ or $x_n \gtrsim y_n$ represents that there exists a constant $C$ such that $|x_n|>C|y_n|$. Finally, we write $x_n \asymp y_n$ if $x_n \lesssim y_n$ and $y_n \lesssim x_n$.

\section{A Hierarchical Stochastic Block Model}\label{sec:model}
% section 2 model

\subsection{Model formulation}
\label{subsec:formulation}
The \emph{Stochastic Block Model} (SBM) proposed by \cite{holland1983stochastic} has been widely used to study the empirical performance and theoretical properties of community detection methods. An SBM can be characterized by a vector $c=\{c_1,\ldots,c_n\}\in \{1, \ldots, K\}^n$ encoding the community membership of each node and a symmetric matrix $\bm{B}\in [0, 1]^{K\times K}$ encoding the connection probabilities between communities. The upper triangular part of the adjacency matrix $\mtx{A}$ has independent entries with $\mtx{A}_{ij} \sim \text{Bernoulli}(\mtx{B}_{c_i c_j})$ for any $i \le j$. By definition, the expected adjacency matrix, or equivalently the matrix of connection probabilities, can be represented as
\[
\E[\mtx{A}] \coloneqq \mtx{P} = \mtx{Z}\mtx{B}\mtx{Z}^\top - \diag(\mtx{Z}\mtx{B}\mtx{Z}^\top)\in \R^{n\times n},
\]
where $\mtx{Z}\in \R^{n \times K}$ denotes the membership matrix with the $i$-th row vector $\mtx{Z}_i=\vec{e}_{c_i}^\top$. %\sout{We allow for self loops for simplicity. } \Lihua{(Lihua: since the unnormalized graph Laplacian we considered does not depend on the diagonal elements so we do not need to make this simplifying assumption.)}

We consider a general Binary Tree Stochastic Block Model (BTSBM), which has been essentially proposed in \cite{clauset2008hierarchical} and \cite{balakrishnan2011noise} with slightly different emphasis. Specifically, given a binary tree $\T$ with $K$ leaf nodes, we represent each non-root node by a binary string recording the moves along the (unique) path from the root to that node, with $0$ denoting a left move and $1$ denoting a right move. For node $s$, let $|s|$ denote the depth of $s$ and $(L(s), R(s))$ its two children nodes. For a pair of nodes $s_1$ and $s_2$, we denote by $\mathcal{A}(s_1, s_2)$ their lowest common ancestor. In such model, each node $s$ on the binary tree $\T$ encodes two pieces of information. The first is a subset of units $\mathcal{G}_s\subset \{1, \ldots, n\}$, with $n_s = |\mathcal{G}_s|$. The model assumes that
$\mathcal{G}_{\emptyset} = \{1, \ldots, n\}$ and $\{\mathcal{G}_{L(s)}, \mathcal{G}_{R(s)}\}$ forms a partition of $\mathcal{G}_{s}$, i.e., $\mathcal{G}_{L(s)}\cap \mathcal{G}_{R(s)}=\emptyset$ and $\mathcal{G}_{L(s)}\cup \mathcal{G}_{R(s)} = \mathcal{G}_{s}$. We refer to the network encoded by a leaf node as a \emph{primitive community} and the network encoded by an internal node as a \emph{mega-community}. The second piece of information is a connection probability $p_s\in [0, 1]$. For each pair of units $i\not = j$, the model assumes that the connection probability between this pair of nodes is 
\[\mtx{P}_{ij} = p_{\mathcal{A}(c(i), c(j))},\]
where $c(i)$ denotes the (unique) leaf node that contains $i$.

\begin{figure}[h]
  \centering
  \includegraphics[width = 0.45\textwidth]{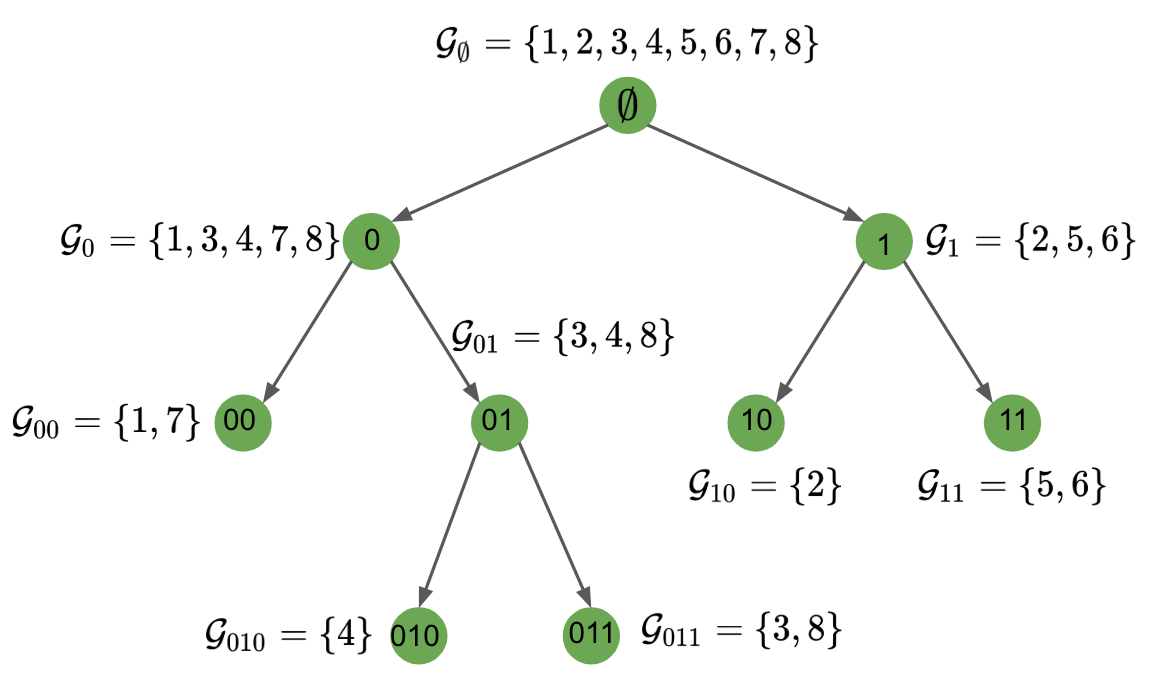}
  \includegraphics[width = 0.45\textwidth]{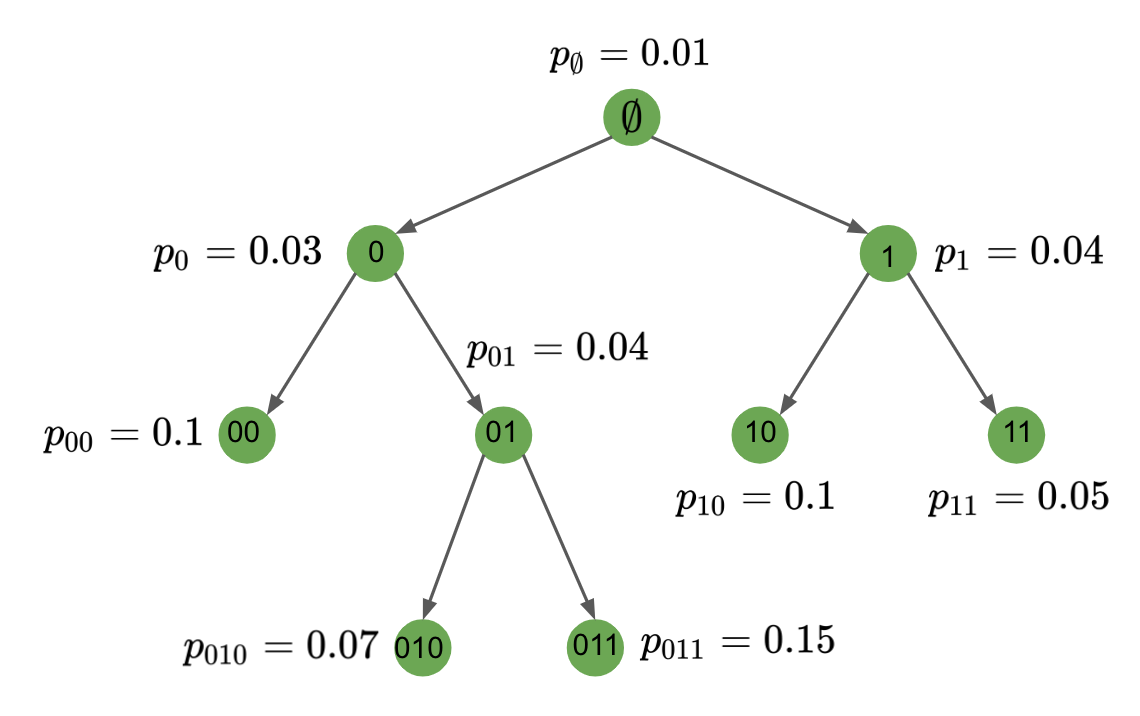}
  \caption{Illustration of a general BTSBM with $n = 8$ units and $K = 5$ primitive communities: (left) sub-network $\mathcal{G}_s$; (right) connection probability $p_s$.}\label{fig:illustrate_hbm}
\end{figure}

We illustrate in Figure \ref{fig:illustrate_hbm} the above definitions by a toy example with $n = 8$ units and a binary tree $\T$ with $K=5$ leaf nodes. The left panel shows the sub-networks encoded by each node. It can be read from the leaf nodes that $c(1) = c(7) = 00$, $c(2) = 10$, $c(3) = c(8) = 011$, $c(4) = 010$, and $c(5) = c(6) = 11$. The right panel shows the connection probabilities. Since $c(3) = 011, c(4) = 010$, and the lowest common ancestor of nodes $011$ and $010$ is $01$, $\mtx{P}_{34} = p_{01} = 0.04$. Similarly, $\mtx{P}_{38} = p_{011} = 0.15, \mtx{P}_{31} = p_{0} = 0.03$ and $\mtx{P}_{35} = p_{\emptyset} = 0.01$.

Clearly, the hierarchical SBM defined above is a special case of SBM with $K = 5$ communities  and the between-community connection probability matrix is
\[
\mtx{B}=\begin{bmatrix}
p_{00} & p_0 & p_0 & p_\es & p_\es\\
p_0 & p_{010} & p_{01}& p_\es & p_\es \\
p_0 & p_{01} & p_{011} & p_\es & p_\es \\
p_\es & p_\es & p_\es & p_{10} & p_1 \\
p_\es & p_\es & p_\es & p_{1} & p_{11}\\
\end{bmatrix}.
\]
% Meanwhile, it is also a special case of the Hierarchical Constant Block Matrix (HCBM) proposed by \cite{balakrishnan2011noise}, which considers a weighted graph. 
Note that the BTSBM discussed in \cite{li2018hierarchical} is a restrictive special case of the model we study here in that $\T$ is assumed to be full and balanced with  $p_{s} = p_{s'}$ and $n_{s} = n_{s'}$ whenever $|s| = |s'|$.

Throughout the paper we focus on general BTSBMs that satisfy the weak-assortativity.
\begin{definition}
\label{ass:weark_assort}
A general BTSBM $(\T, \{p_s: s \in \T\}, c(\cdot))$ is weakly assortative iff $p_{s} < p_{s'}$ whenever $s$ is the parent node of $s'$.
\end{definition}
This is a natural property that is compatible with the intuitive explanations of hierarchies --- nodes that are closer on the hierarchy are more likely to be connected. The example in Figure \ref{fig:illustrate_hbm} is weakly assortative. \cite{li2018hierarchical} also discussed the disassortative setting, though it is less common than the assortative setting in practice. 
 
We close this subsection with some notation that will be used repeatedly: 
\begin{equation}
\label{eq:pop_para}
p^{*} = \max_{1 \leq i, j \leq n}\mtx{P}_{ij}, \quad \bar{p}^{*} = \max_{1 \leq i \leq n}\frac{1}{n}\sum_{j=1}^{n}\mtx{P}_{ij}
, \quad \ubar{p}_{*} = \min_{1 \leq i \leq n}\frac{1}{n}\sum_{j=1}^{n}\mtx{P}_{ij}.   
\end{equation}
By definition, $p^{*}$ is the largest connection probability across the whole network, and $n\bar{p}^{*}$ ($n\ubar{p}_{*}$) is the largest (smallest) expected degree. Obviously, $p^{*}\ge \bar{p}^{*}\ge \ubar{p}_{*}$.

\subsection{Population unnormalized graph Laplacian}
The observed unnormalized graph Laplacian of a network is defined as 
\[\L = \mtx{D} - \A,\] 
where $\A\in \{0, 1\}^{n\times n}$ is the adjacency matrix and $\mtx{D} = \diag(\A \one_{n})$ is the diagonal matrix whose diagonal entries are the node degrees. It is known that $\L$ is positive semidefinite.  Let $\lambda_1 \geq \cdots \geq \lambda_{n -1} \geq \lambda_n = 0$ denote the eigenvalues of $\L$, and $\vec{u}_1,  \ldots \vec{u}_{n-1},\vec{u}_n$ the corresponding unit eigenvectors. Note that we always have $\vct{u}_n = \frac{1}{\sqrt{n}}\vct{1}_n$. In the case $\lambda_{n-1}=0$, we choose $\vct{u}_{n-1}$ to be any eigenvector corresponding to $\lambda_{n-1}=0$ which is orthogonal to $\vct{1}_n$.

As discussed in Section \ref{sec:intro}, our recursive spectral clustering algorithm splits the whole network into two based on the signs of the Fiedler eigenvector $\vec{u}_{n-1}$. Since $\L$ is an approximation of the population unnormalized graph Laplacian $\L^*=\E [\L]$, we shall study the eigenstructure of $\L^*$ as a stepping stone to prove the consistency of hierarchical clustering. Clearly, $\L^* = \diag(\mtx{P}\one) - \mtx{P}$. As with the sample version, the eigenvalues and unit eigenvectors of $\L^*$ are denoted as $\lambda^*_1 \geq \cdots \geq \lambda^*_{n -1} \geq \lambda^*_n = 0$ and $\vec{u}^*_1,  \ldots \vec{u}^*_{n-1},\vec{u}^*_n$. Theorem \ref{thm:eigen} provides an elegant characterization of the eigenstructure of $\L^*$ under weak assortativity. The proof is relegated to Section \ref{sec:proof_eigen}.

\begin{thm}
\label{thm:eigen}
Under a weakly assortative general BTSBM, defined in Section \ref{subsec:formulation},
\begin{enumerate}[(1)]
\item $\lambda_{n-1}^{*} = np_{\es}$ with multiplicity $1$ and the entries of the corresponding eigenvector obeys
\[
\mtx{u}^*_{n-1, i} =  \pm\left\{\begin{array}{ll}
     \sqrt{n_{1} / (n_{0} n)} & (i\in \G_{0})\\
      -\sqrt{n_{0} / (n_{1} n)} & (i\in \G_{1})
    \end{array}
\right.;
\]
\item $\lambda_{n-2}^{*} = \min\{n_{1}p_{1} + n_{0}p_{\es}, n_{0}p_{0} + n_{1}p_{\es}\}$;
\item The number of eigenvalues, accounting for the multiplicity, that are strictly less than $n\ubar{p}_{*}$ is at most $K$, the number of leaf nodes in $\T$;
\item For any $j$ with $\lambda_j^{*} < n\ubar{p}_{*}$, $\|\vec{u}_{j}^{*}\|_{\infty}\le \sqrt{\xi / n}$ where
\begin{equation}\label{eq:xi}
\xi:= \max_{s}\frac{n}{n_s}.
\end{equation}
Here recall that $n_s$ is the number of nodes in the subgragh encoded by the tree node $s$.
% If $\max\{\frac{n_0}{n_1},\frac{n_1}{n_0}\} = O(1)$ \XL{This may not be correct} and $K = O(1)$, then $\max\limits_{j:\lambda_j^{*} < n\ubar{p}_{*}}\|\vec{u}_{j}^{*}\|_{\infty} = O(\frac{1}{\sqrt{n}})$.
\end{enumerate}
\end{thm}
Theorem \ref{thm:eigen} (1) implies that the entry signs of $\vec{u}^*_{n-1}$ encode the first split of the network. Theorem \ref{thm:eigen} (1) and (2) have been proved in \cite{balakrishnan2011noise} and we include them here for completeness. In contrast, Theorem \ref{thm:eigen} (3) and (4) are new to the best of our knowledge. Although the eigenvalues and eigenvectors other than the Fiedler eigenpair appear to be algorithmically irrelevant, they are crucially useful in our theoretical analysis on the strong consistency of the clustering, especially for multiscale networks whose connection probabilities have different scales.

\section{Consistency of Recursive Spectral Clustering}\label{sec:theory}
% section 3 main results
\subsection{Criteria for consistency}\label{subsec:criteria}
Recovering the hierarchy encoded by the tree is equivalent to recovering all mega-communities and primitive communities. When primitive communities are hard to be recovered, owing to either the small size or insufficient gap with others, it is still possible to recover some mega-communities at the top levels of the tree, yielding a partial hierarchy that is informative for downstream analysis. In either case, it is necessary to investigate the consistency of the first split, that is, whether the Fiedler vector of the unnormalized graph Laplacian can partition $\mathcal{G}_{\es}$ accurately into $\mathcal{G}_{0}$ and $\mathcal{G}_1$.

We primarily focus on the strong consistency which requires the partition to be exactly correct with high probability. Without loss of generality, assume that $\vec{u}_{n-1, 1} \ge 0$ and $\vec{u}_{n-1, 1}^{*}\ge 0$. By Theorem \ref{thm:eigen}, the strong consistency of the first split can be formally stated as
\begin{equation}\label{eq:strong_consistency}
\P\lb\sign(\vec{u}_{n-1, i}) \neq \sign(\vec{u}_{n-1, i}^{*}), \,\, \text{for some }i\in \{1, \ldots, n\}\rb = o(1),
\end{equation}
and it can be easily extended to the consistent recovery of the whole hierarchy under the general BTSBM. Indeed, suppose we can identify a set of conditions under which the split at $\mathcal{G}_{\es}$ is strongly consistent, replacing $\mathcal{G}_{\es}$ with $\mathcal{G}_0$ yields the conditions for $\mathcal{G}_{00}$ and $\mathcal{G}_{01}$ to be exactly recovered, because the model for $\mathcal{G}_0$ is still a general BTSBM. Therefore, it is sufficient and necessary to investigate the strong consistency of the first split in order to establish the exact recovery of the full or partial hierarchy.

Another commonly studied criterion is the weak consistency which states that the misclustering error is asymptotically vanishing, i.e.,
\begin{equation}\label{eq:weak_consistency}
\min_{a\in \{-1, 1\}}\frac{1}{n}\sum_{i=1}^{n}I(\sign(\vec{u}_{n-1, i})\cdot \sign(\vec{u}_{n-1, i}^{*}) = a) = o_\P(1).
\end{equation}
In contrast to the strong consistency, the weak consistency of the first split does not carry over to lower splits under a general BTSBM, because the recovered $\mathcal{G}_0$ might involve units from $\mathcal{G}_1$ and  hence the network model is no longer a general BTSBM. Nevertheless, we will still investigate the weak consistency since it relies on weaker conditions, and is also an important stepping stone to obtain strong consistency as we will explain later. 

Both the strong and weak consistencies of the recursive spectral clustering under a general BTSBM can be viewed as extensions of the traditional consistency result for spectral clustering \cite[e.g.][]{lei2015consistency} to a misspecified SBM that mistakenly assumes the number of clusters to be $2$ while $K > 2$ in truth.

\subsection{Main results}\label{subsec:main}

Intuitively, the clustering is consistent if $\vec{u}_{n-1}$ is close to $\vec{u}_{n-1}^{*}$. We will show the perturbation bound in $\ell_{\infty}$ and $\ell_2$ norms under a weakly assortative general BTSBM, implying the strong and weak consistencies, respectively. The detailed proofs are deferred to Section \ref{sec:proofs}. 

\begin{thm}[$\ell_{\infty}$ perturbation] \label{thm:linf_perturb}
In the setting of Theorem \ref{thm:eigen}, and further assume that $\xi = O(1)$, where $\xi$ is defined in Equation \eqref{eq:xi}. Then, for any fixed constant $r > 0$, there exists a constant $C_{\ell_{\infty}}$ that only depends on $r$ and $\xi$, such that
\[
\sqrt{n}\|\vec{u}_{n-1}\sign(\vec{u}_{n-1}^{T}\vec{u}_{n-1}^{*}) - \vec{u}_{n-1}^{*}\|_{\infty} < \min\{\sqrt{n_0/n_1},\sqrt{n_1/n_0}\}
\]
with probability at least $1 - (10K + 4)n^{-r}$, provided the following two conditions:
\begin{align}
&\text{Density gap} &&\min\{n_{0}(p_{0} - p_{\es}), n_{1}(p_{1} - p_{\es})\}\ge C_{\ell_{\infty}}\sqrt{(n_0 p_0 + n_1 p_1)\log n},    
\label{eq:cond_eigen_gap}
\\
&\text{Degree variation} &&
(n(\ubar{p}_{*} - p_{\es}))^{4}\ge C_{\ell_{\infty}}(n\bar{p}^{*})^{3}\log n, \label{eq:cond_pbar_pubar}
\end{align}
where $\ubar{p}_*$ and $\bar{p}^*$ are defined in \eqref{eq:pop_para}.
% \XL{I actually prefer to include $\kappa$ explicitly in these conditions.}\Lihua{Changing $C_2$ to be $C_{\infty}$ that appears to be more meaningful?}
\end{thm}

\begin{rem}[Strong consistency]
Recalling Theorem \ref{thm:eigen} that $\sqrt{n}\min
_{i}|\vec{u}_{n-1, i}^{*}| = \min\{\sqrt{n_0/n_1},\sqrt{n_1/n_0}\}$, Theorem \ref{thm:linf_perturb} implies that the signs of the entries of $\vec{u}_{n-1}^*$ are preserved by $\vec{u}_{n-1}$ with high probability under conditions \eqref{eq:cond_eigen_gap} and \eqref{eq:cond_pbar_pubar}. Therefore, the conditions in Theorem \ref{thm:linf_perturb} imply the strong consistency of the first split. 
\end{rem}

\begin{rem}[Sparse networks]
Theorem \ref{thm:linf_perturb} includes the case of sparse networks. In fact, conditions \eqref{eq:cond_eigen_gap} and \eqref{eq:cond_pbar_pubar} can be simultaneously satisfied if $n_0 \asymp n_1$ and
\[
\ubar{p}_{*} \asymp \bar{p}^{*}\asymp (p_0- p_{\es}) \asymp (p_1 - p_{\es}) = O(\log n / n),
\]
in which case the expected degrees are on the order of $O(\log n)$.
\end{rem}
% In the analysis of a correctly specified SBM with $K = 2$,  between-community probability $p_{\es}$ and within-cluster probabilities $p_0$ and $p_1$, it is often assumed that $p_{\es}, p_0, p_1$ are of the same order $\rho_n$. For strong consistency, $\rho_n$ must be $\Omega(\log n / n)$. In this case, $\ubar{p}_{*} \asymp \bar{p}^{*}\asymp \rho_n$ and thus both  \eqref{eq:cond_eigen_gap} and \eqref{eq:cond_pbar_pubar} are satisfied up to constants. For general BTSBMs with $K > 2$, the condition \eqref{eq:cond_eigen_gap} continues to hold as if the model is the above simplified SBM. If all connection probabilities are all of the same order $O(\rho_n)$,  $\ubar{p}_{*} \asymp \bar{p}^{*}\asymp \rho_n$ and thus the condition \eqref{eq:cond_pbar_pubar} holds up to constants. 

\begin{rem}[Degree variation] 
\label{rem:degree_var}
Here we briefly discuss the degree variation condition \eqref{eq:cond_pbar_pubar}. Consider the SBM with $K=2$ again, where $n_0 \asymp n_1$, $p_{\es}=O(\log n / n)$ and $p_0 = n^{-\gamma_0}$, $p_1 = n^{-\gamma_1}$ with constants $\gamma_0, \gamma_1$ satisfying $0<\gamma_0 < \gamma_1< 1$. Then the condition \eqref{eq:cond_eigen_gap} is satisfied up to a constant if $\gamma_1 < \frac{\gamma_0+1}{2}$ and the condition \eqref{eq:cond_pbar_pubar} will hold up to a constant if $\gamma_1 < \frac{3\gamma_0+1}{4}$. We should admit these conditions may be improvable. For example, the two communities might be recoverable if we let $\gamma_0 \rightarrow 0$ while $\gamma_1 \rightarrow 1$.
\end{rem}

\begin{rem}[Multiscale networks]
While most analyses of SBMs focus on the case where all connection probabilities are of the same order, our Theorem \ref{thm:linf_perturb} can deal with multiscale networks where the maximum degree parameter $\bar{p}^{*} \gg p_{\es}$. For example, when $p_{s} \asymp 1$ for every leaf node $s$, it is clear that the maximum and minimum degree parameters satisfy $\ubar{p}_{*}\asymp \bar{p}^{*}\asymp 1$, and both the density gap condition \eqref{eq:cond_eigen_gap} and the degree variation condition \eqref{eq:cond_pbar_pubar} hold for large $n$ even when $p_{\es}, p_{0}, p_{1}\asymp \log n / n$. Our strong consistency result also guarantees adaptivity of graph Laplacian based spectral clustering to multiscale networks with degree heterogeneity, e.g., $p_{\es}, p_{0}, p_{1}\asymp \log n / n$ and $\bar{p}^{*} = n^{-\gamma_0}$ and $\ubar{p}_{*}=n^{-\gamma_1}$ with $0<\gamma_1 < \frac{3\gamma_0+1}{4}<1$.
% Whereas multiscale networks appear to be easier to handle because of the tighter within-cluster connections, it makes the theoretical analysis more challenging. 
\end{rem}

As we will discuss in Section \ref{subsec:challenge}, the theoretical analysis for multiscale networks is more challenging. In fact, to obtain the same bound in Theorem \ref{thm:linf_perturb}, the off-the-shelf $\ell_{\infty}$ perturbation bound by \cite{lei2019unified} on the Fiedler vector requires 
\[n(\min\{p_0, p_1\} - p_{\es})\gtrsim \sqrt{n\bar{p}^{*}\log n}.\]
When $p_{\es}, p_0, p_1\asymp \log n / n$, $\bar{p}^{*}$ must be $O(\log n / n)$, excluding any multiscale network. For example, consider the general BTSBM with $K = 4$, equal community sizes, and
\begin{equation*}
\mtx{B}=\begin{bmatrix}
p_{00} & p_0 &  p_\es & p_\es\\
p_0  & p_{01} & p_\es & p_\es \\
p_\es & p_\es & p_{10} & p_1 \\
p_\es & p_\es & p_{1} & p_{11} \\
\end{bmatrix}.
\end{equation*}
If we further assume $p_{00} = p_{01} = p_{10} = p_{11}=p^*$, then the maximum degree parameter satisfies $\bar{p}^{*}\asymp p^{*}$ due to weak assortativity, and thus the above eigengap condition implies
$
n\min\{p_0, p_1\}\gtrsim \sqrt{np^{*}\log n}.
$
It does not hold if $p_{\es}, p_{0}, p_{1}\asymp \log n / n$ but $p^{*} \succ \log n / n$. Note that $p^{*}$ measures the connection probability within the primitive communities, the most connected groups on the hierarchy. It is hence disappointing and unrealistic to restrict $p^{*}$ into the same order as $p_{\es}$, $p_0$ and $p_1$. 
% connection probability between the least connected groups on the hierarchy. 
In Section \ref{subsec:challenge}, we will explain why multiscale networks are challenging to work with in theory, and in Section \ref{subsec:ideas} we will explain briefly how these challenges can be addressed. 
% \XL{To add a paragraph to discuss the degree variation condition. We can discuss the matrix $\mtx{B}$, with $p_{00}=p_{01} \asymp \bar{p}^*$ and $p_{01}=p_{11} \asymp \ubar{p}_*$.}

The next result gives an $\ell_2$ perturbation bound for the Fiedler vector. 
\begin{thm}[$\ell_2$ perturbation] 
\label{thm:l2_perturb}
Under the same setting of Theorem \ref{thm:eigen}, for any fixed $r,c > 0$,
\[
\|\vec{u}_{n-1}\sign(\vec{u}_{n-1}^{T}\vec{u}_{n-1}^{*}) - \vec{u}_{n-1}^{*}\|_2 < c
\]
with probability at least $1 - 2n^{-r}$, provided that
\begin{equation}
\label{eq:cond_l2}
\min\{n_{0}(p_{0} - p_{\es}), n_{1}(p_{1} - p_{\es})\}\ge C_{\ell_{2}}\sqrt{(n_0 p_0 + n_1 p_1)\log n}
\end{equation}
where $C_{\ell_{2}}$ is a sufficiently large constant that only depends on $r$ and $c$. 
% \XL{Does is depend on $\kappa$ as well? In fact, I prefer to include $\kappa$ explicitly in this density gap condition.}\Lihua{I don't think this result depends on $n_1 / n_0$. Could you confirm that? Changing $C_0$ to $C_2$ that appears to be more meaningful?}
\end{thm}

\begin{rem}[Beyond hierarchical SBM]
\label{rem:weak_general}
While Theorem \ref{thm:l2_perturb} is stated for general BTSBMs, the result holds for a much broader class of networks such that
\[\mtx{P}_{ij} \left\{\begin{array}{ll}
    = p_{\es} & (i\in \mathcal{G}_0, j\in \mathcal{G}_1) \\
    \ge p_0 & (i,j\in \mathcal{G}_0)\\
    \ge p_1 & (i,j\in\mathcal{G}_1)
\end{array}\right.,\]
where $(p_{\es}, p_{0}, p_{1})$ satisfies the condition \eqref{eq:cond_l2}. The result will be stated formally in the supplement. 
\end{rem}

\begin{rem}[Balancedness]
Theorem \ref{thm:l2_perturb} yields an bound on the misclustering error for the first split. Unlike the $\ell_{\infty}$ perturbation bound in Theorem \ref{thm:linf_perturb}, it does not require $\xi = O(1)$. However, the misclustering rate result may rely on certain balancedness. Assume $\vec{u}_{n-1}^{T}\vec{u}_{n-1}^{*} \ge 0$ without loss of generality and let $\mathcal{M} = \{i: \sign(\vec{u}_{n-1, i})\not = \sign(\vec{u}_{n-1, i}^{*})\}$. Obviously, the misclustering error is $|\mathcal{M}|/n$. By Theorem \ref{thm:eigen} (1), for each $i\in \mathcal{M}$, we have
\[
|\vec{u}_{n-1, i} - \vec{u}_{n-1, i}^{*}|\ge \frac{1}{\sqrt{n}}\min\left\{\sqrt{\frac{n_{1}}{n_{0}}}, \sqrt{\frac{n_{0}}{n_{1}}}\right\}\ge \frac{1}{\sqrt{n\xi}},
\]
where the last inequality uses the fact that $\min\{n_1 / n_0, n_0 / n_1\} \ge \min\{n_1, n_0\} / n\ge 1/\xi$. As a result,
\[
\frac{|\mathcal{M}|}{n}
\le \xi \sum_{i\in \mathcal{M}}(\vec{u}_{n-1, i} - \vec{u}_{n-1, i}^{*})^2 
\le \xi\|\vec{u}_{n-1} - \vec{u}_{n-1}^{*}\|_2^2
\le \xi c^2.
\]
When $\xi = O(1)$, Theorem \ref{thm:l2_perturb} implies that $\xi c^2$ can be arbitrarily small when $C_{\ell_2}$ is sufficiently large. 
\end{rem}
\begin{rem}[Relaxed degree variation condition]
The density gap condition \eqref{eq:cond_l2} is essentially the same as that for strong consistency, i.e., \eqref{eq:cond_eigen_gap}. However, the degree variation condition \eqref{eq:cond_pbar_pubar}, which is required for strong consistency, is not required for weak consistency.
\end{rem}

\begin{rem}[$O(\log n)$ degrees]
For an SBM with $K = 2$ communities, the connection probabilities $p_{\es}, p_{0}, p_{1}$ only need to be $\omega(1/n)$ for weak consistency \citep{zhao2012consistency, abbe2017community}. Unfortunately, this cannot be achieved by spectral clustering based on the adjacency matrix or graph Laplacian \citep{krzakala2013spectral}, for which the condition \eqref{eq:cond_l2} is required up to constants \citep[e.g.][]{lei2015consistency}. Nevertheless, we conjecture that the regularized spectral clustering, with appropriately chosen level of regularization,  works in this regime. For example, one can consider removing the nodes whose degrees are greater than  $C_0np^*$ for some constant $C_0$.  \cite{le2017concentration} proved that the adjacency matrix for the remaining graph has tighter spectral concentration around its population version.  
% Denote by $\tilde{\mtx{L}}$ the graph Laplacian of the remaining network. Then it is likely that we can control the $\ell_2$ distance between the Fiedler vector of $\tilde{\mtx{L}}$ and the original population Fiedler vector $\vct{u}_{n-1}^*$ by a small constant. 
% This is an interesting direction to pursue in future. 
However, their theory does not directly apply to multiscale hierarchical SBMs. Moreover, it is unclear how the truncation threshold should be chosen in practice since $p^{*}$ is unknown. 
% Given Theorem \ref{thm:l2_perturb} mainly serves as an important stepping stone for the strong consistency, we don't pursue the idea of modified Laplacian in this paper.}
We leave this intriguing research question for future work.
\end{rem}
% It suffices if the conditions are not stronger than those for Theorem \ref{thm:linf_perturb}. 

\subsection{Main challenge to handle multiscale networks}\label{subsec:challenge}
% First, we explain the technical difficulty brought on by multiscale networks.
A standard technique to obtain the $\ell_2$ perturbation bound is via the Davis-Kahan $\sin\Theta$ Theorem (Lemma \ref{davis_kahan}), which implies that 
\[\|\vec{u}_{n-1}\sign(\vec{u}_{n-1}^{T}\vec{u}_{n-1}^{*}) - \vec{u}_{n-1}^{*}\|_2 \lesssim \frac{\|\mtx{L} - \mtx{L}^*\|}{\lambda_{n-2}^{*} - \lambda_{n-1}^{*}}.\]
Similarly, a straightforward application of the recently developed $\ell_{\infty}$ perturbation bound for unnormalized graph Laplacians \citep{lei2019unified} implies that, under additional regularity conditions and substantial simplifications, 
\[\|\vec{u}_{n-1}\sign(\vec{u}_{n-1}^{T}\vec{u}_{n-1}^{*}) - \vec{u}_{n-1}^{*}\|_{\infty} \lesssim \frac{\|\mtx{L} - \mtx{L}^*\|}{\lambda_{n-2}^{*} - \lambda_{n-1}^{*}}\|\vec{u}_{n-1}^{*}\|_{\infty}\lesssim \frac{\|\mtx{L} - \mtx{L}^*\|}{\lambda_{n-2}^{*} - \lambda_{n-1}^{*}}\cdot \frac{1}{\sqrt{n}}\]
where the last inequality is implied by Theorem \ref{thm:eigen} (4). As a consequence, to obtain the $O(1/\sqrt{n})$ $\ell_{\infty}$ perturbation bound in Theorem \ref{thm:linf_perturb} and the $O(1)$ $\ell_2$ perturbation bound in Theorem \ref{thm:l2_perturb} based on these techniques straightforwardly, it is required that
\begin{equation}\label{eq:loose_technique_goal}
\lambda_{n-2}^{*} - \lambda_{n-1}^{*}\gtrsim \|\mtx{L} - \mtx{L}^*\|.
\end{equation}
By Theorem \ref{thm:eigen} (1) and (2), 
\[\lambda_{n-2}^{*} - \lambda_{n-1}^{*} = \min\{n_0(p_0 - p_{\es}), n_1(p_1 - p_{\es})\}.\]
The best available matrix perturbation inequality (Lemma \ref{lem:laplacian_concentration}) shows that
\[\|\mtx{L} - \mtx{L}^*\| \lesssim \sqrt{n\bar{p}^{*}\log n}.\]
This bound cannot be improved when the average degrees of all nodes are the same in order, i.e., $\ubar{p}_{*}\asymp \bar{p}^{*}$.  Therefore, to guarantee \eqref{eq:loose_technique_goal}, it requires that
\[\min\{n_0(p_0 - p_{\es}), n_1(p_1 - p_{\es})\}\gtrsim \sqrt{n\bar{p}^{*}\log n}.\]
As discussed in Section \ref{subsec:main}, the above condition is overly stringent, illustrating that the standard techniques fail to handle multiscale networks.

\subsection{Proof ideas}
\label{subsec:ideas}
To overcome the difficulty, we will introduce different novel techniques for the $\ell_{\infty}$ and $\ell_{2}$ perturbation bounds. We start with the $\ell_{\infty}$ perturbation. As illustrated above, the main hurdle brought on by multiscale networks is that the eigengap $\lambda_{n-2}^{*} - \lambda_{n-1}^{*}$ is local while the perturbation $\|\mtx{L} - \mtx{L}^{*}\|$ is global on the hierarchy. When $\bar{p}^{*} \gg \max\{p_0, p_1\}$, the eigengap for $\vec{u}_{n-1}$ is too small to yield a desirable eigenvector perturbation bound. Nevertheless, Theorem \ref{thm:eigen} (3) and (4) imply that there are at most $K$ eigenvalues below $n\ubar{p}_{*}$. By pigeonhole principle, there exists $j \le K$ such that 
\[\lambda_{n-j+1}^{*} - \lambda_{n-j}^{*} \ge \frac{\lambda_{n-K}^{*} - \lambda_{n-1}^{*}}{K}\ge  \frac{n(\ubar{p}_{*} - p_{\es})}{K}\gtrsim n(\ubar{p}_{*} - p_{\es}),\]
where the last inequality uses the fact that $K \le \max_{s}(n / n_{s}) = \xi = O(1)$. Let $\mtx{U}_j$ (resp. $\mtx{U}_j^{*}$) be the $\R^{n\times j}$ matrix including eigenvectors $\{\vec{u}_{n-1}, \ldots, \vec{u}_{n-j}\}$ (resp. $\{\vec{u}_{n-1}^{*}, \ldots, \vec{u}_{n-j}^{*}\}$). Then the generic $\ell_{2\rightarrow\infty}$ bound proposed by \cite{lei2019unified}, with substantial simplifications, implies that
\[\|\mtx{U}_{j}\mtx{O}_{j} - \mtx{U}_{j}^{*}\|_{\tti}\lesssim \frac{n\bar{p}^{*}\sqrt{n\bar{p}^{*}\log n}}{\{n(\ubar{p}_{*} - p_{\es})\}^2}\|\mtx{U}_j^{*}\|_{\tti},\]
where $\mtx{O}_j$ is an orthogonal matrix. The condition \eqref{eq:cond_pbar_pubar} and Theorem \ref{thm:eigen} (4) imply that 
\[\|\mtx{U}_{j}\mtx{O}_{j} - \mtx{U}_{j}^{*}\|_{2\rightarrow\infty}\lesssim \frac{1}{\sqrt{n}}.\]
Assuming $\mtx{O}_j = I$, the above bound implies the desired $\ell_{\infty}$ perturbation bound for $\vec{u}_{n-1}$ since $\|\vec{u}_{n-1} - \vec{u}_{n-1}^{*}\|_{\infty}\le \|\mtx{U}_j - \mtx{U}_j^{*}\|_{\tti}$. This heuristic can be made rigorous by applying the $\ell_2$ perturbation bound given in Theorem \ref{thm:l2_perturb} and Davis-Kahan $\sin\Theta$ Theorem, which show that $\mtx{O}_j\approx I$ in some appropriate sense. In sum, we deduce the $\ell_{\infty}$ perturbation bound for the Fiedler vector from a generic $\ell_{\tti}$ perturbation bound for a larger eigenspace with a large eigengap. This also illustrates why we need the entire eigenstructure of $\mtx{L}^{*}$ even if the algorithm merely uses the Fiedler vector.

As shown above, the $\ell_2$ perturbation bound is key to establish the $\ell_\infty$ perturbation bound. As aforementioned, the direct application of Davis-Kahan $\sin \Theta$ Theorem fails because the matrix perturbation error $\|\L - \L^{*}\|$ is too large compared to the eigengap. However, instead of viewing $\L^{*}$ as the target and $\L - \L^{*}$ as the perturbation, we can replace $\L^{*}$ by any matrix $\td{\L}$ with $\vec{u}_{n-1}(\td{\L}) = \vec{u}_{n-1}^{*}$. The Davis-Kahan $\sin \Theta$ Theorem would imply a tighter bound if our selected $\widetilde{\mtx{L}}$ satisfies
\[\|\L - \td{\L}\| \ll \|\L - \L^{*}\|.\]
Typically, it is difficult to construct an explicit $\td{\L}$ without further structural assumptions on $\L^{*}$. However, we observe an intriguing property of unnormalized graph Laplacians that enables an easy construction of $\td{\L}$.

% \begin{lem}\label{lem:laplacian_property}
% Let $\mtx{\cA}\in \R^{n\times n}$ be the adjacency matrix of any weighted non-negative graph with $\mtx{\cA}_{ij} = p$ for any $i\in \mathcal{G}_0, j\in \mathcal{G}_1$ and $\mtx{\cA}_{ij} > p$ otherwise, where $p\in (0, 1]$ and $(\mathcal{G}_0, \mathcal{G}_1)$ forms a partition of $\{1, \ldots, n\}$. Further, let $\mtx{\cL}$ be the unnormalized graph Laplacian for $\mtx{\cA}$. Then
% \[\lambda_{n-1}(\mtx{\cL}) = np, \quad \text{and }\,\, \vec{u}_{n-1, i}(\mtx{\cL}) = \pm\left\{\begin{array}{ll}
%      \sqrt{|\mathcal{G}_{1}| / |\mathcal{G}_{0}| n} & (i\in \G_{0})\\
%       -\sqrt{|\mathcal{G}_{0}| / |\mathcal{G}_1| n} & (i\in \G_{1})
%     \end{array}
% \right..\]
% \end{lem}

\begin{lem}\label{lem:laplacian_property}
Let $\td{\L}$ be the unnormalized graph Laplacian of the pseudo-adjacency matrix $\widetilde{\A}$ that replaces the between-community edges by their common expectation $p_{\es}$, i.e.,
\[\widetilde{\A}_{ij} = \left\{\begin{array}{ll}
     \A_{ij} & (i, j\in \mathcal{G}_0, \text{ or }i,j\in\mathcal{G}_1)  \\
     p_{\es} & (\text{otherwise})
\end{array}\right..\]
Under the same setting as in Theorem \ref{thm:l2_perturb}, if $C_{\ell_2}$ is a sufficiently large constant that only depends on $r$, with probability at least $1 - n^{-r}$,
\[\lambda_{n-1}(\td{\L}) = np_{\es} \text{ with multiplicity }1, \quad \text{and }\,\, \vec{u}_{n-1, i}(\td{\L}) = \vec{u}_{n-1}^{*}.\]
\end{lem}

% Lemma \ref{lem:laplacian_property} implies that the $(n-1)$-th eigenpair of the unnormalized graph Laplacian does not depend on the within-community connection probabilities if the between-community ones are homogeneous, as in the general BTSBM. As a consequence, we can take $\td{\L}$ as the unnormalized graph Laplacian of the pseudo-adjacency matrix $\td{\A}$ that replaces the between-community edges by their common expectation $p_{\es}$, i.e.,
% \[\td{\A}_{ij} = \left\{\begin{array}{ll}
%      \A_{ij} & (i, j\in \mathcal{G}_0, \text{ or }i,j\in\mathcal{G}_1)  \\
%      p_{\es} & (\text{otherwise})
% \end{array}\right.\]
By Lemma \ref{lem:laplacian_property}, $\td{\L}$ preserves the $(n-1)$-th eigenpair with high probability. Meanwhile, $\L$ and $\td{\L}$ only differ in the between-mega-community entries which are only determined by $p_{\es}$. When $p_{\es} \ll \bar{p}^{*}$, it turns out that
\[\|\L - \td{\L}\|\lesssim \sqrt{np_{\es}\log n} \ll \|\L - \L^{*}\|\asymp \sqrt{n\bar{p}^{*}\log n}.\]
To apply Davis-Kahan $\sin \Theta$ Theorem on the decomposition $\L = \td{\L} + (\L - \td{\L})$, we still need to bound the eigengap of $\td{\L}$. Since $\td{\L}$ is a graph Laplacian, it also preserves the $n$-th eigenvalue. Therefore, it remains to bound $\lambda_{n-2}(\td{\L})$ from below. Note that $\E[\td{\L}] = \L^{*}$. A natural lower bound can be obtained via Weyl's inequality:
\[\lambda_{n-2}(\td{\L})\ge \lambda_{n-2}(\L^{*}) - \|\td{\L} - \L^{*}\|.\]
However, the multiscale issue persists because $\lambda_{n-2}(\L^{*})$ only involves $(p_{\es}, p_{0}, p_{1})$ while $\|\td{\L} - \L^{*}\|$ involves all connection probabilities.

The only exception is when the binary tree $\T$ only involves two leaf nodes $0$ and $1$, in which case the multiscale issue disappears. As a result, it is sufficient to show that a deeper tree always increases $\lambda_{n-2}(\td{\L})$ compared to the simple $2$-leaf tree. To this end, we can generate another graph $\widetilde{\A}'$ by dampening the entries $\widetilde{\A}$, with $\E[\widetilde{\A}']$ being the connection probability matrix corresponding to the $2$-leaf tree. This can be achieved by multiplying $\widetilde{\A}_{ij}$ by a Bernoulli random variable with parameter $p_s / P_{ij}$ for any $i, j\in \mathcal{G}_s$. Recalling that the unnormalized graph Laplacian becomes larger in the positive semidefinite ordering when the adjacency matrix increases entrywise (Lemma \ref{lem:laplacian_monotonicity}), we obtain that $\tilde{\L} \succeq \td{\L}'$, where $\td{\L}'$ is the graph Laplacian for $\widetilde{\A}'$. By Weyl's inequality, 
\[
\lambda_{n-2}(\tilde{\L})\ge \lambda_{n-2}(\td{\L}') + \lambda_{n}(\tilde{\L} - \td{\L}') \ge \lambda_{n-2}(\td{\L}').
\]
Therefore, we reduce the general BTSBM to the simple $2$-leaf case and hence avoid the multiscale issue completely.

\subsection{Comparison with previous theoretical results}

Although consistency of graph Laplacian-based spectral clustering under the general BTSBM has been studied by \cite{balakrishnan2011noise}, their regularity conditions only hold for dense networks. In particular, they consider weighted graphs with sub-Gaussian weights, including the Bernoulli weight in a network as a special case whose sub-Gaussian parameter is $1$. According to Theorem 1 in \cite{balakrishnan2011noise}, to recover the first level of the hierarchy, they require $\gamma^4\sqrt[4]{n / \log n} = \omega(1)$ where $\gamma = \min\{p_0, p_1\}- p_{\es}$.  As a result, the minimal expected degree $n\ubar{p}_{*} = \Omega(n\min\{p_0, p_1\}) =  \omega(n^{15/16})$.
Therefore, their strong consistency guarantee is only valid for very dense networks. In contrast, our Theorem \ref{thm:linf_perturb} holds for sparse networks with $n\bar{p}^{*} = \Omega(\log n)$.

On the other hand, \cite{li2018hierarchical} and \cite{lei2019unified} derive the analogue of Theorem \ref{thm:linf_perturb} for the adjacency matrix under a restrictive BTSBM where $\T$ is a full balanced binary tree and $(p_s, n_s) = (p_{s'}, n_{s'})$ whenever $|s| = |s'|$. Both analyses work for sparse networks; the former allows the expected degree to be $O(\log^{2+\epsilon}n)$ for $\epsilon > 0$ while the latter improves the dependence to the critical regime $O(\log n)$. A crucial property of the restrictive BTSBM is the strict homogeneity in the expected degrees, for which the population Laplacian has the same eigenvectors as the expected adjacency matrix. In addition, both works consider the traditional setting where all connection probabilities are of the same order and hence exclude multiscale networks. Therefore, our analysis can be viewed as a substantial generalization.

\section{Experiments}\label{sec:simulation}
%  section 4 experiments
\subsection{Synthetic networks}

\begin{figure}[h]
     \centering
     \begin{subfigure}[h]{0.244\textwidth}
         \centering
         \includegraphics[width=\textwidth]{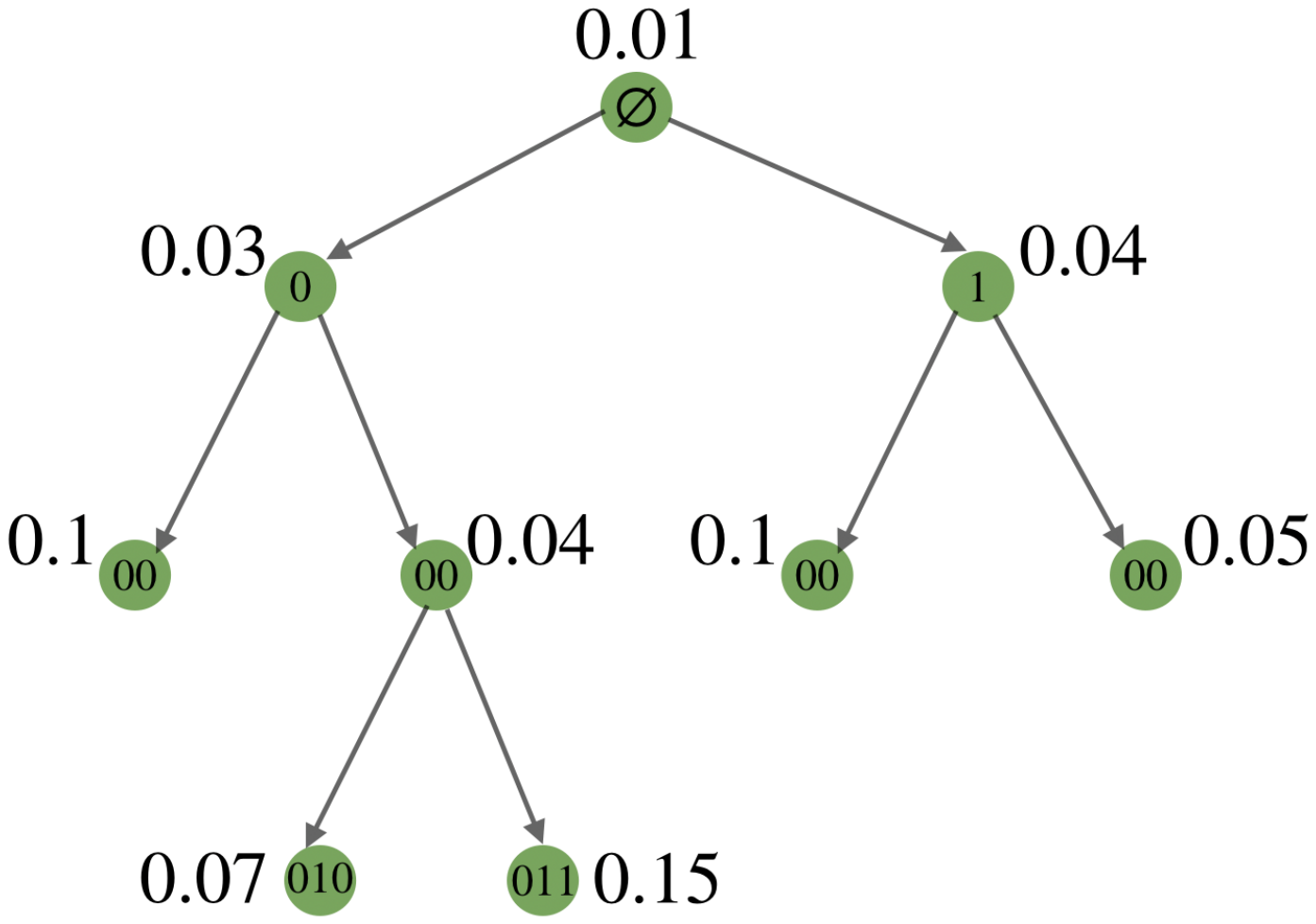}
     \end{subfigure}
     \begin{subfigure}[h]{0.244\textwidth}
         \centering
         \includegraphics[width=\textwidth]{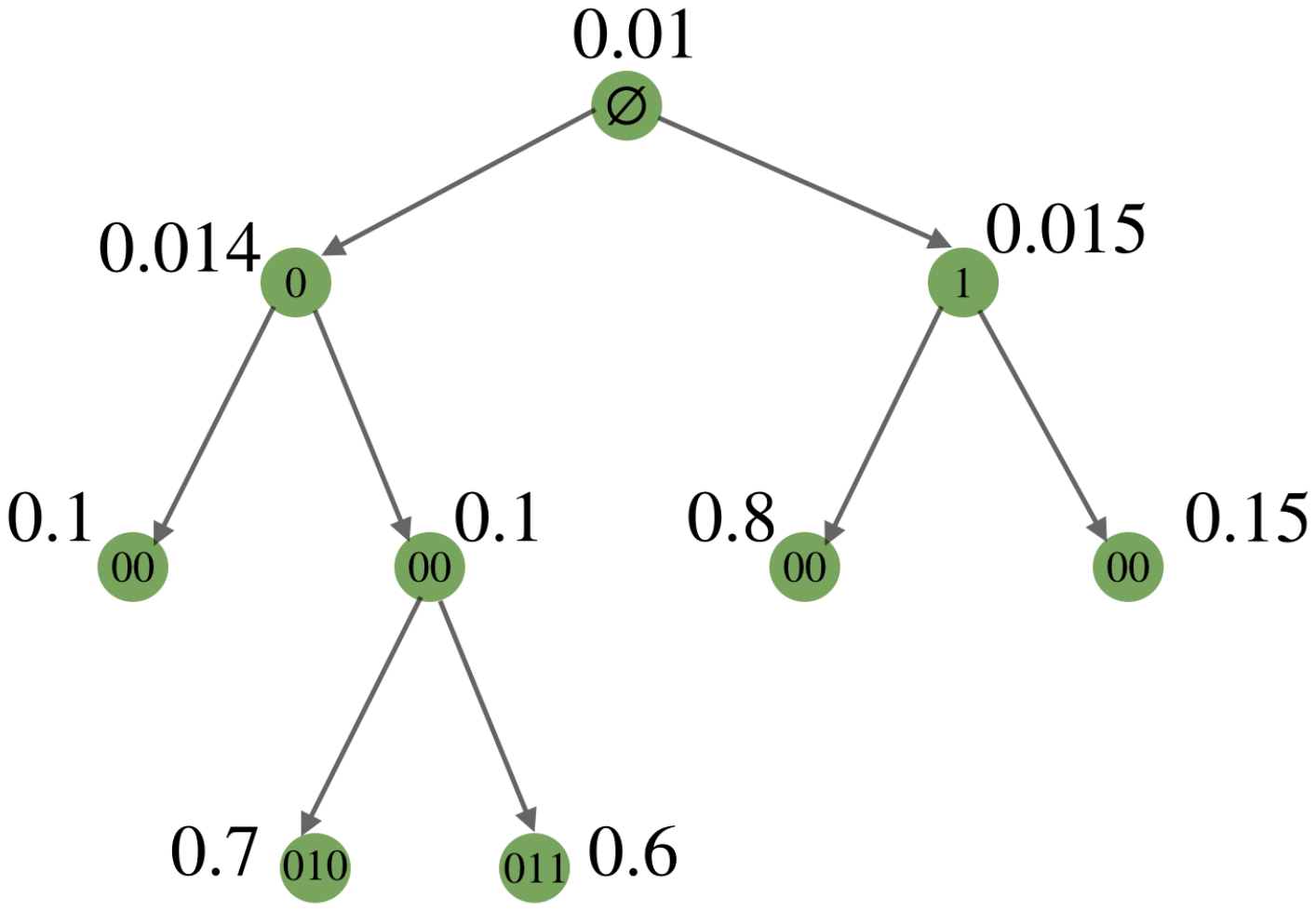}
     \end{subfigure}
     \begin{subfigure}[h]{0.244\textwidth}
         \centering
         \includegraphics[width=\textwidth]{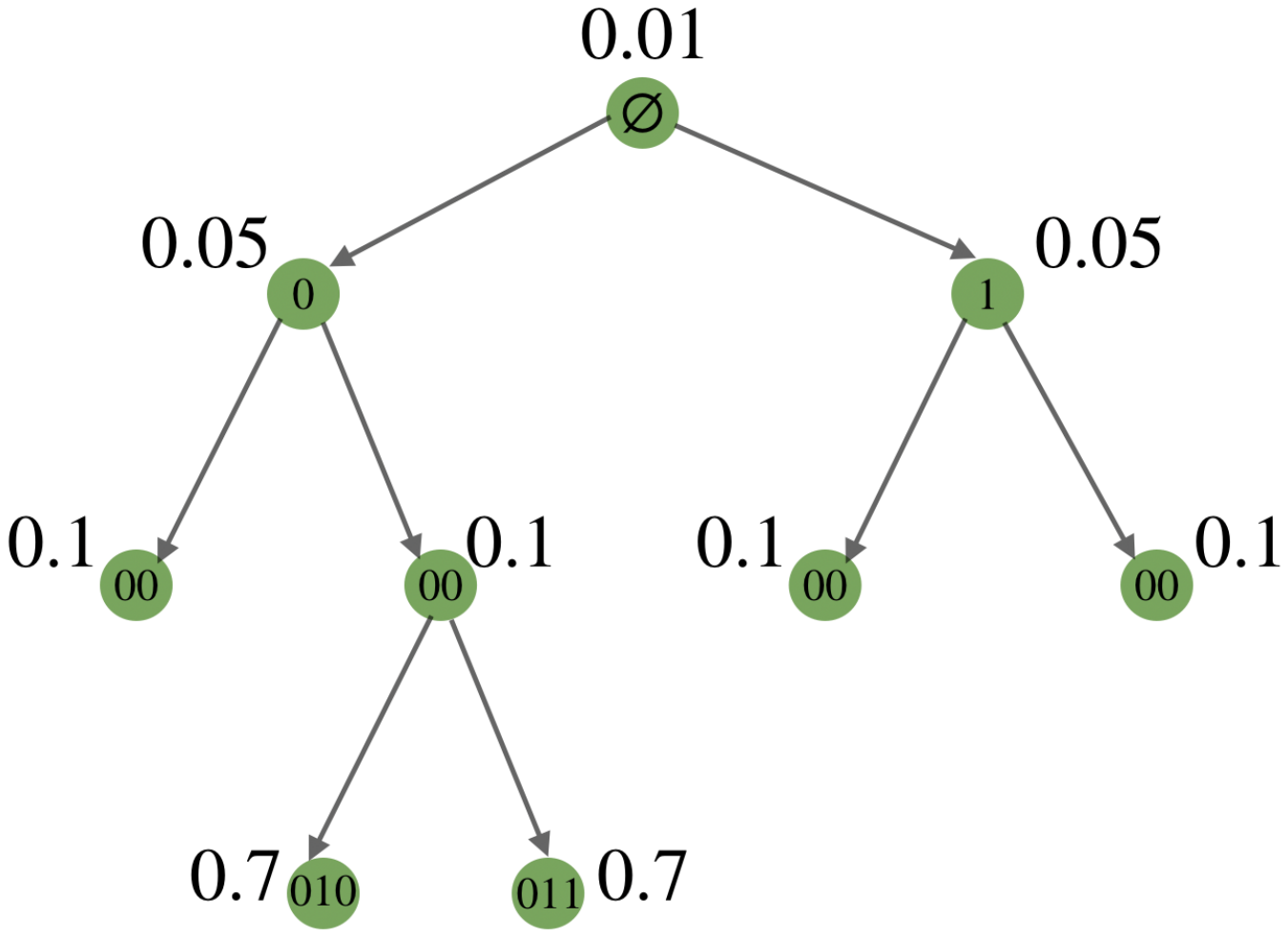}
     \end{subfigure}
     \begin{subfigure}[h]{0.244\textwidth}
         \centering
         \includegraphics[width=\textwidth]{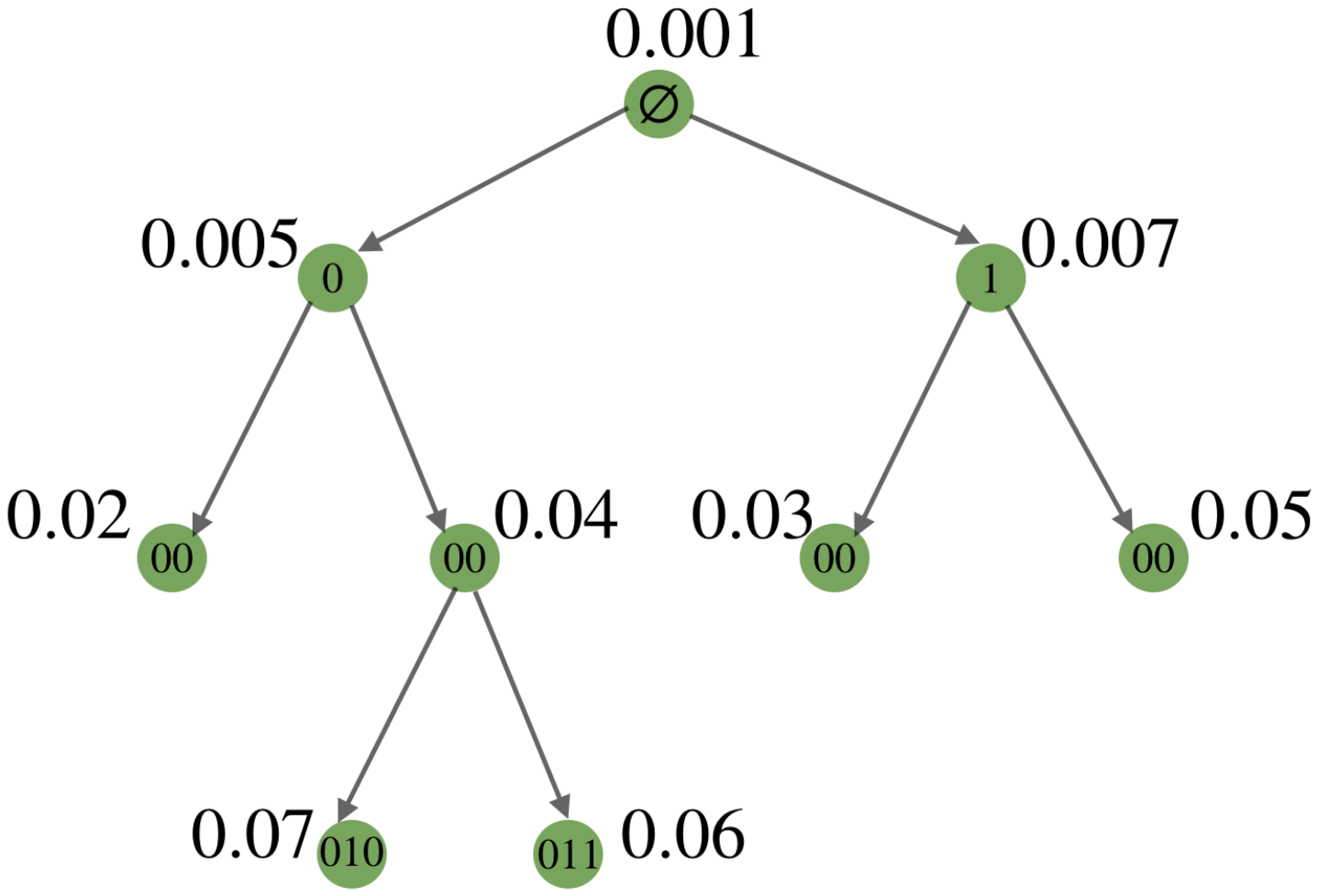}
     \end{subfigure}
     \begin{subfigure}[h]{0.24\textwidth}
         \centering
         \includegraphics[width=\textwidth]{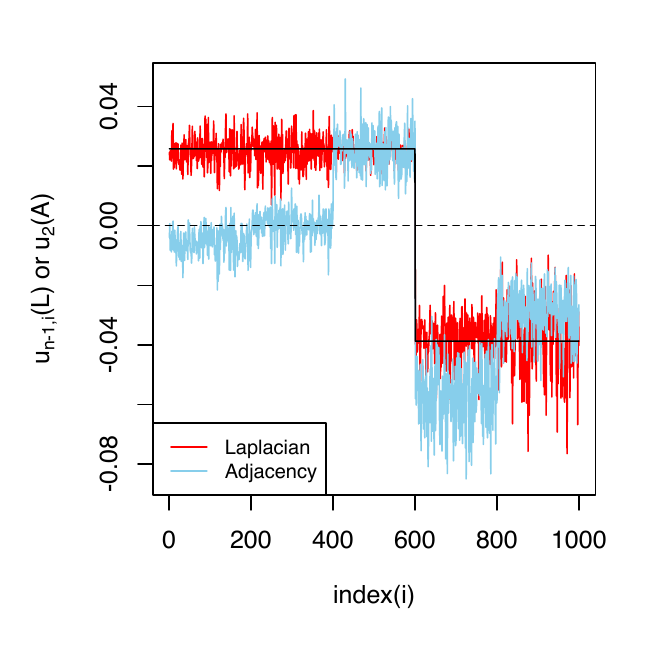}
         \caption{Model 1}
         \label{fig:simu1}
     \end{subfigure}
      \hfill
      \begin{subfigure}[h]{0.24\textwidth}
         \centering
         \includegraphics[width=\textwidth]{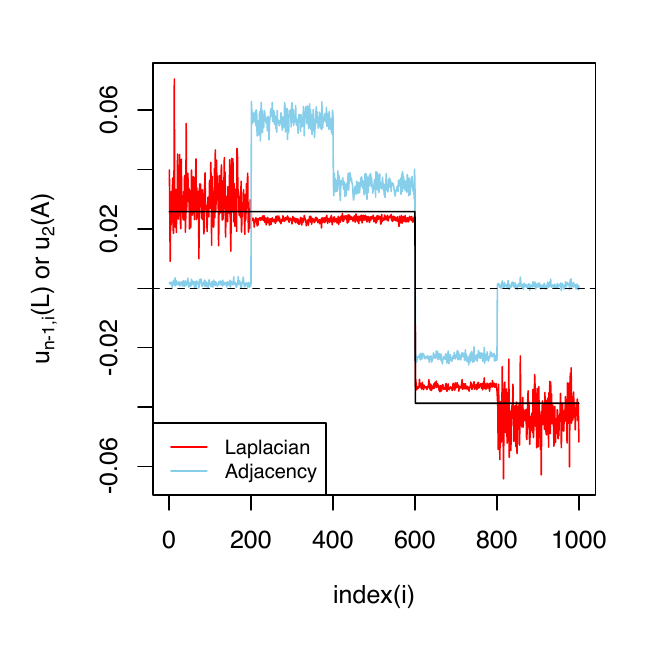}
         \caption{Model 2}
         \label{fig:simu2}
     \end{subfigure}
     \begin{subfigure}[h]{0.24\textwidth}
         \centering
         \includegraphics[width=\textwidth]{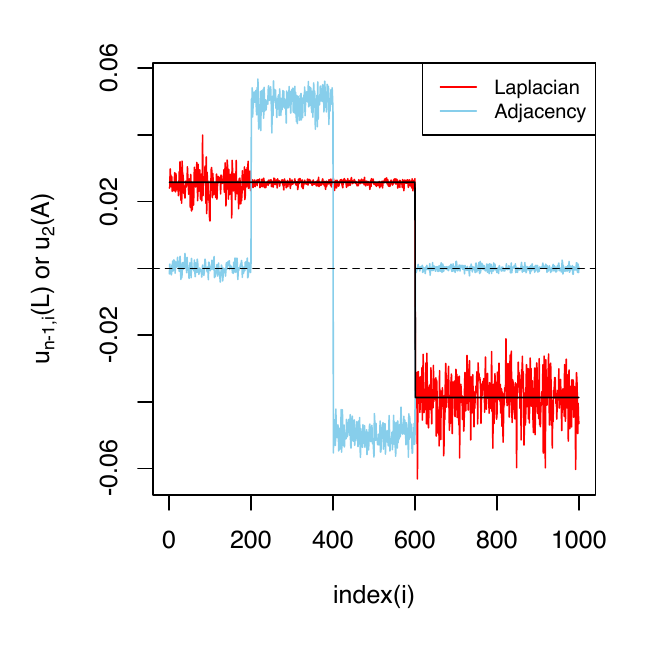}
         \caption{Model 3}
         \label{fig:simu3}
     \end{subfigure}
     \hfill
     \begin{subfigure}[h]{0.24\textwidth}
         \centering
         \includegraphics[width=\textwidth]{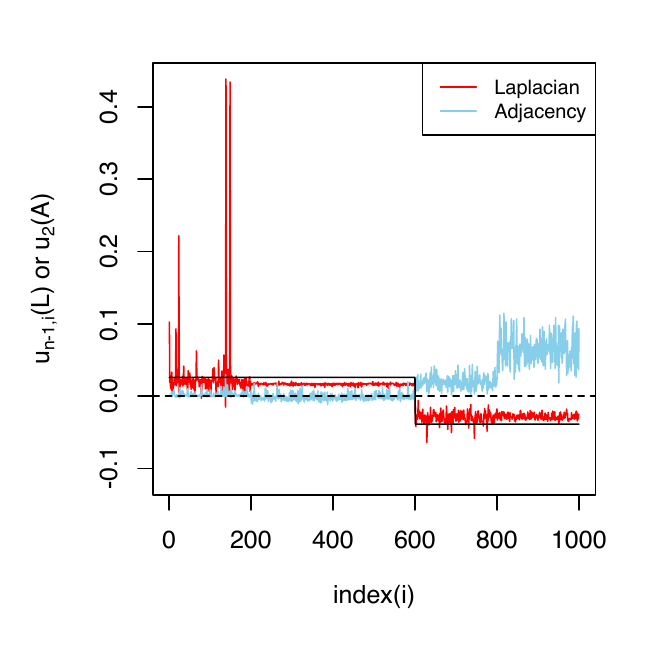}
         \caption{Model 4}
         \label{fig:simu4}
     \end{subfigure}
    \caption{general BTSBMs (top) and the associated eigenvectors (bottom), including $\vec{u}_{n-1}(\L^{*})$ (black), $\vec{u}_{n-1}(\L)$ (red), and the $\vec{u}_2(\A)$ (blue)}
    \label{fig:simu}
\end{figure}

We generate synthetic networks from general BTSBMs with $\T$ being the simple binary tree presented in Figure \ref{fig:illustrate_hbm} and different connection probabilities presented in the top panels of Figure \ref{fig:simu}. Each model has $n = 1000$ units and $200$ units in each primitive cluster. The bottom panels of Figure \ref{fig:simu} compares $\vec{u}_{n-1}(\L)$, the Fiedler eigenvector of the unnormalized graph Laplacian, and $\vec{u}_{n-1}(\L^{*})$, the population Fiedler vector. This can be viewed as an empirical check of the $\ell_{\infty}$ perturbation bound (Theorem \ref{thm:linf_perturb}) and the strong consistency of the first split. As a comparison, we also plot $\vec{u}_2(\A)$, the eigenvector corresponding to the second largest eigenvalue of the adjacency matrix, which is considered in \cite{li2018hierarchical} and \cite{lei2019unified}.

Figure \ref{fig:simu1} shows a setting where the connection probabilities are roughly of the same order. The signs of the Fiedler vector $\vec{u}_{n-1}(\L)$ perfectly align with the mega-community memberships given by the first split, while $\vec{u}_{2}(\A)$ messes up with the mega-community $\mathcal{G}_0$.

Figure \ref{fig:simu2} considers a highly multiscale setting where $\bar{p}^* \gg \ubar{p}_{*}$. This poses a potential threat to the degree variation condition \eqref{eq:cond_pbar_pubar} in Theorem \ref{thm:linf_perturb}. While the $\ell_{\infty}$ error of the Fiedler vector grows substantially compared to Figure \ref{fig:simu1}, the signs of the entries still perfectly identify the first split. An intriguing observation is the asymmetry of entrywise errors; whereas the overall $\ell_{\infty}$ perturbation error is sufficiently large to flip the sign of an entry, it is mainly contributed by ``outbound'' deviations that have no effect on the sign while ``inbound'' deviations that pull entries to the other side of the axis are much smaller. This illustrates the potential suboptimality of our strategy to deduce the strong consistency from a small $\ell_{\infty}$ perturbation error. The phenomenon has been studied by \cite{abbe2017entrywise} and \cite{deng2020strong} for standard SBMs with $K = 2$ and by \cite{lei2019unified}
 for BTSBM. We leave the refined analysis for general BTSBM for future work. 

Figure \ref{fig:simu3} examines a slight misspecification of the balanced BTSBM studied in \cite{li2018hierarchical} where the connection probabilities are identical within each level. Then $\vec{u}_2(\A)$ fails to identify the first split while the Fiedler vector corresponding to the graph Laplacian works as desired. This illustrates the sensitivity of adjacency matrix-based spectral clustering to the model misspecification.

Figure \ref{fig:simu4} presents a setting with tiny connection probabilities, resulting in a very sparse network. In this case, the Fiedler vector exhibits a few spikes, thereby forcing the other values to be small. We observe this pattern frequently in repeated experiments. This suggests that the eigengap condition \eqref{eq:cond_eigen_gap} fails to hold. In this case, it is not surprising that the performance, in terms of strong or weak consistency, degrades drastically. Our observation is in line with \cite{krzakala2013spectral} that spectral clustering based on the adjacency matrix or graph Laplacians cannot handle networks that are too sparse.

\begin{figure}[h]
     \centering
         \includegraphics[width=0.25\textwidth]{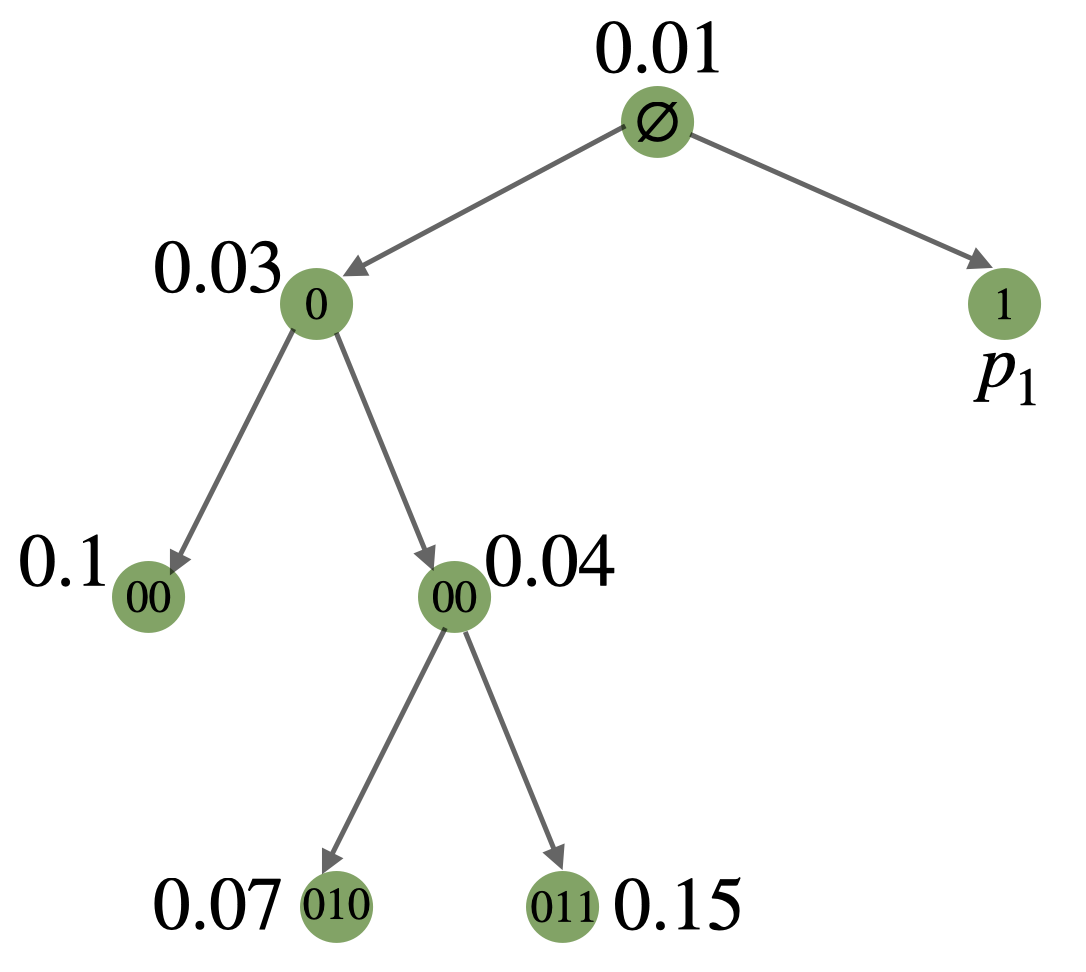}
         \caption{Unbalanced tree}
         \label{fig:untree}
\end{figure}
\begin{figure}[h]
     \centering   
     \begin{subfigure}[h]{0.24\textwidth}
         \centering
         \includegraphics[width=\textwidth,height=3.3cm]{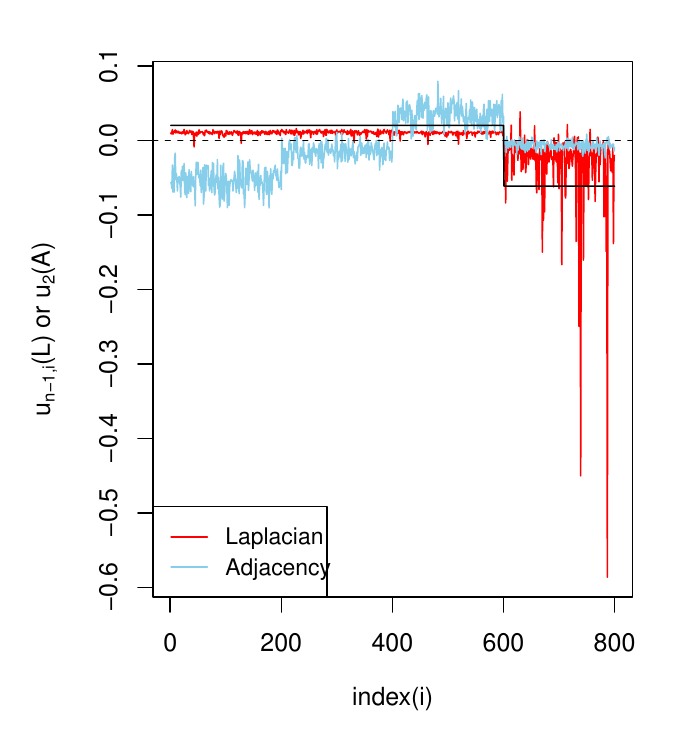}
         \caption{$p_1=0.04$}
         \label{fig:untree1}
     \end{subfigure}
     \hfill
      \begin{subfigure}[h]{0.24\textwidth}
         \centering
         \includegraphics[width=\textwidth,height=3.3cm]{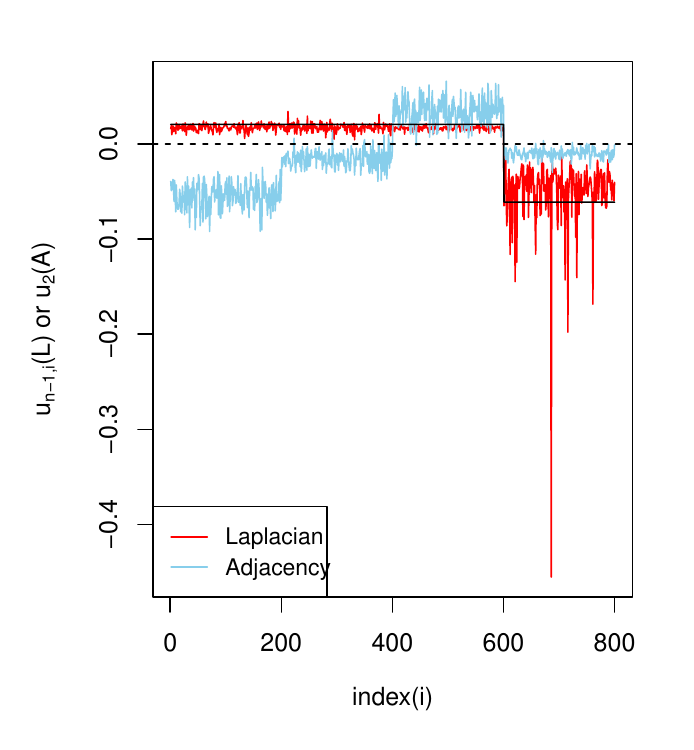}
         \caption{$p_1=0.06$}
         \label{fig:untree2}
     \end{subfigure}
     \hfill
      \begin{subfigure}[h]{0.24\textwidth}
         \centering
         \includegraphics[width=\textwidth,height=3.3cm]{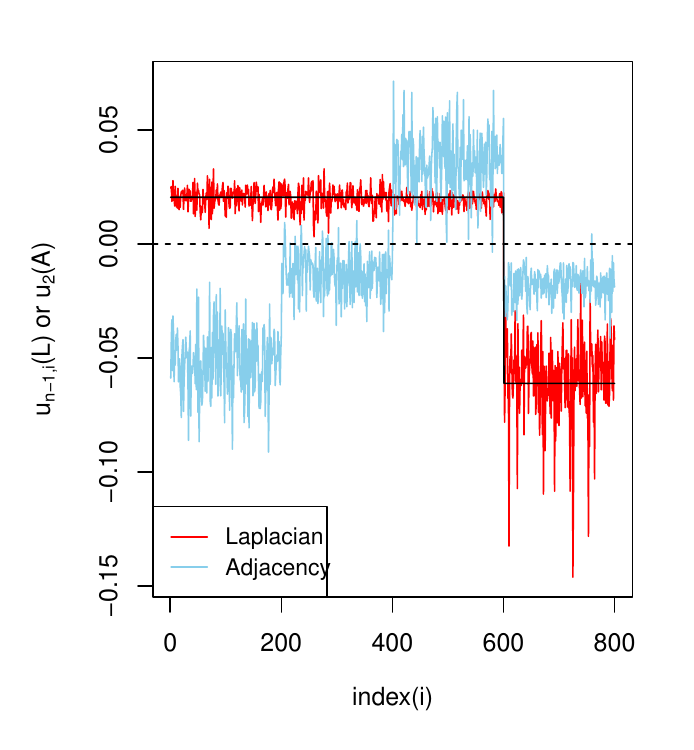}
         \caption{$p_1=0.08$}
         \label{fig:untree3}
     \end{subfigure}
     \hfill
      \begin{subfigure}[h]{0.24\textwidth}
         \centering
         \includegraphics[width=\textwidth,height=3.3cm]{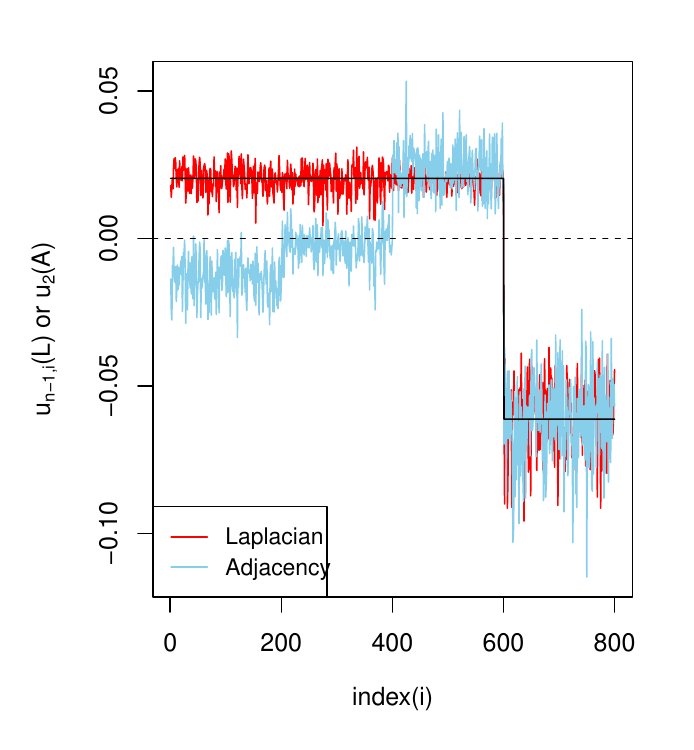}
         \caption{$p_1=0.1$}
         \label{fig:untree4}
     \end{subfigure}
    \caption{Eigenvectors including $\vec{u}_{n-1}(\L^{*})$ (black), $\vec{u}_{n-1}(\L)$ (red), and the $\vec{u}_2(\A)$ (blue) under different values of $p_1$.}
    \label{fig:simulation unbalanced tree}
\end{figure}

We also consider the case when the binary tree in general BTSBM is more unbalanced. We generate synthetic networks from general BTSBMs with $\T$ being a slightly different binary tree with that in Figure \ref{fig:illustrate_hbm}, where we have three primitive clusters versus one primitive clusters in the first split. The tree structure and associated connection probabilities are presented in Figure \ref{fig:untree}. The number of nodes in each primitive cluster is still 200. We compare the aforementioned eigenvectors under different values of $p_1$ in Figure \ref{fig:simulation unbalanced tree}. 

From Figure \ref{fig:simulation unbalanced tree}, we can observe that, as $p_1$ increases, the Fiedler vector is less spiky and the deviation of $\vec{u}_{n-1}(\L)$ from $\vec{u}_{n-1}(\L^{*})$ becomes smaller. In this case, $n_1$ is much smaller than $n_0$, rendering the eigengap condition \eqref{eq:cond_eigen_gap} harder to satisfy. Nevertheless, as long as $p_1$ exceeds certain threshold, the eigengap condition becomes plausible and the Fiedler vector is able to perfectly identify the mega-communities $\mathcal{G}_0$ and $\mathcal{G}_1$. 
% fails to hold. On the other hand, if the eigengap condition \eqref{eq:cond_eigen_gap} is valid, the Fielder vector is able to perfectly bi-partition the network even when the underlying binary tree is unbalanced.

\subsection{Real-world networks}
In this section, we compare recursive spectral bi-clustering algorithms based on the adjacency matrix $\mtx{A}$, the unnormalized graph Laplacian $\L=\mtx{D}-\mtx{A}$ and the normalized graph Laplacian $\mtx{N} = \I-\mtx{D}^{-1/2}\mtx{A}\mtx{D}^{-1/2}$ on seven real-world networks, summarized in Table \ref{tab:summary}. 
 All networks contain explicit information regarding the true community memberships, which we use to evaluate the performance of clustering algorithms; see the references in the second column of Table \ref{tab:summary} for more details. 

\begin{table}[h]
\centering 
\begin{tabular}{l|l|c|c|c|c|c|c} 
\hline
 Dataset & Source & $|V|$ & $|E|$ & $K$ & $d_{\min}$ & $d_{\max}$ & $\bar{d}$ \\
    \hline
    Dolphins        & \cite{lusseau2003bottlenose}  & 62 & 159 & 2 & 1 & 12 & 5.129\\
    Karate          & \cite{zachary1977information}  & 34  & 78 & 2 & 1 & 17 & 4.588\\
    Political books & Krebs (unpublished)  & 92  & 374 & 2 & 1 & 24 & 8.130\\ 
    Political blogs  & \cite{adamic2005political}  & 1222 & 16714 & 2 & 1 & 351 & 27.355\\ 
    UK faculty      & \cite{nepusz2008fuzzy}  & 79  &  552 & 3 & 2 & 39 & 13.975\\
    Football        & \cite{girvan2002community}  & 110 & 570 & 11 & 7 & 13 & 10.364\\
    C. elegans       & \cite{jarrell2012connectome} & 229 &  1085 & 6 & 1 & 34 & 9.585\\
\hline
\end{tabular} 
  \caption{Seven network datasets} 
  \label{tab:summary} 
\end{table}

We evaluate the performance of clustering algorithms via the completeness score \citep{rosenberg2007v}, an external entropy-based cluster evaluation measure. A clustering result satisfies the completeness if all vertices that are members of a true (primitive) community reside in the same estimated community. Equivalently, each estimated community from a complete clustering must be the union of a subset of true (primitive) communities. Grouping all of the vertices into a single community is an extreme example of a complete clustering. The completeness score is designed to measure the proximity to completeness. Suppose the true communities are $V_1,\ldots,V_K$ and the estimated communities are $\hat{V}_1,\ldots,\hat{V}_{\hat{K}}$, where $\hat{K}$ might differ from $K$. The completeness score is defined as
\begin{equation}
\label{def:complete}
c(\hat{V}, V) = \left\{\begin{array}{ll}
1 & \text { if } H(\hat{V})=0 \\
1-\frac{H(\hat{V} | V)}{H(\hat{V})} & \text { otherwise }
\end{array}\right.
\end{equation}
where $H( \hat{V}| V)$ is the conditional entropy of the estimated clusters given the true community assignments and $H(\hat{V})$ is the entropy of the estimated clusters, i.e.,
\[
H( \hat{V}| V) =-\sum_{i=1}^{K} \sum_{j=1}^{\hat{K}} \frac{|V_i \cap \hat{V}_j|}{n} \log \frac{|V_i \cap \hat{V}_j|}{|\hat{V}_j|}, \quad
H(\hat{V}) =-\sum_{j=1}^{\hat{K}} \frac{|\hat{V}_j|}{ n} \log \frac{|\hat{V}_j|}{n}.
\]
Clearly, the completeness score \eqref{def:complete} takes value in $[0,1]$ and a value $1$ implies that the clustering is complete. This metric is invariant to label permutations and asymmetric in $V$ and $\hat{V}$. The asymmetry renders the completeness score a proper metric to evaluate the performance of recovering mega-communities.

\begin{table}[h]
\centering 
\begin{tabular}{l|c|c|c} 
 \hline
 Dataset & $\mtx{A}$  & $\L$ & $\mtx{N}$  \\
    \hline
    Dolphins   & 0.470 & \textbf{1} & 0.883\\
    Karate     & \textbf{1}   & 0.840  & 0.840\\
    Political books & 0.823  & \textbf{0.869}  & \textbf{0.869}\\ 
    Political blogs  & \textbf{0.675} & 0.007  & 0.012\\
    UK faculty      & 0.765  &  0.908 & \textbf{1} \\
    Football        & 0.763  &  \textbf{0.802} & \textbf{0.802} \\
    C. elegans       & 0.416  &  \textbf{0.939}  & 0.807 \\
     \hline
\end{tabular} 
  \caption{Completeness scores for the first split} 
  \label{tab:completeness} 
\end{table} 

For each of the three recursive bi-partitioning algorithms, we compute the completeness score for the first split. The results are reported in Table \ref{tab:completeness}. Notably, the unnormalized graph Laplacian-based algorithm performs well on all networks but Political blogs, which has substantially higher degree variation as shown in Table \ref{tab:summary}. This partly corroborates our theory that the degree variation plays an important role. 

Next, we move deeper into the estimated hierarchy and evaluate the performance of (partial) hierarchy recovery. For illustration, we investigate three networks: UK faculty, Football, and C. elegans. We will examine the performance of the first few splits based on the completeness scores.

\begin{figure}[h]
     \centering   
     \begin{subfigure}[h]{0.3\textwidth}
         \centering
         \includegraphics[width=\textwidth]{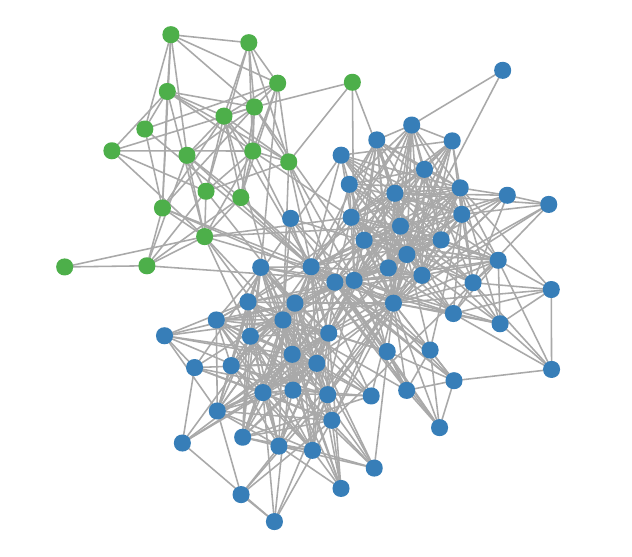}
         \caption{First split (0.908)}
         \label{fig:UK1}
     \end{subfigure}
     \hfill
     \begin{subfigure}[h]{0.3\textwidth}
         \centering
         \includegraphics[width=\textwidth]{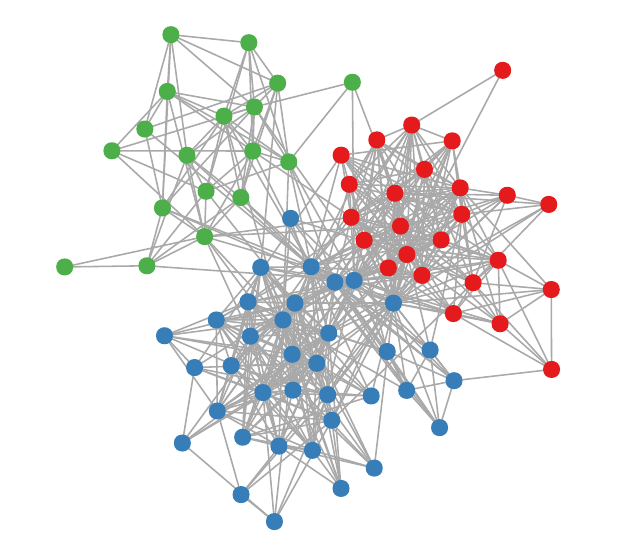}
         \caption{Second split (0.735)}
         \label{fig:UK2}
     \end{subfigure}
      \hfill
     \begin{subfigure}[h]{0.3\textwidth}
         \centering
         \includegraphics[width=\textwidth]{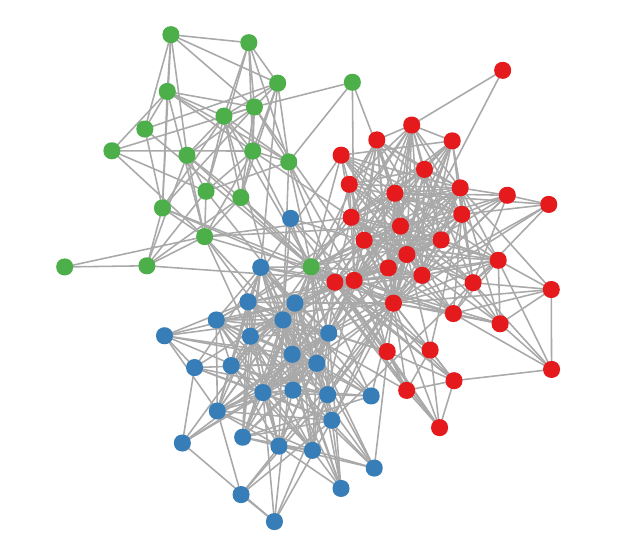}
         \caption{Ground truth}
         \label{fig:UKtruth}
     \end{subfigure}
        \caption{Spectral recursive bi-partition on the UK faculty network, with the completeness scores in the parentheses}
        \label{fig:UK}
\end{figure}

~\\
\textbf{UK faculty} is a personal friendship network of the academic staffs in a UK university, which consists of three separate schools. These three separate schools are treated as the true primitive communities; see Figure \ref{fig:UKtruth}. Figure \ref{fig:UK1} and \ref{fig:UK2} display the first and second splits given by the recursive bi-partitioning algorithm based on the Fiedler vector. The first split separates the green community from the others and achieves a high completeness score 0.908, suggesting that it captures two meaningful mega-communities in this network. Unsurprisingly, The second split has a lower completeness score because the connections within or between the red and blue vertices are similar and thus it is harder to distinguish them. Compared to the ground truth, our algorithm performs reasonably well in recovering both the primitive communities and the hierarchy. 

\begin{figure}[h]
     \centering   
     \begin{subfigure}[h]{0.24\textwidth}
         \centering
         \includegraphics[width=\textwidth]{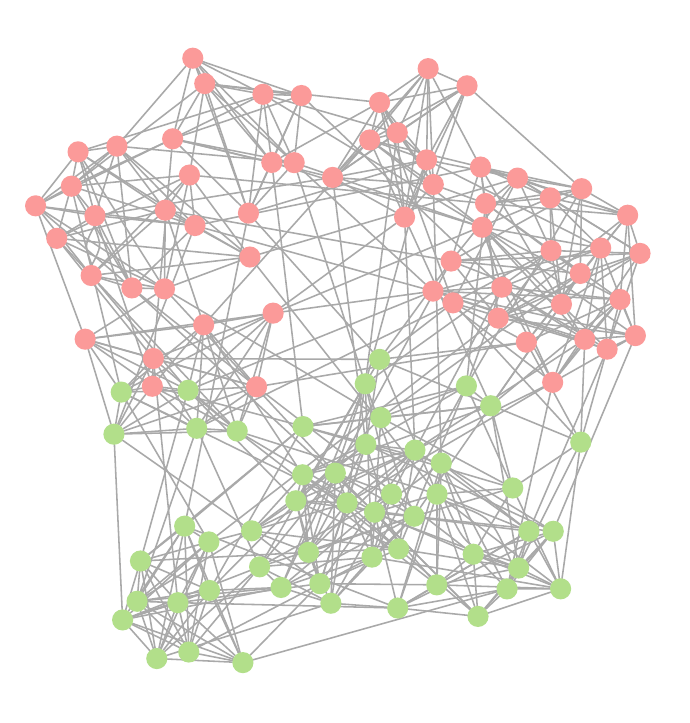}
         \caption{First level (0.802)}
         \label{fig:FB1}
     \end{subfigure}
     \hfill
     \begin{subfigure}[h]{0.24\textwidth}
         \centering
         \includegraphics[width=\textwidth]{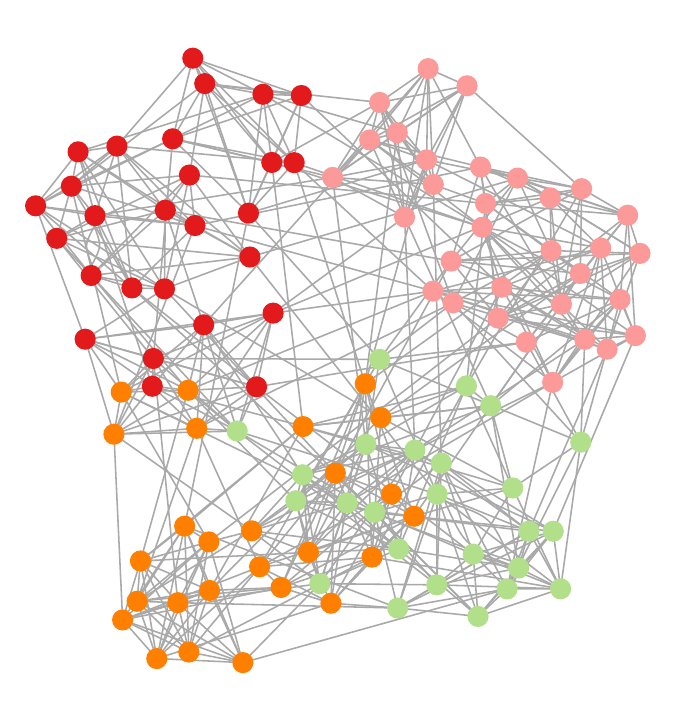}
         \caption{Second level (0.828)}
         \label{fig:FB2}
     \end{subfigure}
      \hfill
    \begin{subfigure}[h]{0.24\textwidth}
         \centering
         \includegraphics[width=\textwidth]{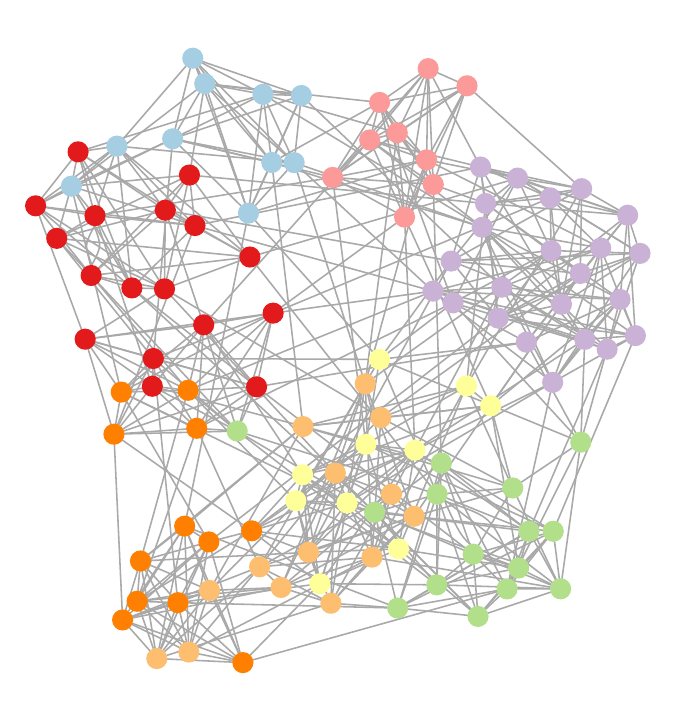}
         \caption{Third level (0.810)}
         \label{fig:FB3}
     \end{subfigure}
      \hfill
     \begin{subfigure}[h]{0.24\textwidth}
         \centering
         \includegraphics[width=\textwidth]{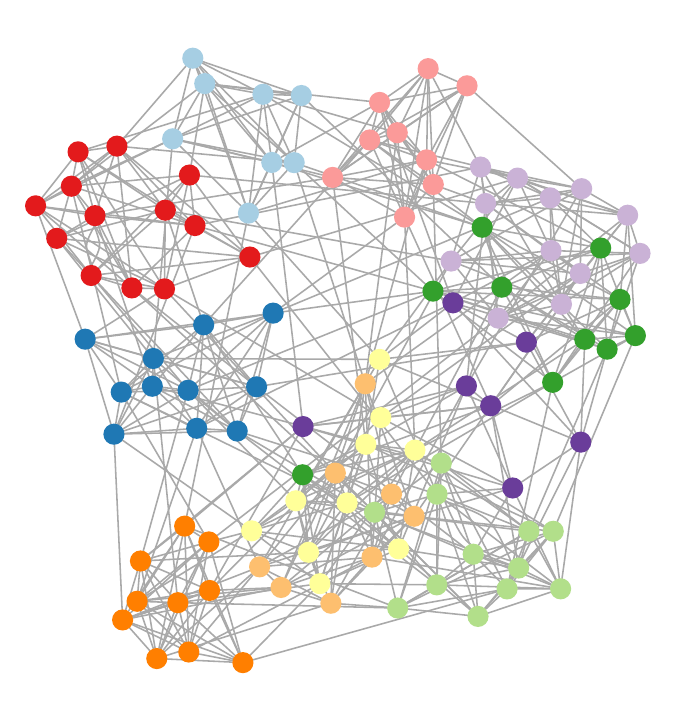}
         \caption{Ground truth}
         \label{fig:FBtruth}
     \end{subfigure}
        \caption{Spectral recursive bi-partition on the Football network, with the completeness scores in the parentheses}
        \label{fig:Football}
\end{figure}

~\\
\textbf{Football} is a network of American college football teams during the regular season Fall 2000. The vertices represent teams and edges represent regular season games between the two teams they connect. The teams are divided into ``conferences'' containing around 8 to 12 teams each, which can be treated as the true primitive communities; see Figure \ref{fig:FBtruth}. Again, we apply the recursive bi-partitioning algorithm based on the Fiedler vector. Here, we build a balanced hierarchy of depth three without resorting to any stopping rule. Figure \ref{fig:FB1} - \ref{fig:FB3} show the first, second and third level of an estimated hierarchy, with 2, 4, and 8 resulting clusters respectively. Each level has a high completeness score, suggesting that the estimated hierarchy is meaningful.

\begin{figure}[h]
     \centering
     \begin{subfigure}[h]{0.24\textwidth}
         \centering
         \includegraphics[width=\textwidth]{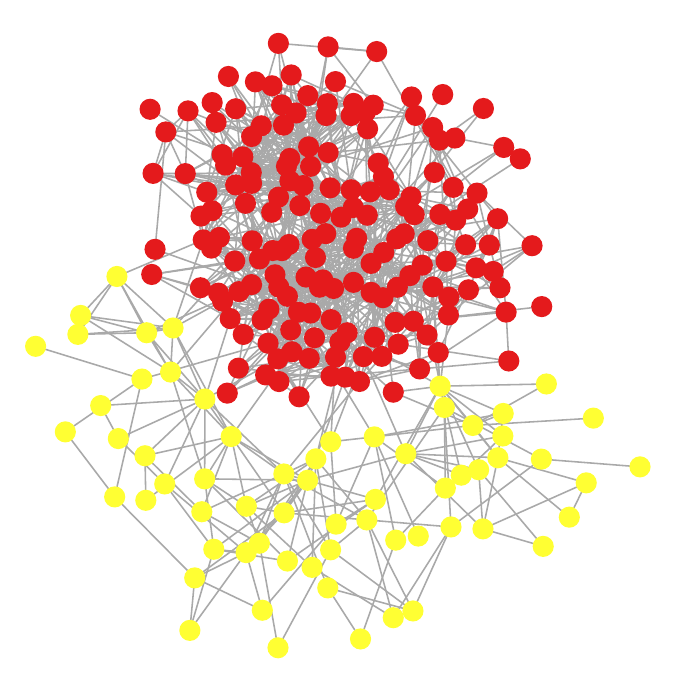}
         \caption{First split (0.939)}
         \label{fig:Ce1}
     \end{subfigure}
     \hfill
     \begin{subfigure}[h]{0.24\textwidth}
         \centering
         \includegraphics[width=\textwidth]{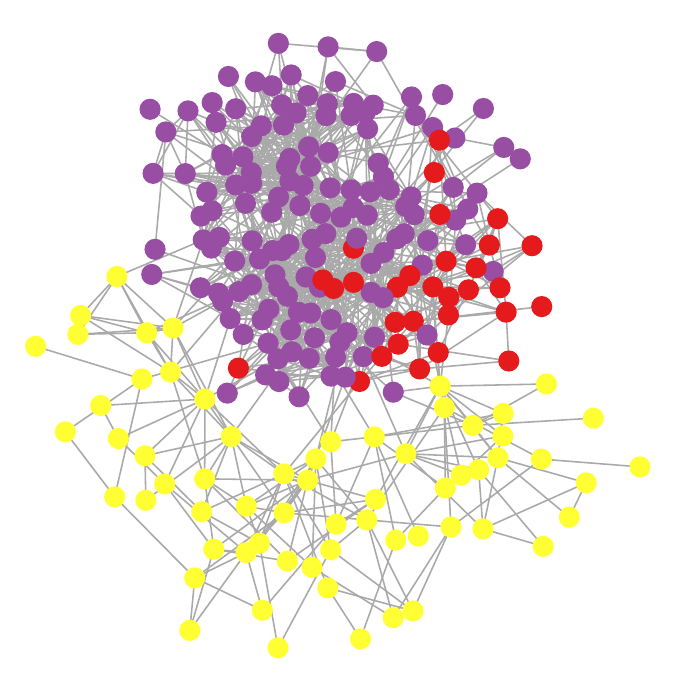}
         \caption{Second split (0.667)}
         \label{fig:Ce2}
     \end{subfigure}
      \hfill
     \begin{subfigure}[h]{0.24\textwidth}
         \centering
         \includegraphics[width=\textwidth]{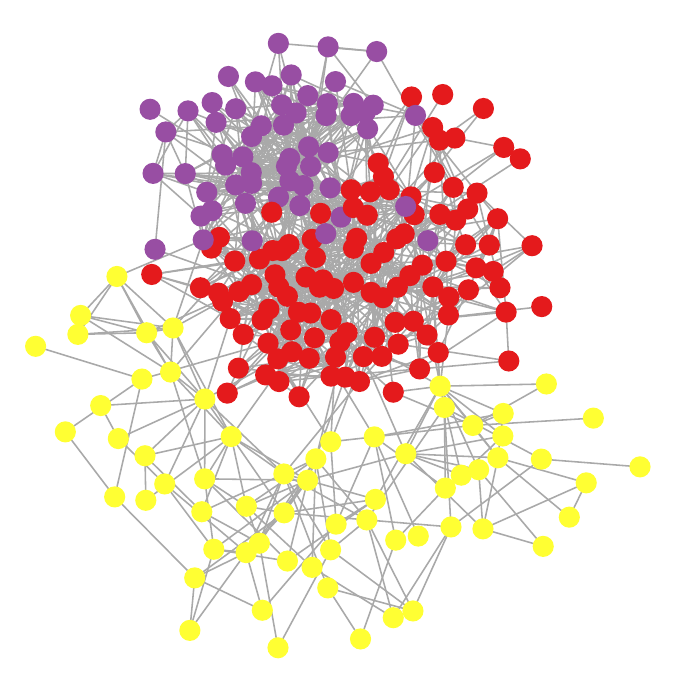}
         \caption{Second split (0.752)}
         \label{fig:Ce3}
     \end{subfigure}
      \hfill
     \begin{subfigure}[h]{0.24\textwidth}
         \centering
         \includegraphics[width=\textwidth]{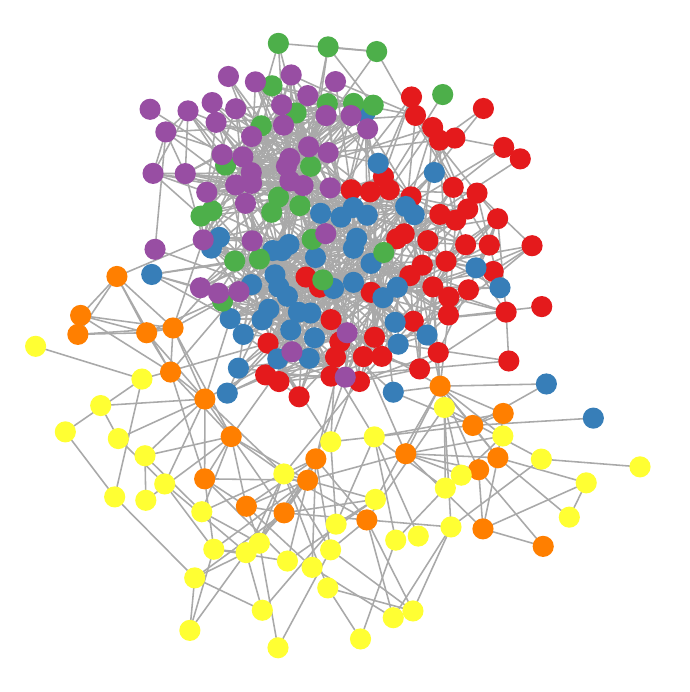}
         \caption{Ground truth}
         \label{fig:Cetruth}
     \end{subfigure}
    \caption{Spectral recursive bi-partition on the C. elegans network, with the completeness scores in the parentheses. (a) and (b) are based on the unnormalized graph Laplacian and (c) is based on the normalized graph Laplacian}
    \label{fig:Celegans}
\end{figure}

~\\
\textbf{C. elegans} is a neural network consisting of gap junctional synaptic connections in the posterior nervous system of a single adult male of \textit{Caenorhabditis elegans}, a primitive worm. The cells are grouped according to the modules and categories described in \cite{jarrell2012connectome}. Specifically, there are six types of cells: sensory neurons, interneurons, gender-shared neurons, command and motor neurons, gender-shared muscle cells, sex-specific muscle cells. We treat the cell types as the true primitive communities; see Figure \ref{fig:Cetruth}. Figure \ref{fig:Ce1} shows the first split given by the spectral bi-partition algorithm based on the Fiedler vector. It performs well in the first split as suggested by the completeness score 0.939. Interestingly, the two mega-communities correspond to neurons and muscle cells precisely. Figure \ref{fig:Ce2} shows the second split, which has a degraded performance.  As a comparison, we also show the second split given by the normalized graph Laplacian in Figure \ref{fig:Ce3}, which has a higher completeness score. They mainly differ in whether the blue nodes are merged with purple or red nodes.

\section{Proofs of main results}\label{sec:proofs}
% section 6 proofs

\subsection{Supporting lemmata}

\begin{lem}[Lemma 3.8 of \cite{lei2019unified}]
Let $\L$ denote the unnormalized graph Laplacian of a random unweighted graph with independent edges, and $\E(\L)=\L^*$. Then for any $r>0$, there exists an absolute constant $C(r)$ that only depends on $r$, such that, with probability at least $1-n^{-r}$,
\begin{equation*}
\nm{\mtx{L}-\mtx{L}^*}_{\op} \leq C(r) \sqrt{\lb\max_{1 \leq i \leq n} L^*_{ii} + \log n\rb \log n}.
\end{equation*}
\label{lem:laplacian_concentration}
\end{lem}

\begin{lem}[Lemma 10 of \cite{balakrishnan2011noise}]\label{lem:laplacian_monotonicity}
Let $\mtx{A}$ and $\td{\mtx{A}}$ be two adjacency matrices with unnormalized Laplacians $\mtx{L}$ and $\td{\mtx{L}}$. If $\td{\mtx{A}}_{ij}\ge \mtx{A}_{ij}$ for any pair $(i, j)$, then $\td{\mtx{L}} - \mtx{L}$ is positive semidefinite.
\end{lem}

% \begin{lem}[Weyl's inequality] Let $\mtx{A}, \mtx{A}^* \in \R^{n\times n}$ be symmetric matrices, with eigenvalues $\lambda_1 \geq\ldots\geq\lambda_n$ and $\lambda^*_1\geq \ldots\geq\lambda^*_n$ respectively. Then the following inequality holds for $j\in [n]$,
% \[ |\lambda_j-\lambda_j^*|\leq \nm{\mtx{A}-\mtx{A}^*}.\]
% \label{weyl}
% \end{lem}

\begin{lem}[Theorem 2 of \cite{yu2014useful}, a variant of Davis-Kahan $\sin \Theta$ Theorem] Let $\mtx{A}, \mtx{A}^* \in \R^{n\times n}$ be symmetric matrices. Fix positive integers $s,d$ and let $\mtx{U}=(\vct{u}_s(\A),\ldots,\vct{u}_{s+d-1}(\A))\in \R^{n\times d}$ and $\mtx{U}^*=(\vct{u}_s(\A^{*}),\ldots,\vct{u}_{s+d-1}(\A^{*}))\in \R^{n\times d}$. Further let $\Theta(\mtx{U}, \mtx{U}^{*})\in \R^{d\times d}$ denote the principal angle matrix between the subspaces spanned by $\mtx{U}$ and $\mtx{U}^{*}$. Then
\begin{align*}
\inf_{\substack{\mtx{O}\in \R^{d\times d}\\\mtx{O}^T\mtx{O} = I_{d}}}\|\mtx{U}\mtx{O} - \mtx{U}^{*}\|_{F} &= 2^{1/2}\nm{\sin \Theta(\mtx{U}, \mtx{U}^*)}_F 
\\
&\leq \frac{2^{3/2}d^{1/2}\nm{\mtx{A}-\mtx{A}^*}_{\op}}{\min\{\lambda_{s-1}(\A^{*})-\lambda_{s}(\A^{*}),\lambda_{s+d-1}(\A^{*})-\lambda_{s+d}(\A^{*})\}}.
\end{align*}
\label{davis_kahan}
\end{lem}

\subsection{Proof of Theorem \ref{thm:eigen}}
\label{sec:proof_eigen}
We first state a more general result that characterizes the entire eigenstructure of $\L^*$ under the general BTSBM.

% For any node $s$ with depth $|s| > i$ and binary representation $b_1 b_2 \ldots b_{|s|}$, let $s_{(i)}$ denote the $i$-th ancestor of $s$: 
% \begin{equation*}
%   \label{eq:si}
%   s_{(i)} \triangleq b_1 b_2 \ldots b_{|s|-i}, \quad \forall i = 0, 1, \ldots, |s| - 1, \mbox{ and } s_{(|s|)} = \emptyset,
% \end{equation*}
% $\bar{s}$ the sibling node:
% \begin{equation*}
%   \label{eq:bars}
%   \bar{s} \triangleq b_1 b_2 \ldots b_{|s|-1}(1 - b_{|s|}),
% \end{equation*}

\begin{lem}
\label{lem:eigen}
Denote
\[
g(s; \T) = \left\{
    \begin{array}{ll}
      1 & (s \mbox{ is an internal node})\\
      n_{s} - 1 & (s \mbox{ is a leaf node})
    \end{array}
\right..
\]
The eigenstructure of $\mtx{L}^*$ has the following properties:
\begin{enumerate}[(1)]
\item $\lambda_n^* = 0$, $\vec{u}_n^* = \frac{1}{\sqrt{n}}[1, 1, \ldots, 1]^\top = \one_{n} / \sqrt{n}$, and $\lambda_{n-1}^* > 0$.
\item For each node $s = b_1 b_2 \ldots b_{|s|}$, let $s_{(i)} = b_1 b_2 \ldots b_{|s| - i}$ and $s_{(|s|)} = \es$. Then
\begin{equation}
\label{eq:eigenvalues}
\lambda^*(s; \T) \triangleq n_{s}p_{s} + \sum_{i=1}^{|s|} (n_{s_{(i)}} - n_{s_{(i-1)}}) p_{s_{(i)}}, 
\end{equation}
is an eigenvalue of $\mtx{L}^*$ with multiplicity 
\[
\sum_{s': \lambda^*(s'; \T) = \lambda^*(s; \T)}g(s'; \T).
\]
\item The eigenspace corresponding to $\lambda^*(s; \T)$ is spanned by 
\[
\bigcup_{s': \lambda^*(s'; \T) = \lambda^*(s; \T)} \colspan\left(\mtx{U}(s'; \T)\right)
\]
where $\mtx{U}(s; \T)\in \R^{n\times g(s; \T)}$ such that
\begin{itemize}
\item if $s$ is an internal node,
\[\mtx{U}_{i}(s; \T) = \left\{
    \begin{array}{ll}
     \sqrt{n_{R(s)} / n_{L(s)} n_s} & i\in \G_{L(s)}\\
      -\sqrt{n_{L(s)} / n_{R(s)} n_s} & i\in \G_{R(s)}\\
      0 & \mbox{otherwise}
    \end{array}
\right.;\]
\item if $s$ is a leaf node, $\mtx{U}_{\G_{s}^{c}} (s; \T) = \mtx{0}_{(n-n_s) \times (n_s -1)}$ and $\mtx{U}_{\G_{s}}(s; \T) \in \mathbb{R}^{n_s \times (n_s-1)}$ is any orthogonal matrix with $\one_{n_{s}}^{\top}\mtx{U}_{\G_{s}}(s; \T) = \vct{0}^\top$.
\end{itemize}
\end{enumerate}
\end{lem}
Lemma \ref{lem:eigen} can be verified through simple algebra. By Lemma \ref{lem:eigen}, it is straightforward to prove Theorem \ref{thm:eigen}.

\begin{proof}[Proof of Theorem \ref{thm:eigen}]
  \begin{enumerate}[(1)]
  \item By Lemma \ref{lem:eigen}, $np_{\es}$ is an eigenvalue, corresponding to the root node and $g(\es; \T) = 1$ since it is an internal node. Under weak assortativity, for any node $s\in \T$,
\[\lambda^{*}(s; \T) > n_{s}p_{\es} + \sum_{i=1}^{|s|}(n_{s_{(i)}} - n_{s_{(i-1)}})p_{\es} = np_{\es}.\]
Therefore, $\lambda_{n-1}^{*} = np_{\es}$ with multiplicity $1$.
  \item By Lemma \ref{lem:eigen}, $n_{1}p_{1} + n_{0}p_{\es}$ and $n_{0}p_{0} + n_{1}p_{\es}$ are both eigenvalues corresponding to node $1$ and $0$, respectively. For all other nodes $s\in \T$, it is easy to show that $\lambda^{*}(s; \T) > n_{0}p_{0} + n_{1}p_{\es}$ if $s$ is a descendant of node $0$,  and $\lambda^{*}(s; \T) > n_{1}p_{1} + n_{0}p_{\es}$ if $s$ is a descendant of node $1$. Therefore, $\lambda_{n-2}^{*}$ must be their minimum.
  \item For each leaf node $s$, it is not hard to see that
\[\lambda^{*}(s; \T) = \sum_{j=1}^{n}p_{ij}, \quad \forall i \mbox{ in the community }s.\]
  By definition, 
\[\lambda^{*}(s; \T)\ge n\ubar{p}_{*}.\]
  The number of eigenvalues that are at least $n\ubar{p}_{*}$, accounting for multiplicity, is at least
\[\sum_{\mbox{{\scriptsize leaf node}~} s}g(s; \T) = \sum_{\mbox{{\scriptsize leaf node}~} s}(n_{s} - 1) = n - K.\]
\item Based on the previous part in this proof, if the eigenvalue $\lambda^*_{j} < n\ubar{p}_{*}$, then $\lambda^*_{j}$ must correspond to an internal node. Then the part (3) of Lemma \ref{lem:eigen} implies that
\[\|\vec{u}_j^{*}\|_{\infty} = \frac{1}{\sqrt{n}}\max\left\{\sqrt{\frac{n_{R(s)}}{n_{L(s)}}}, \sqrt{\frac{n_{L(s)}}{n_{R(s)}}}\right\}\le \sqrt{\frac{\xi}{n}}.\]
  \end{enumerate}
\end{proof}

\subsection{Proof of Lemma \ref{lem:laplacian_property}}
Let $\widetilde{\A}$ and $\widetilde{\A}'$ be defined as in Section \ref{subsec:ideas}, i.e.,
\[\widetilde{\A}_{ij} = \left\{\begin{array}{ll}
     \A_{ij} & (i, j\in \mathcal{G}_0, \text{ or }i,j\in\mathcal{G}_1)  \\
     p_{\es} & (\text{otherwise})
\end{array}\right., \quad \widetilde{\A}'_{ij} = \left\{\begin{array}{ll}
     \A_{ij}\mtx{B}_{ij} & (i, j\in \mathcal{G}_0, \text{ or }i,j\in\mathcal{G}_1)  \\
     p_{\es} & (\text{otherwise})
\end{array}\right.,\]
where $\{\mtx{B}_{ij}: i < j\}$ are independent Bernoulli random variables that are independent of $\A$ with $\E[\mtx{B}_{ij}] = p_{0} / \E[\A_{ij}]$ if $i,j\in \mathcal{G}_0$ and $\E[\mtx{B}_{ij}] = p_{1} / \E[\A_{ij}]$ if $i,j\in \mathcal{G}_1$. Then $\widetilde{\A}'$ is the adjacency matrix given by a general BTSBM with a 2-leaf tree and parameters $(p_{\es}, p_{0}, p_{1})$. Further let $\td{\L}$ and $\td{\L}'$ be their unnormalized graph Laplacians. Since $\widetilde{\A}_{ij} \ge \widetilde{\A}'_{ij}$ for all pairs $(i, j)$, by Lemma \ref{lem:laplacian_monotonicity},
\begin{equation}\label{eq:lambda_n-2}
\lambda_{n-2}(\td{\L})\ge \lambda_{n-2}(\td{\L}').
\end{equation}
Note that 
\[\max_{1\le i \le n}\E[\td{\L}_{ii}'] = \max\{n_0 p_0 + n_1 p_{\es}, n_1 p_1 + n_0 p_{\es}\}\le n_0 p_0 + n_1 p_1.\]
By Weyl's inequality and Lemma \ref{lem:laplacian_concentration}, with probability $1 - n^{-r}$,
\[\lambda_{n-2}(\td{\L}')\ge \lambda_{n-2}(\E[\td{\L}']) - \|\td{\L}' - \E[\td{\L}']\|\ge \lambda_{n-2}(\E[\td{\L}']) - C(r)\sqrt{(n_0 p_0 + n_1 p_1 + \log n)\log n}.\]
By Lemma \ref{lem:eigen},
\[\lambda_{n-2}(\E[\td{\L}']) = \min\{n_0 p_0 + n_1 p_{\es}, n_1 p_1 + n_0 p_{\es}\} = np_{\es} + \min\{n_0 (p_0 - p_{\es}), n_1 (p_1 - p_{\es})\}.\]
The condition \eqref{eq:cond_l2} implies
\begin{equation}\label{eq:cond_l2_implication}
n_0 p_0 + n_1 p_1  \ge C_{\ell_2}\sqrt{(n_0 p_0 + n_1 p_1)\log n}\Longrightarrow n_0 p_0 + n_1 p_1\ge C_{\ell_2}^2 \log n.
\end{equation}
When $C_{\ell_2} > \max\{1, 3C(r)\}$,
\begin{align}
   \lambda_{n-2}(\td{\L}')& \ge np_{\es} + \left\{C_{\ell_2} - C(r)\sqrt{1 + \frac{1}{C_{\ell_2}^2}}\right\}\sqrt{(n_0 p_0 + n_1 p_1)\log n}\nonumber\\
   & \ge np_{\es} + \frac{C_{\ell_2}}{2}\sqrt{(n_0 p_0 + n_1 p_1)\log n}\label{eq:lambda_n-2'}
\end{align}
By \eqref{eq:lambda_n-2}, 
\begin{equation}\label{eq:lambda_n-2_lower}
\lambda_{n-2}(\td{\L}) > np_{\es}, \quad \text{with probability }1 - n^{-r}.
\end{equation}
On the other hand, for any $i\in \mathcal{G}_0$, 
\[(\td{\L} \vec{u}_{n-1}^{*})_{i} = \sum_{j=1}^{n}\td{\L}_{ij}\vec{u}^{*}_{n-1, j} = \sqrt{\frac{n_1}{n_0 n}}\sum_{j\in \mathcal{G}_0}\td{\L}_{ij} - \sqrt{\frac{n_0}{n_1 n}}\sum_{j\in \mathcal{G}_1}\td{\L}_{ij}.\]
By definition, 
\[\sum_{j\in \mathcal{G}_0}\td{\L}_{ij} = \td{\L}_{ii} - \sum_{j\in \mathcal{G}_0\setminus \{i\}}\mtx{A}_{ij} = \sum_{j\in \mathcal{G}_1\cup \{i\}}\mtx{A}_{ij} = n_1 p_{\es} = -\sum_{j\in \mathcal{G}_1}\td{\L}_{ij}.\]
Thus, 
\[(\td{\L} \vec{u}^{*}_{n-1})_{i} = \lb\sqrt{\frac{n_1}{n_0 n}} + \sqrt{\frac{n_0}{n_1 n}}\rb n_1 p_{\es} = \sqrt{\frac{n_1}{n_0 n}}np_{\es} = (np_{\es})\vec{u}_{n-1, i}^{*}.\]
Similarly, for any $i\in \mathcal{G}_1$,
\[(\td{\L} \vec{u}^{*}_{n-1})_{i} = (np_{\es})\vec{u}^{*}_{n-1, i}.\]
As a result, $(np_{\es}, \vec{u}^{*}_{n-1})$ is an eigenpair of $\td{\L}$.  Since $\td{\L}$ is an unnormalized Laplacian, $0$ is always an eigenvalue. Therefore, on the event \eqref{eq:lambda_n-2_lower}, which occurs with probability at least $1 - n^{-r}$, $\lambda_{n-1}(\td{\L}) = np_{\es}$ with multiplicity $1$ and $\vec{u}_{n-1}(\td{\L}) = \vec{u}^{*}_{n-1}$.

\subsection{Proof of Theorem \ref{thm:l2_perturb}}
\label{sec:proof_l2_perturb}
Let $\mathcal{E}$ denote the event given by \eqref{eq:lambda_n-2'}. As shown in the proof of Lemma \ref{lem:laplacian_property}, 
\[\P(\mathcal{E})\ge 1 - n^{-r},\]
if $C_{\ell_2} > \max\{1, 3C(r)\}$. 
On the event $\mathcal{E}$, by \eqref{eq:lambda_n-2} and Lemma \ref{lem:laplacian_property},
\[\lambda_{n-1}(\td{\L}) = np_{\es},\quad \lambda_{n-2}(\td{\L}) - \lambda_{n-1}(\td{\L})\ge \frac{C_{\ell_2}}{2}\sqrt{(n_0 p_0 + n_1 p_1)\log n}, \quad \vec{u}_{n-1}(\td{\L}) = \vec{u}_{n-1}^{*}.\]
Fix any $\nu > 0$. Let 
\begin{equation*}
  \L_{\nu} = \L + \nu \J, \quad \td{\L}_{\nu} = \td{\L} + \nu \J, \quad \mbox{where }\J = \I - \frac{1}{n}\one_{n}\one_{n}^{T}.
\end{equation*}
Since both $\L$ and $\td{\L}$ are unnormalized graph Laplacians, $\lambda_{n}(\L) = \lambda_{n}(\td{\L}) = 0, \vec{u}_{n}(\L) = \vec{u}_{n}(\td{\L}) = \one_{n}/\sqrt{n}$. In addition, for any $\nu > 0$,
\[\vec{u}_{n-1}(\L_{\nu}) = \vec{u}_{n-1}(\L), \quad \vec{u}_{n-1}(\td{\L}_{\nu}) = \vec{u}_{n-1}(\td{\L}),\]
and for $j\ge 1$,
\[\lambda_{n-j}(\L_{\nu}) = \lambda_{n-j}(\L) + \nu, \quad \lambda_{n-j}(\td{\L}_{\nu}) = \lambda_{n-j}(\td{\L}) + \nu.\]
On the event $\mathcal{E}$,
\begin{equation}\label{eq:eigen_gap_tdL_full}
    \min\{\lambda_{n-2}(\td{\L}) - \lambda_{n-1}(\td{\L}), \lambda_{n-1}(\td{\L}) - \lambda_{n}(\td{\L})\}\ge \min\left\{\frac{C_{\ell_2}}{2}\sqrt{(n_0 p_0 + n_1 p_1)\log n}, np_{\es} + \nu\right\}.
\end{equation}
By Davis-Kahan Theorem (Lemma \ref{davis_kahan}),
\begin{align*}
  \|\vec{u}_{n-1}\sign(\vec{u}_{n-1}^T\vec{u}_{n-1}^{*}) - \vec{u}_{n-1}^{*}\|_{2}&\le \frac{2^{3/2}\|\L_{\nu} - \td{\L}_{\nu}\|}{\min\{\lambda_{n-2}(\td{\L}) - \lambda_{n-1}(\td{\L}), \lambda_{n-1}(\td{\L}) - \lambda_{n}(\td{\L})\}}\\
  & \le \frac{2^{3/2}\|\L - \td{\L}\|}{\min\left\{\frac{C_{\ell_2}}{2}\sqrt{(n_0 p_0 + n_1 p_1)\log n}, np_{\es} + \nu\right\}}
\end{align*}
By definition, $\L - \td{\L}$ can be formulated as $\L_{0} - \E[\L_{0}]$ where $\L_{0}$ is the unnormalized graph Laplacian for an adjacency matrix $\A_{0}$ with 
\[\A_{0, ij} = \left\{\begin{array}{ll}
     0 & (i, j\in \mathcal{G}_0, \text{ or }i,j\in\mathcal{G}_1)  \\
     \A_{ij} & (\text{otherwise})
\end{array}\right..\]
By Lemma \ref{lem:laplacian_concentration}, with probability at least $1 - n^{-r}$,
\[\|\L - \td{\L}\|\le C(r)\sqrt{(np_{\es} + \log n)\log n}.\]
Let $\nu = n\max\{p_1, p_0\}$. Then
\[ \|\vec{u}_{n-1}\sign(\vec{u}_{n-1}^T\vec{u}_{n-1}^{*}) - \vec{u}_{n-1}^{*}\|_{2}\le \frac{2^{5/2}C(r)}{C_{\ell_2}}\frac{\sqrt{(np_{\es} + \log n)\log n}}{\sqrt{(n_0 p_0 + n_1 p_1)\log n}}.\]
By weak assortativity and \eqref{eq:cond_l2_implication},
\[\sqrt{(np_{\es} + \log n)\log n}\le \sqrt{2(n_0 p_0 + n_1 p_1)\log n}.\]
Thus, with probability $1 - 2n^{-r}$,
\[\|\vec{u}_{n-1}\sign(\vec{u}_{n-1}^T\vec{u}_{n-1}^{*}) - \vec{u}_{n-1}^{*}\|_{2}\le \frac{8C(r)}{C_{\ell_2}}.\]
The proof is completed by setting $C_{\ell_2} = \max\{1, 8C(r)/ c\}$ and replacing $r$ with $r + 1$.

\subsection{Proof of Theorem \ref{thm:linf_perturb}}
\label{sec:proof_linf_perturb}
Without loss of generality we assume $\mtx{u}_{n-1}^{T}\mtx{u}_{n-1}^{*} \ge 0$. For each $j\in \{1, \ldots, n\}$ such that $\lambda_{n-j+1}^{*} < \lambda_{n-j}^{*}$, let $\mtx{U}_{j}^{*}\in \R^{n\times j}$ and $\mtx{U}_{j}\in \R^{n\times j}$ denote the eigenvector matrices $(\mtx{u}_{n-1}^{*}, \ldots, \mtx{u}_{n-j}^{*})$ and $(\mtx{u}_{n-1}, \ldots, \mtx{u}_{n-j})$, respectively. Define
\[\mtx{O}_{j} = \sign(\mtx{U}_{j}^{T}\mtx{U}_{j}^{*}),\]
where $\sign(\mtx{M})$ denotes the matrix sign. Specifically, if $U\Sigma V^{T}$ is the singular value decomposition of $\mtx{M}$, then $\sign(\mtx{M}) = UV^{T}$. Since $\mtx{O}_{j}\in \R^{j\times j}$ is an orthogonal matrix, we have
\[\sqrt{n}\|\mtx{U}_{j}\mtx{O}_{j} - \mtx{U}_{j}^{*}\|_{\tti} = \sqrt{n}\|\mtx{U}_{j} - \mtx{U}_{j}^{*}\mtx{O}_{j}^{T}\|_{\tti}.\]
Let $O_{j, 1, i}$ denote the entry of $\mtx{O}_{j}$ in the first row and $i$-th column and $\mtx{U}_{j, -1}^{*}$ denote the matrix $\mtx{U}_{j}^{*}$ with the first column $\vec{u}_{n-1}^{*}$ removed. Then
\begin{align*}
&\sqrt{n}\|\mtx{U}_{j} - \mtx{U}_{j}^{*}\mtx{O}_{j}^{T}\|_{\tti}\ge \sqrt{n}\bigg\|\vec{u}_{n-1} - \sum_{i=1}^{j}O_{j, 1, i}\vec{u}_{n-i}^{*}\bigg\|_{\infty}\\
& \quad \ge \sqrt{n}\|\vec{u}_{n-1} - O_{j,1,1}\vec{u}_{n-1}^{*}\|_{\infty} - \sqrt{n}\sqrt{\sum_{i\not = 1}O_{j,1,i}^{2}} \|\mtx{U}_{j, -1}^{*}\|_{\tti}.
\end{align*}
Furthermore, we know that
\begin{align*}
\sqrt{n}\|\vec{u}_{n-1} - O_{j,1,1}\vec{u}_{n-1}^{*}\|_{\infty} & = \sqrt{n}\|\vec{u}_{n-1} - \vec{u}_{n-1}^{*} + \vec{u}_{n-1}^{*} - O_{j,1,1}\vec{u}_{n-1}^{*}\|_{\infty}\\
& \ge \sqrt{n}\|\vec{u}_{n-1} - \vec{u}_{n-1}^{*}\|_{\infty} - \sqrt{n}|1 - O_{j,1,1}|\|\vec{u}_{n-1}^{*}\|_{\infty}\\
& = \sqrt{n}\|\vec{u}_{n-1} - \vec{u}_{n-1}^{*}\|_{\infty} - |1 - O_{j,1,1}|\max\left\{\sqrt{\frac{n_1}{n_0}}, \sqrt{\frac{n_0}{n_1}}\right\}\\
& \ge \sqrt{n}\|\vec{u}_{n-1} - \vec{u}_{n-1}^{*}\|_{\infty} - |1 - O_{j,1,1}|\sqrt{\xi}.
\end{align*}
The second last equality invokes Theorem \ref{thm:eigen}. Since $\mtx{O}_{j}$ is orthogonal,
\[\sum_{i=1}^{j}O_{j,1,i}^{2} = 1\Longrightarrow \sum_{i\not = 1}O_{j,1,i}^{2} = 1 - O_{j, 1, 1}^{2}\le 2(1 - O_{j, 1, 1}).\]
Also notice that $\|\mtx{U}_{j, -1}^{*}\|_{\tti}\le \|\mtx{U}_{j}^{*}\|_{\tti}$. As a result,
\begin{align}
  &\sqrt{n}\|\vec{u}_{n-1} - \vec{u}_{n-1}^{*}\|_{\infty}\nonumber\\
  \le &\sqrt{n}\|\mtx{U}_{j}\mtx{O}_{j} - \mtx{U}_{j}^{*}\|_{\tti} + \sqrt{2(1 - O_{j, 1, 1})}\lb\sqrt{n}\|\mtx{U}_{j}^{*}\|_{\tti}\rb+ |1 - O_{j,1,1}|\sqrt{\xi}.  \label{eq:sketch_key0}
\end{align}
To further simplify the second and the third terms, let $\mtx{H}_{j}= \mtx{U}_{j}^{T}\mtx{U}_{j}^{*}$ with singular value decomposition $\mtx{H}_{j} = \bar{\mtx{U}}_{j}(\cos \Theta(\mtx{U}_{j}, \mtx{U}_{j}^{*}))\bar{\mtx{V}}_{j}^{T}$, where $\cos\Theta(\mtx{U}_{j}, \mtx{U}_{j}^{*}) = \diag(\cos\theta_{j1}, \ldots, \cos \theta_{jj})$ and $\theta_{ji}$'s are the principal angles between $\vec{u}_{n-i}$ and $\vec{u}^{*}_{n-i}$. By definition, $\mtx{O}_{j}=  \bar{\mtx{U}}_{j}\bar{\mtx{V}}_{j}^{T}$. As a result,
\[\|\mtx{H}_{j} - \mtx{O}_{j}\| = \|\bar{\mtx{U}}_{j}(\I - \cos \Theta_{j})\bar{\mtx{V}}_{j}^{T}\|\le \|\I - \cos \Theta(\mtx{U}_{j}, \mtx{U}_{j}^{*})\|.\]
For any $\theta\le \pi / 2$, 
\[1 - \cos \theta \le 1 - \cos^{2}\theta = \sin^{2}\theta.\]
Therefore,
\[\|\mtx{H}_{j} - \mtx{O}_{j}\|\le \|\sin \Theta(\mtx{U}_{j}, \mtx{U}_{j}^{*}))\|^{2}\]
% By Davis-Kahan Theorem,
% \[\|H_{j} - O_{j}\|\le \frac{\|E\|^{2}}{\Delta_{j}^{*2}}.\]
On the other hand,
\[|1 - H_{j, 1, 1}| = |1 - \vec{u}_{n-1}^{T}\vec{u}_{n-1}^{*}| = |(\vec{u}_{n-1} - \vec{u}_{n-1}^{*})^{T}\vec{u}_{n-1}^{*}|\le \|\vec{u}_{n-1} - \vec{u}_{n-1}^{*}\|_{2}.\]
As a consequence, 
\begin{align*}
  &|1 - O_{j, 1, 1}|\le |1 - H_{j, 1, 1}| + |H_{j, 1, 1} - O_{j, 1, 1}|\\
\le &\|\vec{u}_{n-1} - \vec{u}_{n-1}^{*}\|_{2} + \|\mtx{H}_{j} - \mtx{O}_{j}\| \le \|\vec{u}_{n-1} - \vec{u}_{n-1}^{*}\|_{2} + \|\sin \Theta(\mtx{U}_{j}, \mtx{U}_{j}^{*})\|^{2}.
\end{align*}
This together with \eqref{eq:sketch_key0} imply
\begin{equation}
\label{eq:key_ineq}
\begin{split}
   &\sqrt{n}\|\vec{u}_{n-1} - \vec{u}_{n-1}^{*}\|_{\infty}\\
   \le & \sqrt{n}\|\mtx{U}_{j}\mtx{O}_{j} - \mtx{U}_{j}^{*}\|_{\tti} + \lb\sqrt{2\|\vec{u}_{n-1} - \vec{u}_{n-1}^{*}\|_{2}} + \sqrt{2}\|\sin \Theta(\mtx{U}_{j}, \mtx{U}_{j}^{*})\|\rb\lb\sqrt{n}\|\mtx{U}_{j}^{*}\|_{\tti}\rb  \\
  & + \sqrt{\xi}\lb\|\vec{u}_{n-1} - \vec{u}_{n-1}^{*}\|_{2} + \|\sin \Theta(\mtx{U}_{j}, \mtx{U}_{j}^{*})\|^{2}\rb.
\end{split}
\end{equation} 
Let $\td{K}$ denote the number of eigenvalues that are strictly smaller than $n\ubar{p}_{*}$. By Theorem \ref{thm:eigen}, $\td{K}\le K$. Since $\xi = O(1)$, $K = O(1)$. Then for any $j\le \td{K}$, by part (4) of Theorem \ref{thm:eigen}, 
\[\|\mtx{U}_{j}^{*}\|_{\tti}\lesssim\frac{1}{\sqrt{n}}.\]
Thus, to prove $\sqrt{n}\|\vec{u}_{n-1} - \vec{u}_{n-1}^{*}\|_{\infty} < \min\{\sqrt{n_0/n_1},\sqrt{n_1/n_0}\}$, it remains to prove
\begin{equation}
  \label{eq:sketch_goal}
   \|\vec{u}_{n-1} - \vec{u}_{n-1}^{*}\|_{2}\le c, \quad \|\sin \Theta(\mtx{U}_{j}, \mtx{U}_{j}^{*})\|\le c,\quad \sqrt{n}\|\mtx{U}_{j}\mtx{O}_{j} - \mtx{U}_{j}^{*}\|_{\tti} \le c,
\end{equation}
for a sufficiently small constant $c$ that only depends on $\xi$  for some $2\le j\le \td{K}$ with high probability. The first bound $\|\vec{u}_{n-1} - \vec{u}_{n-1}^{*}\|_{2}\le c$ has been proved in Theorem \ref{thm:l2_perturb} if $C_{\ell_{\infty}}\ge C_{\ell_{2}}$. We will show the other two bounds in the following subsections.

\subsubsection{Choice of $j$ via the pigeonhole principle}\label{subsubsec:choice_j}
% We start by some notation and a few simple facts that will be used repeatedly. Let
% \[p^{*} = \max_{ij}p_{ij}, \quad \bar{p}^{*} = \max_{i}\frac{1}{n}\sum_{j=1}^{n}p_{ij}
% , \quad \ubar{p}_{*} = \min_{i}\frac{1}{n}\sum_{j=1}^{n}p_{ij}, \quad \bar{p}_{2}^{*} = \max_{i}\sqrt{\frac{1}{n}\sum_{j=1}^{n}p_{ij}^{2}}.\]
% Then
% \begin{equation}\label{eq:barppstar}
% p^{*}\ge \bar{p}_{2}^{*}\ge\bar{p}^{*}\ge \frac{n_{s}}{n}p^{*}\Longrightarrow p^{*}, \bar{p}_{2}^{*} = O(\bar{p}^{*}).
% \end{equation}
By definition of $\td{K}$, $np_{\es} = \lambda_{n-1}^{*}\le \cdots \lambda_{n-\td{K}}^{*} < n\ubar{p}_{*} \le \lambda^{*}_{n-\td{K}-1}$. Let 
\begin{equation}\label{eq:tdk}
\delta_{j}^{*} = \min\{n\ubar{p}_{*}, \lambda_{n-j-1}^{*}\} - \lambda_{n-j}^{*}, \quad \td{j} = \argmax_{j\le \td{K}}\delta_{j}^{*}.
\end{equation}
By Theorem \ref{thm:eigen}, $\td{K}\le K$ and thus
\begin{equation}
  \label{eq:deltaJ}
  \delta_{\td{j}}^{*} \ge \frac{1}{\td{K}}\sum_{j=1}^{\td{K}}\delta_{j}^{*} = \frac{n(\ubar{p}_{*} - p_{\es})}{\td{K}}.
\end{equation}
Throughout the rest of the proof we will fix $j = \td{j}$ and depress the subscript $j$ when no confusion can arise. This option guarantees a sufficiently large eigengap so that the off-the-shelf technical tools can be applied directly to obtain meaningful perturbation bounds. 

For notational convenience, denote by $\mtx{\Lambda}$ the diagonal matrix of the $\td{K}$ smallest eigenvalues, i.e.,
\[\mtx{\Lambda} = \diag(\lambda_{n-1}, \ldots, \lambda_{n - \td{j} }), \quad \mtx{\Lambda}^{*} = \diag(\lambda_{n-1}^{*}, \ldots, \lambda_{n - \td{j}}^{*}).\]
We write $\mtx{U}, \mtx{O}$ and $\mtx{U}^{*}$ for $\mtx{U}_{\td{j}}, \mtx{O}_{\td{j}}$ and $\mtx{U}^{*}_{\td{j}}$ throughout the rest of the section. Similar to the proof of Theorem \ref{thm:l2_perturb}, let 
\begin{equation}\label{eq:Lnu}
  \L_{\nu} = \L + \nu \J, \quad \L_{\nu}^{*} = \L^{*} + \nu \J, \quad \mbox{where }\J = \I - \frac{1}{n}\one_{n}\one_{n}^{T}.
\end{equation}
Since $\one_{n}^{T}\mtx{U} = \one_{n}^{T}\mtx{U}^{*} = 0$,
\[\L_{\nu}\mtx{U} = \mtx{U}(\mtx{\Lambda} + \nu \I), \quad \L_{\nu}^{*}\mtx{U}^{*} = \mtx{U}^{*}(\mtx{\Lambda}^{*} + \nu \I).\]
Throughout we take 
\begin{equation}
  \label{eq:nu}
  \nu = n\bar{p}^{*}.
\end{equation}
% Moreover, by assumption \textbf{C}2, 
% \begin{equation}
%   \label{eq:barp}
%   (n\bar{p}^{*})^{4} \ge (n(\ubar{p}_{*} - p_{\es})^{4})\ge C_1(n\bar{p}^{*})^{3}\log n \Longrightarrow n\bar{p}^{*}\gtrsim \log n.
% \end{equation}
% Finally, since $\td{k}\le K \lesssim 1$, 
% \begin{equation}
%   \label{eq:mnormU}
%   \sqrt{n}\mnorm{\mtx{U}^{*}}\lesssim 1.
% \end{equation}

\subsubsection{Bounding $\|\sin \Theta(U, U^{*})\|$}
Applying Lemma \ref{lem:laplacian_concentration}, we have 
\begin{equation}
  \label{eq:Eop2}
  \|\mtx{\L} - \mtx{\L}^{*}\|\le C(r)\lb\sqrt{n\bar{p}^{*}\log n} + \log n\rb\quad \mbox{with probability }1 - n^{-r}
\end{equation}
Since $\delta_{\td{j}}^{*}$, the eigengap defined in \eqref{eq:deltaJ}, is invariant to $\nu$, by Davis-Kahan Theorem (Lemma \ref{davis_kahan}) and \eqref{eq:nu},
\begin{align*}
\|\sin \Theta(\mtx{U}, \mtx{U}^{*})\|&\le \frac{2\|\mtx{L}_{\nu} - \mtx{\L}^{*}_{\nu}\|}{\min\left\{\delta_{\td{j}}^{*}, np_{\es} + \nu\right\}} = \frac{2\|\mtx{L} - \mtx{\L}^{*}\|}{\delta_{\td{j}}^{*}} \le 2KC(r)\frac{\lb\sqrt{n\bar{p}^{*}\log n} + \log n\rb}{n(\ubar{p}_{*} - p_{\es})}
\end{align*}
By the condition \eqref{eq:cond_pbar_pubar},
\begin{equation}
  \label{eq:barp}
  (n\bar{p}^{*})^{4} \ge (n(\ubar{p}_{*} - p_{\es}))^{4}\ge C_{\ell_{\infty}}(n\bar{p}^{*})^{3}\log n \Longrightarrow n\bar{p}^{*}\ge C_{\ell_{\infty}}\log n.
\end{equation}
Thus, with probability $1 - n^{-r}$,
\begin{align*}
\|\sin \Theta(\mtx{U}, \mtx{U}^{*})\| &\le 2KC(r)\frac{\lb\sqrt{n\bar{p}^{*}\log n} + \log n\rb}{C_{\ell_{\infty}}^{1/4}(n\bar{p}^{*})^{3/4}(\log n)^{1/4}}
\\
&\le 2KC(r)\frac{1 + C_{\ell_{\infty}}^{-1/2}}{C_{\ell_{\infty}^{1/2}}} = 2KC(r)\lb C_{\ell_{\infty}}^{-1/2} + C_{\ell_{\infty}}^{-1}\rb.
\end{align*}
Therefore, when $C_{\ell_{\infty}} \ge \max\{1, 4KC(r)/c\}$,
\[\|\sin \Theta(\mtx{U}, \mtx{U}^{*})\|\le c,\]
with probability $1 - n^{-r}$.

\subsubsection{Bounding $\sqrt{n}\|UO - U^{*}\|_{\tti}$}
We will apply Theorem 2.6 of \cite{lei2019unified} on $\L_{\nu}$ and $\L_{\nu}^{*}$, defined in \eqref{eq:Lnu} with $\nu$ defined in \eqref{eq:nu}. 
%\footnote{Indeed we are using a weaker version of Theorem 2.6 discussed in Remark 2.3 of \cite{lei2019unified} that bounds $\mnorm{U\sign(U^{T}U^{*}) - U^{*}}$ instead of $d_{\tti}(U, U^{*})$. This result replaces the effective condition number $\bar{\kappa}^{*}$ by the condition number $\kappa^{*}$. It is proved in Step I-IV in the proof of Theorem 2.5 in Appendix B.} 
To be self-contained, we state the theorem in the supplement together with all necessary definitions. We can prove the following result.

\begin{lem}
\label{lem:multiscale_linfty}
Let $\mtx{U}, \mtx{O}, \mtx{U}^{*}$ be defined as in Section \ref{subsubsec:choice_j}. Fix any constant $r, c > 0$. Under the same setting as in Theorem \ref{thm:linf_perturb}, if $C_{\ell_{\infty}}$ is a sufficiently large constant that only depends on $r$ and $c$, with probability at least $1 - (10K + 1)n^{-r}$, 
\[\sqrt{n}\mnorm{\mtx{U}\mtx{O} - \mtx{U}^{*}}\le c.\]
\end{lem}
The proof is lengthy and hence relegated to the supplement.

\section{Conclusion and Discussion}\label{sec:discuss}
% section 5 discussion
In this paper, we present a novel analysis of an unnormalized graph Laplacian-based recursive spectral clustering algorithm for sparse networks. Under a broad class of hierarchical network models, we show that the proposed algorithm is effective in both community detection and hierarchy estimation. Both weak and strong consistencies for mega-communities are established based on novel $\ell_2$ and $\ell_{\infty}$ perturbation bounds of the Fiedler vector. Compared to earlier works on hierarchical and non-hierarchical community detection, our result substantially relaxes the constraints on connection probabilities, degree heterogeneity, and the hierarchical structure to handle sparse and multiscale networks with an unbalanced hierarchy. 

% Our main theoretical contribution is that the result applies to sparse networks and gets rid of the constraint that the connection probabilities are at the same scale, which is required in previous works, by taking full advantage of the inherent structure of the graph Laplacian matrix.

One limitation of our model is that the hierarchy is restricted to be a binary tree. Theoretically, a binary tree cannot encode, for example, a $3$-block SBM with equal between-community connection probabilities due to the indistinguishability of the three primitive communities. However, the strict homogeneity as above is arguably a rare corner case in practice. When the between-community connection probabilities are mutually distinct, there exists a meaningful binary hierarchy in this case. Indeed, if $B_{12} = \max\{B_{12}, B_{13}, B_{23}\}$, the two mega-communities on the first level of the binary tree can be defined as $\{1, 2\}$ and $\{3\}$. We can verify that all between-mega-community connection probabilities are smaller than all within-mega-community connection probabilities in this case. Although it is beyond a general BTSBM, some of our results can be possibly extended to this case; see Remark \ref{rem:weak_general} for instance. It would be interesting to investigate what hierarchical structures can be equivalently formulated as a binary one.

Our main results can also be interpreted as the adaptivity of graph Laplacian based spectral clustering to inhomogeneous connection probabilities. The general BTSBM studied in this work is a special case of inhomogeneous block models in which within-community connection probabilities are greater than $p$ while between-community connection probabilities are smaller than $q<p$. Note that it is known that convex optimization approach would be adaptive to such inhomogeneous model \citep{moitra2016robust}, which is usually computationally expensive and sometimes involves sensitive tuning parameters. It would be interesting to see whether certain spectral method also has such adaptivity. In particular, we are interested in studying such adaptivity under either multiscale regime or the very sparse regime where the connection probabilities are on the order of $\omega(1/n)$.

Another important question is how degree heterogeneity degrades the performance of a clustering algorithm. Degree heterogeneity is ubiquitous in real-world networks while most theoretical works restrict the degrees to be on the same order. This work takes one step in mitigating the gap for a specific spectral clustering algorithm. As discussed in Remark \ref{rem:degree_var}, the degree variation condition \eqref{eq:cond_pbar_pubar} in Theorem \ref{thm:linf_perturb} appears to be a mathematical artifact, but to relax it requires novel techniques beyond our proof strategies summarized in Section \ref{subsec:ideas}. We leave this as an open problem and look for affirmation or negation of our conjecture.

% Generally, it is worth deriving similar results on the robustness to multiscale networks for other clustering algorithms. 

In terms of the technical proofs, as alluded to in Section \ref{sec:simulation}, the $\ell_{\infty}$ perturbation bound is only sufficient yet not necessary for the strong consistency of mega-communities recovery. Previous works \citep[e.g.][]{abbe2017entrywise, lei2019unified, deng2020strong} suggest that the sign consistency is strictly weaker than $\ell_{\infty}$ consistency for eigenvectors under certain special SBMs or BTSBMs. It is mathematically intriguing to explore how those advanced techniques can be adapted to the more heterogeneous BTSBMs.

\section*{Acknowledgment}
X. Li and X. Lou acknowledge support from the NSF via the Career Award DMS-1848575. X. Lou also acknowledges support from the NSF via the Grant CCF-1934568.

%\stitle{Supplemental materials to ``Consistency of Spectral Clustering on Hierarchical Stochastic Block Models''}
%\sdescription{We give in the supplement %\cite{lei2020supplement} proofs to Lemma %\ref{lem:multiscale_linfty} and an extension of Theorem %\ref{thm:l2_perturb} as mentioned in Remark %\ref{rem:weak_general}.}
%\end{supplement}

\bibliographystyle{apalike} 
\bibliography{ref}

\begin{thebibliography}{}

\bibitem[Abbe, 2017]{abbe2017community}
Abbe, E. (2017).
\newblock Community detection and stochastic block models: recent developments.
\newblock {\em The Journal of Machine Learning Research}, 18(1):6446--6531.

\bibitem[Abbe et~al., 2015]{abbe2015exact}
Abbe, E., Bandeira, A.~S., and Hall, G. (2015).
\newblock Exact recovery in the stochastic block model.
\newblock {\em IEEE Transactions on Information Theory}, 62(1):471--487.

\bibitem[Abbe et~al., 2017]{abbe2017entrywise}
Abbe, E., Fan, J., Wang, K., and Zhong, Y. (2017).
\newblock Entrywise eigenvector analysis of random matrices with low expected
  rank.
\newblock {\em arXiv preprint arXiv:1709.09565}.

\bibitem[Adamic and Glance, 2005]{adamic2005political}
Adamic, L.~A. and Glance, N. (2005).
\newblock The political blogosphere and the 2004 us election: divided they
  blog.
\newblock In {\em Proceedings of the 3rd international workshop on Link
  discovery}, pages 36--43.

\bibitem[Balakrishnan et~al., 2011]{balakrishnan2011noise}
Balakrishnan, S., Xu, M., Krishnamurthy, A., and Singh, A. (2011).
\newblock Noise thresholds for spectral clustering.
\newblock In {\em Advances in Neural Information Processing Systems}, pages
  954--962.

\bibitem[Cape et~al., 2019]{cape2019signal}
Cape, J., Tang, M., and Priebe, C.~E. (2019).
\newblock Signal-plus-noise matrix models: eigenvector deviations and
  fluctuations.
\newblock {\em Biometrika}, 106(1):243--250.

\bibitem[Clauset et~al., 2008]{clauset2008hierarchical}
Clauset, A., Moore, C., and Newman, M.~E. (2008).
\newblock Hierarchical structure and the prediction of missing links in
  networks.
\newblock {\em Nature}, 453(7191):98.

\bibitem[Damle and Sun, 2020]{damle2020uniform}
Damle, A. and Sun, Y. (2020).
\newblock Uniform bounds for invariant subspace perturbations.
\newblock {\em SIAM Journal on Matrix Analysis and Applications},
  41(3):1208--1236.

\bibitem[Dasgupta et~al., 2006]{dasgupta2006spectral}
Dasgupta, A., Hopcroft, J., Kannan, R., and Mitra, P. (2006).
\newblock Spectral clustering by recursive partitioning.
\newblock In {\em European Symposium on Algorithms}, pages 256--267. Springer.

\bibitem[Dasgupta, 2016]{dasgupta2016cost}
Dasgupta, S. (2016).
\newblock A cost function for similarity-based hierarchical clustering.
\newblock In {\em Proceedings of the forty-eighth annual ACM symposium on
  Theory of Computing}, pages 118--127.

\bibitem[Deng et~al., 2020]{deng2020strong}
Deng, S., Ling, S., and Strohmer, T. (2020).
\newblock Strong consistency, graph laplacians, and the stochastic block model.
\newblock {\em arXiv preprint arXiv:2004.09780}.

\bibitem[Eldridge et~al., 2017]{eldridge2017unperturbed}
Eldridge, J., Belkin, M., and Wang, Y. (2017).
\newblock Unperturbed: spectral analysis beyond davis-kahan.
\newblock {\em arXiv preprint arXiv:1706.06516}.

\bibitem[Fiedler, 1975]{fiedler1975property}
Fiedler, M. (1975).
\newblock A property of eigenvectors of nonnegative symmetric matrices and its
  application to graph theory.
\newblock {\em Czechoslovak Mathematical Journal}, 25(4):619--633.

\bibitem[Girvan and Newman, 2002]{girvan2002community}
Girvan, M. and Newman, M.~E. (2002).
\newblock Community structure in social and biological networks.
\newblock {\em Proceedings of the national academy of sciences},
  99(12):7821--7826.

\bibitem[Holland et~al., 1983]{holland1983stochastic}
Holland, P.~W., Laskey, K.~B., and Leinhardt, S. (1983).
\newblock Stochastic blockmodels: First steps.
\newblock {\em Social networks}, 5(2):109--137.

\bibitem[Jarrell et~al., 2012]{jarrell2012connectome}
Jarrell, T.~A., Wang, Y., Bloniarz, A.~E., Brittin, C.~A., Xu, M., Thomson,
  J.~N., Albertson, D.~G., Hall, D.~H., and Emmons, S.~W. (2012).
\newblock The connectome of a decision-making neural network.
\newblock {\em Science}, 337(6093):437--444.

\bibitem[Jin, 2015]{jin2015fast}
Jin, J. (2015).
\newblock Fast community detection by score.
\newblock {\em The Annals of Statistics}, 43(1):57--89.

\bibitem[Krzakala et~al., 2013]{krzakala2013spectral}
Krzakala, F., Moore, C., Mossel, E., Neeman, J., Sly, A., Zdeborov{\'a}, L.,
  and Zhang, P. (2013).
\newblock Spectral redemption in clustering sparse networks.
\newblock {\em Proceedings of the National Academy of Sciences},
  110(52):20935--20940.

\bibitem[Le et~al., 2017]{le2017concentration}
Le, C.~M., Levina, E., and Vershynin, R. (2017).
\newblock Concentration and regularization of random graphs.
\newblock {\em Random Structures \& Algorithms}, 51(3):538--561.

\bibitem[Lei and Rinaldo, 2015]{lei2015consistency}
Lei, J. and Rinaldo, A. (2015).
\newblock Consistency of spectral clustering in stochastic block models.
\newblock {\em The Annals of Statistics}, 43(1):215--237.

\bibitem[Lei, 2019]{lei2019unified}
Lei, L. (2019).
\newblock Unified $\ell_{2\rightarrow\infty}$ eigenspace perturbation theory
  for symmetric random matrices.
\newblock {\em arXiv preprint arXiv:1909.04798}.

\bibitem[Li et~al., 2018]{li2018hierarchical}
Li, T., Lei, L., Bhattacharyya, S., Sarkar, P., Bickel, P.~J., and Levina, E.
  (2018).
\newblock Hierarchical community detection by recursive partitioning.
\newblock {\em arXiv preprint arXiv:1810.01509}.

\bibitem[Lusseau et~al., 2003]{lusseau2003bottlenose}
Lusseau, D., Schneider, K., Boisseau, O.~J., Haase, P., Slooten, E., and
  Dawson, S.~M. (2003).
\newblock The bottlenose dolphin community of doubtful sound features a large
  proportion of long-lasting associations.
\newblock {\em Behavioral Ecology and Sociobiology}, 54(4):396--405.

\bibitem[Lyzinski et~al., 2016]{lyzinski2016community}
Lyzinski, V., Tang, M., Athreya, A., Park, Y., and Priebe, C.~E. (2016).
\newblock Community detection and classification in hierarchical stochastic
  blockmodels.
\newblock {\em IEEE Transactions on Network Science and Engineering},
  4(1):13--26.

\bibitem[Mao et~al., 2017]{mao2017estimating}
Mao, X., Sarkar, P., and Chakrabarti, D. (2017).
\newblock Estimating mixed memberships with sharp eigenvector deviations.
\newblock {\em arXiv preprint arXiv:1709.00407}.

\bibitem[McSherry, 2001]{mcsherry2001spectral}
McSherry, F. (2001).
\newblock Spectral partitioning of random graphs.
\newblock In {\em Proceedings 42nd IEEE Symposium on Foundations of Computer
  Science}, pages 529--537. IEEE.

\bibitem[Moitra et~al., 2016]{moitra2016robust}
Moitra, A., Perry, W., and Wein, A.~S. (2016).
\newblock How robust are reconstruction thresholds for community detection?
\newblock In {\em Proceedings of the forty-eighth annual ACM symposium on
  Theory of Computing}, pages 828--841.

\bibitem[Nepusz et~al., 2008]{nepusz2008fuzzy}
Nepusz, T., Petr{\'o}czi, A., N{\'e}gyessy, L., and Bazs{\'o}, F. (2008).
\newblock Fuzzy communities and the concept of bridgeness in complex networks.
\newblock {\em Physical Review E}, 77(1):016107.

\bibitem[Rohe et~al., 2011]{rohe2011spectral}
Rohe, K., Chatterjee, S., and Yu, B. (2011).
\newblock Spectral clustering and the high-dimensional stochastic blockmodel.
\newblock {\em The Annals of Statistics}, 39(4):1878--1915.

\bibitem[Rosenberg and Hirschberg, 2007]{rosenberg2007v}
Rosenberg, A. and Hirschberg, J. (2007).
\newblock V-measure: A conditional entropy-based external cluster evaluation
  measure.
\newblock In {\em Proceedings of the 2007 joint conference on empirical methods
  in natural language processing and computational natural language learning
  (EMNLP-CoNLL)}, pages 410--420.

\bibitem[Yu et~al., 2014]{yu2014useful}
Yu, Y., Wang, T., and Samworth, R.~J. (2014).
\newblock A useful variant of the davis--kahan theorem for statisticians.
\newblock {\em Biometrika}, 102(2):315--323.

\bibitem[Zachary, 1977]{zachary1977information}
Zachary, W.~W. (1977).
\newblock An information flow model for conflict and fission in small groups.
\newblock {\em Journal of anthropological research}, 33(4):452--473.

\bibitem[Zhao et~al., 2012]{zhao2012consistency}
Zhao, Y., Levina, E., and Zhu, J. (2012).
\newblock Consistency of community detection in networks under degree-corrected
  stochastic block models.
\newblock {\em The Annals of Statistics}, 40(4):2266--2292.

\end{thebibliography}


\begin{thebibliography}{11}
\bibitem{lei2019unified_S} Lei, L. (2019). Unified $\ell_{2\rightarrow\infty}$  Eigenspace Perturbation Theory for Symmetric Random Matrices, \emph{arXiv preprint arXiv:1909.04798}.
\end{thebibliography}

\newpage
\pagebreak
\begin{center}
\textbf{\large Supplemental materials to ``Consistency of Spectral Clustering on Hierarchical Stochastic Block Models''}
\end{center}
%%%%%%%%%% Merge with supplemental materials %%%%%%%%%%
%%%%%%%%%% Prefix a "S" to all equations, figures, tables and reset the counter %%%%%%%%%%
\setcounter{equation}{0}
\setcounter{figure}{0}
\setcounter{table}{0}
\setcounter{section}{0}
\setcounter{page}{1}
\makeatletter
\renewcommand{\theequation}{S\arabic{equation}}
\renewcommand{\thefigure}{S\arabic{figure}}
\renewcommand{\thesection}{S\arabic{section}}
\renewcommand{\bibnumfmt}[1]{[S#1]}
\renewcommand{\citenumfont}[1]{S#1}

\section{$\ell_2$ Perturbation Theory for Fiedler Vector Under Unrestricted Heterogeneity Within Mega-Communities}\label{app:l2}

\begin{thm}
\label{thm:l2_general}
Let $(\G_0, \G_1)$ denote a partition of $\{1, \ldots, n\}$, with $n_0 = |\G_0|$ and $n_1 = |\G_1|$, and $\mtx{A}\in\{0,1\}^{n\times n}$ the adjacency matrix such that
\[\E[\A_{ij}] \left\{\begin{array}{ll}
    = p_{\es} & (i\in \mathcal{G}_0, j\in \mathcal{G}_1) \\
    \ge p_0 & (i,j\in \mathcal{G}_0)\\
    \ge p_1 & (i,j\in\mathcal{G}_1)
\end{array}\right.,\]
where $p_{\es} < \min(p_0,p_1)$. Further let $\vec{u}_{n-1}$ and $\vec{u}_{n-1}^{*}$ denote the Fiedler vector of $\A$ and $\E[\A]$, respectively. Then, for any fixed $r > 0$, there exists a constant $B_{\ell_2}(r)$ that only depends on $r$ such that, with probability at least $1 - 2n^{-r}$,
\[
\|\vec{u}_{n-1}\sign(\vec{u}_{n-1}^{T}\vec{u}_{n-1}^{*}) - \vec{u}_{n-1}^{*}\|_2 \le B_{\ell_2}(r)\frac{\sqrt{(np_{\es} + \log n)\log n}}{\min\{n_{0}(p_{0} - p_{\es}), n_{1}(p_{1} - p_{\es})\}}.
\]
\end{thm}

\begin{proof}
Note that the proofs for Lemma \ref{lem:laplacian_property} and Theorem \ref{thm:l2_perturb} do not rely on the modelling assumptions within the mega-communities, except that $\mtx{P}_{ij} \ge p_0$ if $i, j\in \G_0$ and $\mtx{P}_{ij}\ge p_1$ if $i, j\in \G_1$. Therefore, the proofs carry over to Theorem \ref{thm:l2_general}.
\end{proof}

\section{$\ell_{\tti}$ Perturbation Theory for Unnormalized Laplacians}

\subsection{A Generic $\ell_{\tti}$ Perturbation Bound}
\label{app:generic}
In this subsection, we rephrase the weaker version Theorem 2.6 of \cite{lei2019unified_S}, discussed in their Remark 2.3, by only keeping the parts that are relevant to our purpose. Throughout this section, we consider two generic symmetric real matrices $\mtx{G}$ and $\mtx{G}^{*}$ with 
\begin{equation}\label{eq:E}
  \mtx{E} = \mtx{G} - \mtx{G}^{*}.
\end{equation}
 Let $\lambda_{1}\ge \lambda_{2}\ge \ldots \ge \lambda_{n}$ and $\lambda_{1}^{*}\ge \lambda_{2}^{*}\ge \ldots \ge \lambda_{n}^{*}$ be the eigenvalues of $\mtx{G}$ and $\mtx{G}^{*}$, respectively. Given positive integers $s$ and $r$, let
\begin{equation}\label{eq:LambdaLambda*}
\mtx{\Lambda} = \diag(\lambda_{s + 1}, \lambda_{s + 2}, \ldots, \lambda_{s + r}), \quad \mtx{\Lambda}^{*} = \diag(\lambda_{s + 1}^{*}, \lambda_{s + 2}^{*}, \ldots, \lambda_{s + r}^{*}).
\end{equation}
Let $\mtx{U}, \mtx{U}^{*}\in \R^{n\times r}$ be a matrix of eigenvectors such that 
\begin{equation}\label{eq:UU*}
  \mtx{G}\mtx{U} = \mtx{U}\mtx{\Lambda}, \quad \mtx{G}^{*}\mtx{U}^{*} = \mtx{U}^{*}\mtx{\Lambda}^{*}.
\end{equation}
To state the generic bound, we define the following quantities. 
\begin{itemize}
\item \emph{Modified perturbation matrix} $\td{\mtx{E}}$:
\[\td{\mtx{E}} = \mtx{G} - \mtx{\Sigma} - (\mtx{G}^{*} - \mtx{\Sigma}^{*})\]
where
\[\mtx{\Sigma} = \diag(\mtx{G}), \quad \mtx{\Sigma}^{*} = \diag(\mtx{G}^{*}).\]
\item \emph{Condition number} $\kappa^{*}$:
  \begin{equation}
    \label{eq:barkappa}
    \lambda_{\max}^{*} = \lambda_{\max}(\mtx{\Lambda}^{*}), \quad \lambda_{\min}^{*} = \lambda_{\min}(\mtx{\Lambda}^{*}), \quad \kappa^{*} = \lambda_{\max}^{*} / \lambda_{\min}^{*}.
  \end{equation}
\item \emph{Effective eigengap} $\Delta^{*}$: 
  \begin{equation}
    \label{eq:gap}
    \Delta^{*} \triangleq \min\{\mathrm{sep}^{*}, \lambda_{\min}^{*}\},
  \end{equation}
where $\mathrm{sep}^{*} = \min\{\lambda_{s}^{*} - \lambda_{s+1}^{*}, \lambda_{s+r}^{*} - \lambda_{s+r+1}^{*}\}$ and $\lambda_{0}^{*} = \infty, \lambda_{n+1}^{*} = -\infty$.
\end{itemize}

The assumptions for the generic bound are stated below.
\begin{enumerate}[\textbf{A}1]\label{enumi:assumptions}
\item For any $\delta \in (0, 1)$, 
\[\frac{\min_{j\in [s + 1, s + r]}|\Lambda_{jj}^{*}|}{\min_{j\in [s + 1, s + r], k\in [n]}|\Lambda_{jj}^{*} - \Sigma_{kk}|}\le \Theta(\delta),\]
with probability at least $1 - \delta$ for some deterministic function $\Theta(\delta) > 0$.
\item For any $\delta \in (0, 1)$, there exists a random matrix $\mtx{G}^{(k)}\in \R^{n\times n}$ such that
v\[d_{TV}\lb\P_{(\td{\mtx{E}}_{k}, \mtx{G}^{(k)})}, \P_{\td{\mtx{E}}_{k}}\times \P_{\mtx{G}^{(k)}}\rb\le \delta / n.\]
where $d_{TV}$ denotes the total variation distance and it holds simultaneously for all $k$ and all contiguous subsets $S\subset [r]$ that
\[\|\mtx{G}^{(k)} - \mtx{G}\|_{\op}\le L_{1}(\delta), \quad \frac{\|(\mtx{G}^{(k)} - \mtx{G})\mtx{U}\|_{\op}}{\lambda_{\min}^{*}}\le \lb \kappa(\mtx{\Lambda}^{*})L_{2}(\delta) + L_{3}(\delta)\rb \mnorm{\mtx{U}},\]
with probability at least $1 - \delta$ for some deterministic functions $L_{1}(\delta), L_{2}(\delta), L_{3}(\delta)$.
\item There exists deterministic functions $\lambda_{-}(\delta), E_{+}(\delta), \Ep_{\infty}(\delta)$, such that for any $\delta \in (0, 1)$, the following event holds with probability at least $1 - \delta$:
\[\maxnorm{\mtx{\Lambda} - \mtx{\Lambda}^{*}}\le \lambda_{-}(\delta), \quad \|\mtx{E} \mtx{U}^{*}\|_{\op}\le E_{+}(\delta), \quad \mnorm{\td{\mtx{E}}}\le \Ep_{\infty}(\delta).\]
\item There exists deterministic functions $\bp_{\infty}(\delta), \bp_{2}(\delta) > 0$, such that for any $\delta \in (0, 1)$, $k\in [n]$, and fixed matrix $\mtx{W}\in \R^{n\times \td{j}}$,
  \begin{align*}
    \|\td{\mtx{E}}_{k}^{T}\mtx{W}\|_{2} &\le \bp_{\infty}(\delta)\mnorm{\mtx{W}} + \bp_{2}(\delta)\|\mtx{W}\|_{\op}, \,\, \mbox{with probability at least } 1 - \delta / n.
  \end{align*}
\item $\Delta^{*}\ge 4\lb\Theta(\delta)\sigmap(\delta) + L_{1}(\delta) + \lambda_{-}(\delta) + E_{+}(\delta)\rb$ where
  \begin{equation}
    \label{eq:etap}
\etap(\delta)  = \Ep_{\infty}(\delta) + \bp_{\infty}(\delta) + \bp_{2}(\delta), \quad \sigmap(\delta) = \{\kappa^{*}L_{2}(\delta) + L_{3}(\delta) + 1\}\etap + E_{+}(\delta).
  \end{equation}
\end{enumerate}

\begin{thm}[Theorem 2.6 of \cite{lei2019unified_S}, Remark 2.3]\label{thm:generic}
  Under assumptions \textbf{A}1 - \textbf{A}5, 
\begin{align}
 \mnorm{\mtx{U}\sign(\mtx{U}^{T}\mtx{U}^{*})- \mtx{U}^{*}}&\le C\bigg\{\frac{\Theta(\delta)}{\lambda_{\min}^{*}}\mnorm{\mtx{E} \mtx{U}^{*}} + \lb \frac{E_{+}(\delta)^{2}}{(\Delta^{*})^{2}} + \frac{\Theta(\delta)\sigmap(\delta)}{\Delta^{*}}\rb \mnorm{\mtx{U}^{*}}\nonumber\\
& \qquad\quad + \frac{\Theta(\delta)(\bp_{2}(\delta) + \mnorm{\mtx{G}^{*} - \mtx{\Sigma}^{*}})E_{+}(\delta)}{\lambda_{\min}^{*}\Delta^{*}}\bigg\},\nonumber
\end{align}
with probability at least $1 - B(r)\delta$, where $C$ is a universal constant (that can be chosen as $136$) and 
\begin{equation}
  \label{eq:Br}
  B(r) = 10\min\{r, 1 + \log_{2}\kappa^{*}\}.
\end{equation}
\end{thm}

\subsection{Proof of Lemma \ref{lem:multiscale_linfty}}
Let $\mtx{G} = \L + \nu \J$ and $\mtx{G}^{*} = \L^{*} + \nu \J$. Now we verify each of assumptions \textbf{A}1 - \textbf{A}5. We add the subscript $\nu$ into all quantities defined in Supplement \ref{app:generic} to highlight their dependence on $\nu$, including $\Lambda_{\nu}, \Sigma_{\nu}, \Theta_{\nu}$ (or $\Lambda_{\nu}^{*}, \Sigma_{\nu}^{*}, \Theta_{\nu}^{*}$), and $\mtx{E}_{\nu}, \td{\mtx{E}}_{\nu}, L_{1, \nu}, L_{2, \nu}, L_{3, \nu}, \td{E}_{\infty, \nu}, E_{+, \nu}, \lambda_{-, \nu}, \td{b}_{\infty, \nu}, \td{b}_{2, \nu}, \kappa_{\nu}^{*}, \Delta_{\nu}^{*}, \td{\eta}_{\nu}, \td{\sigma}_{\nu}$. We remove the subscript $\nu$ when $\nu = 0$. Moreover, we let
\begin{equation}
  \label{eq:Mdelta}
  M(\delta) = \sqrt{n\bar{p}^{*}\log (n / \delta)} + \log (n / \delta), \quad R(\delta) = \log (n / \delta) + \td{j}.
\end{equation}
By definition \eqref{eq:tdk}, $\td{j} \le \td{K}\le K = O(1)$. In each of the following steps, $\delta$ is always set to be $n^{-r}$. Unless otherwise specified, $a \gtrsim b$ ($a \lesssim b$) iff $a \ge Cb$ ($a\le Cb$) for some constant $C$ that only depends on $r$ and $\xi$. To apply Theorem \ref{thm:generic} in Supplement \ref{app:generic}, we need to verify Assumptions \textbf{A}1 - \textbf{A}5. 

~\\
\noindent\textbf{Checking Assumption A1}: We recall Lemma 3.12 of \cite{lei2019unified_S}, rephrased for our purpose. 
\begin{lem}\label{lem:A1}
Let $\Theta(\delta)$ be defined in assumption \textbf{A}1 in Supplement \ref{app:generic}. Further let
  \begin{equation*}
    \Theta^{*} = \frac{\min_{j\in [n-\td{j}, n-1]}|\Lambda^{*}_{jj}|}{\min_{j\in [n-\td{j}, n-1], k\in [n]}|\Lambda^{*}_{jj} - \Sigma^{*}_{kk}|}.
  \end{equation*}
Then $\Theta(\delta)\le 5\Theta^{*}$ if  
\[\min_{j\in [n-\td{j}, n-1], k\in [n]}|\Lambda^{*}_{jj} - \Sigma^{*}_{kk}|\ge 5M(\delta).\]
\end{lem}
In this case, $\Lambda_{jj}^{*} < n\ubar{p}_{*}$ for all $j\in [n-\td{j}, n-1]$ and $\Sigma_{kk}^{*} \ge n\ubar{p}_{*}$ for all $k\in [n]$. Thus,
\[\min_{j\in [n-\td{j}, n-1], k\in [n]}|\Lambda^{*}_{jj} - \Sigma^{*}_{kk}| = \min_{k\in [n]}\Sigma_{kk}^{*} - \max_{j\in [n - \td{j}, n-1]}\Lambda_{jj}^{*} = n\ubar{p}_{*} - \lambda^{*}_{n-\td{j}}.\]
By definition of $\td{j}$ and equation \eqref{eq:deltaJ},
\[n\ubar{p}_{*} - \lambda^{*}_{n-\td{j}}\ge \frac{n(\ubar{p}_{*} - p_{\es})}{K}.\]
By definition,
\[\Lambda^{*}_{\nu, jj} = \Lambda^{*}_{jj} + \nu, \quad \Sigma^{*}_{\nu, kk} = \Sigma^{*}_{kk} + \frac{n-1}{n}\nu.\]
Thus, 
\[\min_{j\in [n-\td{j}, n-1], k\in [n]}|\Lambda^{*}_{\nu, jj} - L^{*}_{\nu, kk}|\ge \frac{n(\ubar{p}_{*} - p_{\es})}{K} - \frac{\nu}{n}.\]
By the condition \eqref{eq:cond_pbar_pubar} and \eqref{eq:barp},
\[\frac{n(\ubar{p}_{*} - p_{\es})}{K}\ge \frac{C_{\ell_{\infty}}^{1/4}(n\bar{p}^{*})^{3/4}(\log n)^{1/4}}{K}\ge \frac{C_{\ell_{\infty}}}{K}\log n.\]
If $C_{\ell_{\infty}}\ge 3\xi\ge 3K$,
\[\frac{n(\ubar{p}_{*} - p_{\es})}{K}\ge 2 \ge 2\bar{p}^{*}.\]
By the definition \eqref{eq:nu} of $\nu$, 
\[\frac{\nu}{n} = \bar{p}^{*}\le \frac{n(\ubar{p}_{*} - p_{\es})}{2K}.\]
As a result,
\[\min_{j\in [n-\td{j}, n-1], k\in [n]}|\Lambda^{*}_{\nu, jj} - L^{*}_{\nu, kk}|\ge \frac{n(\ubar{p}_{*} - p_{\es})}{2K}.\]
On the other hand, by \eqref{eq:barp}, 
\[M(n^{-r}) \le (r+1)\lb 1 + C_{\ell_{\infty}}^{-1}\rb\sqrt{n\bar{p}^{*}\log n}.\]
By the condition \eqref{eq:cond_pbar_pubar} and \eqref{eq:barp} again, if $C_{\ell_{\infty}}^{1/2}\ge 20\xi(r+1)\ge 20K(r+1)$,
\[\frac{n(\ubar{p}_{*} - p_{\es})}{2K}\ge \frac{C_{\ell_{\infty}}^{1/4}(n\bar{p}^{*})^{3/4}(\log n)^{1/4}}{2K}\ge  \frac{C_{\ell_{\infty}}^{1/2}\sqrt{(n\bar{p}^{*})\log n}}{2K}\ge 5M(n^{-r}).\]
Therefore, 
\begin{equation}
  \label{eq:Theta}
  \Theta_{\nu}(n^{-r})\le 5 \Theta_{\nu}^{*} \le \frac{5(np_{\es} + \nu)}{n(\ubar{p}_{*} - p_{\es}) / 2K}\lesssim \frac{n\bar{p}^{*}}{n(\ubar{p}_{*} - p_{\es})}.
\end{equation}

~\\
\textbf{Checking Assumption A2}: We recall Lemma 3.10 of \cite{lei2019unified_S}.
\begin{lem}\label{lem:A2}
There exists $\mtx{G}^{(1)}, \ldots, \mtx{G}^{(n)}$ satisfying \textbf{A}2 for $\mtx{G} = \L$ with
\[L_{1}(\delta) \lesssim M(\delta), \quad L_{2}(\delta) = 1, \quad L_{3}(\delta) \lesssim \frac{n\bar{p}^{*} + \log (n / \delta)}{\lambda_{\min}(\mtx{\Lambda}^{*})},\]
where $\lesssim$ only hides absolute constants and $L_{1}(\delta), L_{2}(\delta) , L_{3}(\delta)$ are defined in Assumption \textbf{A}2 in Supplement \ref{app:generic}.
\end{lem}

\noindent In this case, let 
\[\mtx{G}^{(k)}_{\nu} = \mtx{G}^{(k)} + \nu \J.\]
Then it is easy to see that 
\[\mtx{G}^{(k)}_{\nu} - \L_{\nu} = \mtx{G}^{(k)} - \L.\]
Therefore, Lemma \ref{lem:A2} holds for any $\nu > 0$. Let $\delta = n^{-r}$, we have
\begin{equation}
  \label{eq:L1L2L3}
  L_{1, \nu}(n^{-r})\lesssim \sqrt{n\bar{p}^{*}\log n}, \quad L_{2, \nu}(n^{-r})\lesssim 1, \quad L_{3, \nu}(n^{-r})\lesssim \frac{n\bar{p}^{*} + \log n}{np_{\es} + \nu}\lesssim \frac{n\bar{p}^{*} + \log n}{n\bar{p}^{*}}\lesssim 1,
\end{equation}
where the last inequality uses \eqref{eq:barp}.

~\\
\textbf{Checking Assumption A3}: We recall Lemma 3.8 of \cite{lei2019unified_S}. 
\begin{lem}\label{lem:A3}
  Assumption \textbf{A}3 is satisfied for $\mtx{G} = \L$ with
\[\Ep_{\infty}(\delta)\lesssim  \sqrt{n\bar{p}^{*}} + \sqrt{\log (n / \delta)},\quad E_{+}(\delta), \lambda_{-}(\delta)\lesssim M(\delta),\]
where $\lesssim$ only hides absolute constants and $\Ep_{\infty}(\delta), E_{+}(\delta), \lambda_{-}(\delta)$ are defined in Assumption \textbf{A}3 in Supplement \ref{app:generic}.
\end{lem}
Since $\mtx{\Lambda}_{\nu} = \mtx{\Lambda} + \nu \I$ and $\mtx{\Lambda}_{\nu}^{*} = \mtx{\Lambda}^{*} + \nu \I$, $\mtx{\Lambda}_{\nu} - \mtx{\Lambda}_{\nu}^{*} = \mtx{\Lambda} - \mtx{\Lambda}^{*}$ is invariant to $\nu$. Similarly, $\mtx{E}_{\nu} = \mtx{E}$ and $\td{\mtx{E}}_{\nu} = \td{\mtx{E}}$. Thus Lemma \ref{lem:A3} holds for any $\nu > 0$. By \eqref{eq:barp},
\begin{equation}
  \label{eq:Eop}
  \Ep_{\infty, \nu}(n^{-r})\lesssim \sqrt{n\bar{p}^{*}}, \quad E_{+, \nu}(n^{-r}), \lambda_{-, \nu}(n^{-r})\lesssim \sqrt{n\bar{p}^{*}\log n}.
\end{equation}

~\\
\textbf{Checking Assumption A4}: We recall Lemma 3.7 of \cite{lei2019unified_S}.
\begin{lem}\label{lem:A4}
Assumption \textbf{A}4 is satisfied for $\mtx{G} = \L$ with
\[\bp_{\infty}(\delta) \lesssim \frac{R(\delta)}{\alpha \log R(\delta)}, \quad \bp_{2}(\delta) \lesssim \frac{\sqrt{p^{*}}R(\delta)^{(1 + \alpha) / 2}}{\alpha \log R(\delta)},\]
where $\lesssim$ only hides absolute constants and $\bp_{\infty}(\delta), \bp_{2}(\delta)$ are defined in Assumption \textbf{A}4 in Supplement \ref{app:generic}.
\end{lem}
As with Assumption \textbf{A}3, $\Ep$ is invariant to $\nu$. Thus Lemma \ref{lem:A4} holds for any $\nu > 0$. Let $\alpha = 1 / \log R(\delta)$. Since $\td{j} \le K = O(1)$,
\begin{equation}
  \label{eq:binftyb2}
  \bp_{\infty, \nu}(n^{-r})\lesssim R(n^{-r})\lesssim \log n, \quad \bp_{2, \nu}(\delta)\lesssim \sqrt{R(\delta)p^{*}}\lesssim \sqrt{(\log n)p^{*}} \lesssim \sqrt{(\log n)\bar{p}^{*}}.
\end{equation}
where the last inequality uses the fact that $\bar{p}^{*}\ge (\max n_{s})p^{*} / n$.

~\\
\textbf{Checking Assumption A5}: We first refer the readers to Supplement \ref{app:generic} for the definitions of $\kappa^{*}, \Delta^{*}, \etap(\delta)$ and $\sigmap(\delta)$. By definition, 
\begin{equation}
  \label{eq:kappa}
  \kappa_{\nu}^{*} = \frac{\lambda_{\max}(\Lambda_{\nu}^{*})}{\lambda_{\min}(\Lambda_{\nu}^{*})} = \frac{\lambda_{n-\td{j}}^{*} + \nu}{np_\es+\nu}\lesssim 1,
\end{equation}
and 
\begin{equation}
  \label{eq:Delta}
  \Delta_{\nu}^{*} = \min\{\nu, n\ubar{p}_{*} - \lambda_{n-\td{j}+1}\}\ge \min \left\{\nu, \frac{n(\ubar{p}_{*} - p_{\es})}{K}\right\} =  \frac{n(\ubar{p}_{*} - p_{\es})}{K}.
\end{equation}
By definition of $\etap$, \eqref{eq:Eop} and \eqref{eq:binftyb2},
\begin{equation*}
  \etap_{\nu}(n^{-r}) \lesssim \sqrt{n\bar{p}^{*}} + \log n.
\end{equation*}
By \eqref{eq:L1L2L3}, \eqref{eq:Eop} and \eqref{eq:barp},
\begin{equation}
  \label{eq:sigma}
  \sigmap_{\nu}(n^{-r})\lesssim \etap(n^{-r}) + \sqrt{n\bar{p}^{*}\log n} \lesssim \sqrt{n\bar{p}^{*}} + \log n + \sqrt{n\bar{p}^{*}\log n}\lesssim \sqrt{n\bar{p}^{*}\log n}.
\end{equation}
By \eqref{eq:Theta}, \eqref{eq:L1L2L3} and \eqref{eq:Eop},
\begin{align*}
  &\Theta_{\nu}(n^{-r})\sigmap_{\nu}(n^{-r}) + L_{1, \nu}(n^{-r}) + \lambda_{-, \nu}(n^{-r}) + E_{+, \nu}(n^{-r})\\
\lesssim &\frac{n\bar{p}^{*}}{n(\ubar{p}_{*} - p_{\es})}\sqrt{n\bar{p}^{*}\log n} + \sqrt{n\bar{p}^{*}\log n}\lesssim \frac{n\bar{p}^{*}}{n(\ubar{p}_{*} - p_{\es})}\sqrt{n\bar{p}^{*}\log n}.
\end{align*}
% Equivalently, there exists a constant $C_{1}$ that only depends on $q, K$ and $\{s\in \T: n_{L(s)} / n_{R(s)}\}$,
% \[\Theta_{\nu}(n^{-r})\sigmap_{\nu}(n^{-r}) + L_{1, \nu}(n^{-r}) + \lambda_{-, \nu}(n^{-r}) + E_{+, \nu}(n^{-r})\le C_{1}\frac{n\bar{p}^{*}}{n(\ubar{p}_{*} - p_{\es})}\sqrt{n\bar{p}^{*}\log n}\]
By \eqref{eq:Delta} and the condition \eqref{eq:cond_pbar_pubar}, if $C_{\ell_{\infty}}$ is sufficiently large,
\begin{equation}\label{eq:A5}
\Delta_{\nu}^{*}\ge 4\lb\Theta_{\nu}(n^{-r})\sigmap_{\nu}(n^{-r}) + L_{1, \nu}(n^{-r}) + \lambda_{-, \nu}(n^{-r}) + E_{+, \nu}(n^{-r})\rb,
\end{equation}
Thus, Assumption \textbf{A}5 is satisfied.

~\\
\textbf{Final Result}:\\
In the previous five steps, we show that Assumption \textbf{A}1 - \textbf{A}5 are satisfied under the condition \eqref{eq:cond_pbar_pubar}, if $C_{\ell_{\infty}}$ is sufficiently large. By Theorem \ref{thm:generic}, with probability $1 - B(\td{j})n^{-r}$, 
\begin{align}
 &\mnorm{\mtx{U}\sign(\mtx{U}^{T}\mtx{U}^{*})- \mtx{U}^{*}}\\
 &\lesssim \frac{\Theta_{\nu}(n^{-r})}{\lambda_{\min}(\mtx{\Lambda}_{\nu}^{*})}\mnorm{\mtx{E}_{\nu} \mtx{U}^{*}} + \lb \frac{E_{+, \nu}^{2}(n^{-r})}{(\Delta_{\nu}^{*})^{2}} + \frac{\Theta_{\nu}(n^{-r})\sigmap_{\nu}(n^{-r})}{\Delta_{\nu}^{*}}\rb \mnorm{\mtx{U}^{*}}\nonumber\\
& \qquad\quad + \frac{\Theta_{\nu}(n^{-r})E_{+, \nu}(n^{-r})}{\Delta_{\nu}^{*}}\frac{\bp_{2, \nu}(n^{-r}) + \mnorm{\L_{\nu}^{*} - \mtx{\Sigma}_{\nu}^{*}}}{\lambda_{\min}(\mtx{\Lambda}_{\nu}^{*})}.\nonumber
\end{align}
To bound $\mnorm{\mtx{E}_{\nu}\mtx{U}^{*}}$, we recall Lemma 3.9 of \cite{lei2019unified_S}. 
\begin{lem}\label{lem:EU}
Let $M(\delta)$ and $R(\delta)$ be defined in \eqref{eq:Mdelta}. Then with probability $1 - \delta$,
\[\mnorm{\mtx{E} \mtx{U}^{*}}\lesssim (M(\delta) + \td{j})\mnorm{\mtx{U}^{*}} + \sqrt{R(\delta)p^{*}}.\]
\end{lem}
Note that $\mtx{E}_{\nu} = \mtx{E}$. When $\delta = n^{-r}$, by \eqref{eq:barp}, 
\begin{equation}
  \label{eq:EU}
  \mnorm{\mtx{E}_{\nu} \mtx{U}^{*}}\lesssim \sqrt{n\bar{p}^{*}\log n}\mnorm{\mtx{U}^{*}} + \sqrt{(\log n)\bar{p}^{*}}\lesssim \sqrt{n\bar{p}^{*}\log n}\mnorm{\mtx{U}^{*}},
\end{equation}
where the last line uses the fact that $\sqrt{n}\mnorm{\mtx{U}^{*}}\ge 1$. 

Now we derive bounds for other terms. By \eqref{eq:Theta} and the definition \eqref{eq:nu} of $\nu$,
\begin{equation}\label{eq:term1}
\frac{\Theta_{\nu}(n^{-r})}{\lambda_{\min}(\Lambda_{\nu}^{*})}\le \frac{\Theta_{\nu}(n^{-r})}{\nu}\lesssim \frac{1}{n(\ubar{p}_{*} - p_{\es})}  
\end{equation}
Furthermore, by \eqref{eq:Theta}, \eqref{eq:Eop}, \eqref{eq:Delta} and \eqref{eq:sigma}, 
\begin{align}
  &\frac{E_{+, \nu}^{2}(n^{-r})}{(\Delta_{\nu}^{*})^{2}} + \frac{\Theta_{\nu}(n^{-r})\sigmap_{\nu}(n^{-r})}{\Delta_{\nu}^{*}}\lesssim \frac{n\bar{p}^{*}\log n}{(n(\ubar{p}_{*} - p_{\es}))^{2}} + \frac{n\bar{p}^{*}\sqrt{n\bar{p}^{*}\log n}}{(n(\ubar{p}_{*} - p_{\es}))^{2}}\lesssim \frac{n\bar{p}^{*}\sqrt{n\bar{p}^{*}\log n}}{(n(\ubar{p}_{*} - p_{\es}))^{2}},\label{eq:term2}
\end{align}
where the last inequality uses \eqref{eq:barp}. For the third term, note that
% \[\L_{\nu, ii}^{*} - \Sigma_{\nu, ii}^{*} = 0, \quad \L_{\nu, ij}^{*} - \Sigma_{\nu, ij}^{*} = \L_{ij}^{*} - \Sigma_{ij}^{*} - \frac{\nu}{n}.\]
\[\L_{\nu}^{*} - \mtx{\Sigma}_{\nu}^{*} = \L^{*} - \mtx{\Sigma}^{*} +  \nu\J - \frac{n - 1}{n}\nu \I.\]
Thus,
\begin{align}
  \sqrt{n}\mnorm{\L_{\nu}^{*} - \mtx{\Sigma}_{\nu}^{*}} \le & \sqrt{n}\mnorm{\L^{*} - \mtx{\Sigma}^{*}} + \sqrt{n}\nu\mnorm{\J - \frac{n - 1}{n}\I}\nonumber\\
\le &\sqrt{n}\mnorm{\L^{*} - \mtx{\Sigma}^{*}} + \nu\le n\bar{p}^{*} + \nu \lesssim n\bar{p}^{*}.\label{eq:L-Sigma}
\end{align}
 Furthermore, by \eqref{eq:binftyb2}, 
\begin{equation}
  \label{eq:b2}
  \sqrt{n}\bp_{2, \nu}(n^{-r})\lesssim \sqrt{n\bar{p}^{*}\log n}.
\end{equation}
Putting \eqref{eq:EU} - \eqref{eq:b2} together and using the fact that $\sqrt{n}\mnorm{\mtx{U}^{*}}\preceq 1$, we obtain that 
\begin{align}
 &\sqrt{n}\mnorm{\mtx{U}\sign(\mtx{U}^{T}\mtx{U}^{*})- \mtx{U}^{*}}\nonumber\\
  \lesssim &\frac{\sqrt{n\bar{p}^{*}\log n}}{n(\ubar{p}_{*} - p_{\es})} + \frac{n\bar{p}^{*}\sqrt{n\bar{p}^{*}\log n}}{(n(\ubar{p}_{*} - p_{\es}))^{2}} + \frac{n\bar{p}^{*}\sqrt{n\bar{p}^{*}\log n}}{(n(\ubar{p}_{*} - p_{\es}))^{2}}\frac{n\bar{p}^{*}}{n\bar{p}^{*}}\nonumber\\
\lesssim & \frac{\sqrt{n\bar{p}^{*}\log n}}{n(\ubar{p}_{*} - p_{\es})} + \frac{n\bar{p}^{*}\sqrt{n\bar{p}^{*}\log n}}{(n(\ubar{p}_{*} - p_{\es}))^{2}}\nonumber\\
\stackrel{(i)}{\lesssim} & \frac{(n\bar{p}^{*})^{3/4}(\log n)^{1/4}}{n(\ubar{p}_{*} - p_{\es})} + \frac{n\bar{p}^{*}\sqrt{n\bar{p}^{*}\log n}}{(n(\ubar{p}_{*} - p_{\es}))^{2}}\nonumber\\
\stackrel{(ii)}{\lesssim} & \frac{n\bar{p}^{*}\sqrt{n\bar{p}^{*}\log n}}{(n(\ubar{p}_{*} - p_{\es}))^{2}}\nonumber,
\end{align}
where (i) uses \eqref{eq:barp} and (ii) uses the condition \eqref{eq:cond_pbar_pubar}. As a consequence, there exists a constant $C$ that only depends on $r$ and $\xi$ such that
\[\sqrt{n}\mnorm{\mtx{U}\sign(\mtx{U}^{T}\mtx{U}^{*})- \mtx{U}^{*}}\le C\frac{n\bar{p}^{*}\sqrt{n\bar{p}^{*}\log n}}{(n(\ubar{p}_{*} - p_{\es}))^{2}}.\]
By the condition \eqref{eq:cond_pbar_pubar} again, 
\[\sqrt{n}\mnorm{\mtx{U}\sign(\mtx{U}^{T}\mtx{U}^{*})- \mtx{U}^{*}}\le \frac{C}{\sqrt{C_{\ell_{\infty}}}}.\]
If $C_{\ell_{\infty}}\ge C^2 / c^2$, 
\[\sqrt{n}\mnorm{\mtx{U}\sign(\mtx{U}^{T}\mtx{U}^{*})- \mtx{U}^{*}}\le c,\]
with probability $1 - (B(\td{j}) + 1)n^{-r}\ge 1 - (10K + 1)n^{-r}$.
% By \eqref{eq:A5}
% \[\frac{E_{+, \nu}^{2}(n^{-r})}{(\Delta_{\nu}^{*})^{2}} + \frac{\Theta_{\nu}(n^{-r})\sigmap_{\nu}(n^{-r})}{\Delta_{\nu}^{*}}\le c_{0}^{2} + c_{0} \le 2c_{0},\]
% where we use the fact that $c_{0}\le 1/4$.

\end{document}